\documentclass[11pt]{amsart}
\usepackage{amssymb}
\usepackage{amscd}
\usepackage{amsmath}
\usepackage{amsmath}
\usepackage[all]{xy}
\usepackage{mathrsfs}
\usepackage{setspace}
 \usepackage{color}

\theoremstyle{plain}
\newtheorem{theorem}{Theorem}[section]
\newtheorem{lemma}[theorem]{Lemma}
\newtheorem{corollary}[theorem]{Corollary}
\newtheorem{proposition}[theorem]{Proposition}

\theoremstyle{definition}

\newtheorem{hypothesis}[theorem]{Hypothesis}
\theoremstyle{definition}
\newtheorem{definition}[theorem]{Definition}
\newtheorem{remark}[theorem]{Remark}

\font\russ=wncyr10  1
\def\sha{\hbox{\russ\char88}}

\DeclareMathOperator{\Ext}{Ext}
\DeclareMathOperator{\Gal}{Gal}

\DeclareMathOperator{\Hom}{Hom}

\DeclareMathOperator{\Spec}{Spec}
\DeclareMathOperator{\N}{N}

\DeclareMathOperator{\im}{im}
\DeclareMathOperator{\coker}{coker}

\newcommand{\CC}{\mathbb{C}}

\newcommand{\GG}{\mathbb{G}}

\newcommand{\QQ}{\mathbb{Q}}

\newcommand{\RR}{\mathbb{R}}

\newcommand{\ZZ}{\mathbb{Z}}

\newcommand{\cA}{\mathcal{A}}

\newcommand{\cF}{\mathcal{F}}

\newcommand{\cG}{\mathcal{G}}
\newcommand{\cH}{\mathcal{H}}
\newcommand{\cI}{\mathcal{I}}
\newcommand{\cJ}{\mathcal{J}}
\newcommand{\cK}{\mathcal{K}}

\newcommand{\cN}{\mathcal{N}}
\newcommand{\co}{\mathcal{O}}
\newcommand{\cP}{\mathcal{P}}
\newcommand{\cQ}{\mathcal{Q}}
\newcommand{\cR}{\mathcal{R}}

\newcommand{\fa}{\mathfrak{a}}
\newcommand{\fb}{\mathfrak{b}}
\newcommand{\fc}{\mathfrak{c}}
\newcommand{\fd}{\mathfrak{d}}
\newcommand{\fm}{\mathfrak{m}}
\newcommand{\fn}{\mathfrak{n}}

\newcommand{\fq}{\mathfrak{q}}
\newcommand{\fr}{\mathfrak{r}}

\newcommand{\id}{\mathrm{id}}

\newcommand{\bz}{\mathbb{Z}}

\begin{document}

\title[]{On the theory of higher rank\\
Euler, Kolyvagin and Stark systems}

\author{David Burns and Takamichi Sano}

\begin{abstract}
Mazur and Rubin have recently developed a theory of higher rank Kolyvagin and Stark systems over principal artinian rings and discrete valuation rings. In this article we describe a natural extension of (a slightly modified version of) their theory to systems over more general coefficient rings. We also construct unconditionally, and for general $p$-adic representations, a canonical, and typically large, module of higher rank Euler systems and show that for $p$-adic representations satisfying standard hypotheses the image under a natural higher rank Kolyvagin-derivative type homomorphism of each such system is a higher rank Kolyvagin system that originates from a Stark system.
\end{abstract}

\address{King's College London,
Department of Mathematics,
London WC2R 2LS,
U.K.}
\email{david.burns@kcl.ac.uk}

\address{Osaka City University,
Department of Mathematics,
3-3-138 Sugimoto\\Sumiyoshi-ku\\Osaka\\558-8585,
Japan}
\email{sano@sci.osaka-cu.ac.jp}

\maketitle

\tableofcontents

\section{Introduction} \label{Intro}

\subsection{Some background}\label{background} Let $K$ be a number field, $p$ an {\it odd} prime number, $O$ the ring of integers of a finite extension $Q$ of $\QQ_p$ and $M$ a non-zero element of $O$. Let $T$ be a free $O$-module of finite rank that is endowed with a continuous $O$-linear action of the absolute Galois group $G_K$ of $K$ that is ramified at only finitely many places and set $\mathcal{A} := T/MT$.

An Euler system for $T$ is a collection of cohomology classes in the groups $H^1(F,T^\ast(1))$ for finite extensions $F$ of $K$ whose norm projections satisfy a specific family of relations (involving multiplication by Euler factors) as $F$ varies. This notion was introduced by Kolyvagin in the late 1980's and has played a vital role in the proof of many results in arithmetic geometry related to the structure of Selmer groups.

Mazur and Rubin subsequently showed that suitable `Kolyvagin systems' of families of cohomology classes in  $H^1(K,\mathcal{A}^\ast(1))$ controlled the structure of Selmer groups and that a certain modification of Kolyvagin's `derivative operator' gave a canonical homomorphism from the module of Euler systems of $T$ to the module of Kolyvagin systems of $\mathcal{A}$ (see \cite[Theorem 3.2.4]{MRkoly}).

In this way, it became clear that, under suitable hypotheses, Kolyvagin systems play the key role in obtaining structural results about Selmer groups and that the link to Euler systems is pivotal for the supply of Kolyvagin systems that are related to the special values of $L$-series.

For many representations $T$, however, families of cohomology classes are not sufficient to control the relevant Selmer groups and in such `higher rank' cases authors have considered various collections of elements in exterior powers of cohomology groups.

The theory of `higher rank Euler systems' has been well-understood for some time by now, with the first general approach being described by Perrin-Riou in \cite{PR} and significant contributions in an important special case made by Rubin \cite{rubinstark}.

More recently, in \cite{MRselmer} Mazur and Rubin have developed a theory of `higher rank Kolyvagin systems', and an associated notion of `higher rank Stark systems' (these are collections of cohomology classes generalizing the units predicted by Stark-type conjectures), and showed that, under suitable hypotheses, such systems can control Selmer groups.

However, the technical difficulties involved in working with higher exterior powers impacted on their theory in two important ways. Firstly, coefficient rings were restricted to being either principal artinian local rings or discrete valuation rings and, secondly, any link to the theory of higher rank Euler systems remained `mysterious'.

In addition, Mazur and Rubin explicitly point out that allowing more general coefficient rings in the theory is essential if one is to deal effectively with questions arising in deformation theory whilst one can expect that a concrete link between (higher rank) Euler systems on the one hand and higher rank Kolyvagin and Stark systems on the other would have a key role to play in any subsequent arithmetic applications. The present article was motivated by these problems.

Speaking broadly, we shall introduce algebraic techniques that allow a resolution of the first problem. In addition, we shall define a canonical (equivariant) higher Kolyvagin derivative homomorphism and use an analysis of Galois cohomology complexes to construct unconditionally a large family of Euler systems (of arbitrary rank) for which one can provide a canonical link to the corresponding equivariant theories of higher rank Kolyvagin and Stark systems.

In the course of obtaining these results we shall also, amongst other things, answer an explicit question of Kato (see Remark \ref{kato question}),  show that a large family of Euler systems validate a natural extension of congruence conjectures of Darmon, of Mazur and Rubin and of the second-named author (see Remark \ref{gen darmon}) and provide a new interpretation of the `finite-singular comparison map' of Mazur and Rubin in terms of natural Bockstein homomorphisms (see Remark \ref{finite-sing rem}).

\subsection{Modified Kolyvagin and Stark systems} 

The key idea in overcoming the algebraic difficulties faced by Mazur and Rubin in \cite{MRselmer} is to adopt a slight change of approach.

To do this we define for any commutative ring $R$, any $R$-module $X$ and any non-negative integer $r$ the `$r$-th exterior bidual' of $X$ to be the $R$-module obtained by setting

\begin{equation*}\label{bidual def} {\bigcap}_R^r X:= \Hom_R\left({\bigwedge}_R^r \Hom_R(X,R),R\right),\end{equation*}
and we note that there is a natural homomorphism of $R$-modules
\begin{equation}\label{can homo} {\bigwedge}_R^r X \to {\bigcap}_R^r X \end{equation}
that sends each element $a$ to the homomorphism $\Phi \mapsto \Phi(a)$.

Whilst exterior biduals already played a role in our earlier work with Kurihara \cite{bks2-2} we only now develop their basic theory in an appendix and find that, in several key ways, it is better behaved than the corresponding theory of exterior powers (this is possible since the maps (\ref{can homo}) are in general neither injective nor surjective). In fact, since Gorenstein rings are particularly important in arithmetic applications (arising, for example, as either integral group rings or Hecke algebras), we often assume that coefficient rings are of this type when developing the theory.

In particular, for any finite self-injective ring (i.e. zero-dimensional Gorenstein ring) $R$ and any finite $G_K$-representation $\cA$ that is free over $R$ (and satisfies suitable hypotheses) the theory can be combined with  the approach of \cite{MRselmer} to define canonical modules of $R$-equivariant rank $r$ Stark systems ${\rm SS}_r(\mathcal{A})$ and $R$-equivariant rank $r$ Kolyvagin systems ${\rm KS}_r(\cA)$ for $\cA$ together with a canonical homomorphism of $R$-modules
\[ {\rm Reg}_r : {\rm SS}_r(\mathcal{A}) \to {\rm KS}_r(\mathcal{A}). \]
(See Definitions \ref{def unit}, \ref{def koly}, and \ref{def stub}.) This homomorphism can be naturally regarded as a `regulator' and we denote its image by ${\rm KS}_r^{{\rm reg}}(\mathcal{A})$.

If $R$ is any ring as in \cite{MRselmer}, then any Stark or Kolyvagin system constructed in loc. cit. gives rise (via the canonical homomorphisms (\ref{can homo})) to a Stark or Kolyvagin system in our generalized sense. In this way, the module ${\rm KS}_r^{{\rm reg}}(\mathcal{A})$ of `regulator Kolyvagin systems' corresponds to the module ${\rm KS}_r'(\mathcal{A})$ of `stub Kolyvagin systems' considered in loc. cit. and is in general strictly smaller than ${\rm KS}_r(\mathcal{A})$.

Under certain additional restrictions on $\cA$ and $R$ we can also prove equivariant refinements of the link between Stark systems and the structure of Selmer modules that are established by Mazur and Rubin in \cite[Proposition 8.5]{MRselmer}. In this way, for example, we obtain results such as the following (for details of which see Theorem \ref{theorem higher fitt}).

\begin{theorem}\label{equiv MR} Under certain natural hypotheses on $R$ and $\mathcal{A}$ the $R$-module ${\rm SS}_r(\mathcal{A})$ is free of rank one and any basis of this module explicitly determines all (higher) Fitting ideals of the $R$-module
$\sha^2(\co_{K,S},\cA^\ast(1))$. \end{theorem}

In a similar way we will also later prove natural equivariant refinements of both \cite[Theorem 4.5.6]{MRkoly} and \cite[Corollary 13.1]{MRselmer} (see Theorem \ref{bounding selmer}).

\subsection{Derivable Euler systems}

To discuss the link between higher rank Euler systems and higher rank Stark and Kolyvagin systems we fix a representation $T$ of $G_K$ as in \S\ref{background}. We assume that $T$ is a free module with respect to an action of a commutative Gorenstein $O$-order $\mathcal{R}$ (that commutes with the given action of $G_K$). We note that in this case for any non-zero element $M$ of $O$ the ring $R := \mathcal{R}/M\mathcal{R}$ is self-injective and the module $\mathcal{A} := T/MT$ is free over $R$.

We also fix an abelian pro-$p$-extension $\cK$ of $K$ (that one should typically think of as `large' - see, for example, Hypothesis \ref{hyp last}).

In this setting we define a module ${\rm ES}_{r}(T,\cK)$ of $\mathcal{R}$-equivariant Euler systems of rank $r$ for the pair $(T,\mathcal{K})$ that has a natural structure as $\mathcal{R}[[\Gal(\cK/K)]]$-module (see Definition \ref{def euler}).

\subsubsection{}To discuss `higher Kolyvagin derivatives' in this setting we first define a canonical $R$-module ${\rm KC}_r(\mathcal{A})$ of `Kolyvagin collections of rank $r$' for $\mathcal{A}$ that contains ${\rm KS}_r(\mathcal{A})$ as a canonical $R$-submodule (see Definition \ref{def KC}).

For large enough fields $\cK$ we can then define a natural `higher rank Kolyvagin derivative homomorphism'
\[ \mathcal{D}_r : {\rm ES}_{r}(T,\cK) \to {\rm KC}_r(\mathcal{A}) \]
and we say that an Euler system in ${\rm ES}_{r}(T,\cK)$ is `Kolyvagin-derivable', respectively `Stark-derivable', for the $R$-module $\mathcal{A}$, if its image under $\mathcal{D}_r$ belongs to ${\rm KS}_r(\mathcal{A})$, respectively to ${\rm KS}_r^{{\rm reg}}(\mathcal{A})$ (see Definition \ref{derivable}).

In this context Theorem \ref{equiv MR} leads to results such as the following.

\begin{theorem}[{Theorem \ref{stark der thm 2}}]\label{stark der thm}
Fix a system $c$ in ${\rm ES}_r(T,\mathcal{K})$ that is Stark-derivable for the $R$-module $\mathcal{A}$.
Then, under certain natural hypotheses on $T$, there exists a canonical (non-zero) ideal $J=J_{c,T,M}$ of $R$ for which
 $ \langle \mathcal{D}_r( c) \rangle_{R}= J\cdot {\rm KS}_r^{{\rm reg}}(\mathcal{A})$ and the associated module of Stark systems
\[ \{ \epsilon \in {\rm SS}_r(\mathcal{A}) \mid {\rm Reg}_r(\epsilon) \in  \langle \mathcal{D}_r( c) \rangle_{R}\}\]
explicitly determines $J\cdot{\rm Fitt}_{R}^i(\sha^2(\co_{K,S},\cA^\ast(1)))$ for every non-negative integer $i$.
\end{theorem}


%
%

\subsubsection{}If $r=1$ and $\mathcal{R} = O$ (so $R = O/MO$) then the map $\mathcal{D}_r$ defined here coincides with the homomorphism ${\rm ES}_1(T,\cK) \to {\rm KS}_1(\mathcal{A})$ considered by Mazur and Rubin in \cite{MRkoly}, one has ${\rm KS}_1(\mathcal{A}) = {\rm KS}_1^{{\rm reg}}(\mathcal{A})$ (by \cite[Theorem 4.4.1]{MRkoly}) and every element of ${\rm ES}_1(T,\cK)$ is Kolyvagin (and hence Stark)-derivable for the $R$-module $\mathcal{A}$ (by \cite[Theorem 3.2.4]{MRkoly}).

More generally, if $r=1$, then for any $\mathcal{R}$ as above the approach of \cite{MRkoly} can be used to show that every element of ${\rm ES}_1(T,\mathcal{K})$ is Kolyvagin-derivable for the $\mathcal{R}$-module $\mathcal{A}$.

However, if either $r > 1$ or $R$ is not a principal ideal ring, then ${\rm KS}_r^{{\rm reg}}(\mathcal{A})$ is often strictly smaller than ${\rm KS}_r(\mathcal{A})$ and so Theorem \ref{stark der thm} leads to the problem of deciding which Euler systems are Stark-derivable. (Note also that if $R$ is a principal ideal ring, then for $r > 1$ our theory differs slightly from that developed by Mazur and Rubin in that we consider elements in exterior power biduals rather than exterior powers).

Whilst we do not present a complete solution to this problem, the following result shows (unconditionally) that Stark-derivable Euler systems are common in very general contexts and also that, modulo certain standard conjectures, all Euler systems that arise from the special values of motivic $L$-functions should be Stark-derivable.

%


In this result we write $S_\infty$ for the set of archimedean places of $K$. 

\begin{theorem}\label{basic euler system thm}
 Assume that $H^0(K,T^*(1))$ vanishes. Assume also that the $\cR$-module $\bigoplus_{v \in S_\infty}H^0(K_v, T)$ is free and write $r$ for its rank.
 
Then there exists a canonical (non-zero) cyclic $\mathcal{R}[[\Gal(\cK/K)]]$-submodule ${\mathcal{E}}^{\rm b}( T,\mathcal{K})$ of ${\rm ES}_r( T,\mathcal{K})$ that is in most cases free of rank one. Further, if, in addition,  $\cK$ is large enough and $T$ satisfies a variety of standard hypotheses, then every Euler system in $\mathcal{E}^{\rm b}( T,\cK)$ is Stark-derivable for $\mathcal{A}$, i.e. one has $\mathcal{D}_r({\mathcal{E}}^{\rm b}( T,\mathcal{K})) \subseteq  {\rm KS}_r^{{\rm reg}}(\mathcal{A})$, and there is an Euler system $c$ in ${\mathcal{E}}^{\rm b}( T,\mathcal{K})$ for which one has $J_{c, T,M}=R$ in terms of the notation in Theorem \ref{stark der thm}.
\end{theorem}

The first sentence in the second paragraph of Theorem \ref{basic euler system thm} follows directly from Theorem \ref{det euler fitt} (and Remark \ref{remark free}) below and the second sentence from Theorem \ref{main}. We note also that the rank $r$ of $\bigoplus_{v \in S_\infty}H^0(K_v, T)$ that occurs in Theorem \ref{basic euler system thm} corresponds to the notion of `core rank' introduced by Mazur and Rubin (see Remark \ref{remark core rank}).

In the sequel we will call an Euler system `basic' (for $T$ and $\cK$) if it lies in $\mathcal{E}^{\rm b}(T,\cK)$.


\begin{remark}\label{kato question} We construct the module ${\mathcal{E}}^{\rm b}( T,\mathcal{K})$ by means of a detailed analysis of complexes in Galois cohomology (see Theorem \ref{theorem vertical} and Definition \ref{def basic}). This construction is motivated (broadly speaking) by an approach of Kato to the study of arithmetic zeta elements and the formulation of generalized main conjectures in Iwasawa theory. In particular, Remark \ref{how big rem}(ii) shows that Theorem \ref{basic euler system thm} provides an effective answer to the question explicitly raised by Kato in \cite[Remark 4.14]{kato} as to whether the classical techniques of Kolyvagin could be applied in the general setting of arithmetic zeta elements. In a subsequent article we will explain how this approach extends to the context of non-abelian extensions of number fields and also how, in certain circumstances, it leads to the construction of new Euler systems with prescribed properties.
\end{remark}

\begin{remark}\label{gen darmon} As a key part of the proof of Theorem \ref{basic euler system thm} we show that every basic Euler system satisfies a natural extension to general $p$-adic representations of a conjecture of Darmon \cite[Conjecture 4.3]{D} concerning congruences for cyclotomic units. This result also constitutes a natural extension to general $p$-adic representations of the explicit congruence conjectures on Rubin-Stark elements proposed by Mazur-Rubin in \cite[Conjecture 5.2]{MRGm} and by the second author in \cite[Conjecture 3]{sano}. See Remark \ref{remark mrs} for more details.
\end{remark}

\begin{remark}\label{how big rem} It is natural to ask how common it is for an Euler system to be basic. We will address this problem elsewhere but, for the moment, record the following facts showing that ${\mathcal{E}}^{\rm b}( T,\mathcal{K})$ is a large, and arithmetically significant, submodule of ${\rm ES}_r( T,\mathcal{K})$. 

\noindent{}(i) In many cases $\mathcal{E}^{\rm b}(T,\cK)$ is free of rank one over $\cR[[\Gal(\cK/K)]]$ (see Remark \ref{remark free}).

\noindent{}(ii) Standard conjectures imply that all Euler systems arising from the leading terms of motivic $L$-functions are basic (see Remark \ref{motivating exam}).

\noindent{}(iii) Suppose $\cR=O$, fix a $\ZZ_p$-power extension $K_\infty$ of $K$ in $\cK$, set $\Lambda := O[[\Gal(K_\infty/K)]]$ and write $Q(\Lambda)$ for the quotient field of $\Lambda$. Then one can show that the $\Lambda$-module ${\varprojlim}_{F\in \Omega(K_\infty/K)}H^2(\mathcal{O}_{F,S},T^*(1))$ is torsion if and only if the projection ${\mathcal{E}}^{\rm b}(T,\mathcal{K})_\infty$ of ${\mathcal{E}}^{\rm b}(T,\mathcal{K})$ to $ Q(\Lambda)\otimes_\Lambda ({\bigcap}_\Lambda^rH^1(\mathcal{O}_{K_\infty,S},T^*(1)))$ is non-zero in which case both ${\bigcap}_\Lambda^rH^1(\mathcal{O}_{K_\infty,S},T^*(1))$ and ${\mathcal{E}}^{\rm b}(T,\mathcal{K})_\infty$ span the same dimension one space over $Q(\Lambda)$.

\noindent{}(iv) The proof of Theorem \ref{basic euler system thm} will make clear that the key inclusion $\mathcal{D}_r({\mathcal{E}}^{\rm b}( T,\mathcal{K})) \subseteq  {\rm KS}_r^{{\rm reg}}(\mathcal{A})$ remains valid if one replaces ${\mathcal{E}}^{\rm b}( T,\mathcal{K})$ by the (a priori)  larger submodule ${\mathcal{E}}^{\rm b}( T,\mathcal{K})'$ of ${\rm ES}_r(T,\mathcal{K})$ comprising all Euler systems that agree at the $K$-th level with a system in ${\mathcal{E}}^{\rm b}( T,\mathcal{K})$.  If $r=1$, $K =\QQ$ and $\cR=O$, then Kato's argument in \cite{katoselmer} can be used to show that, under certain standard hypotheses, 
the $O[[\Gal(\cK/K)]]$-module
${\rm ES}_r( T,\mathcal{K})/{\mathcal{E}}^{\rm b}( T,\mathcal{K})'$ is both elementary in nature and `small' in size, being annihilated by a product of local Euler factors.

\end{remark}

\begin{remark}\label{finite-sing rem} The `finite-singular comparison map' defined by Mazur and Rubin plays a key role in their theory of (higher rank) Kolyvagin systems \cite{MRkoly, MRselmer}. As an
important step in the proof of Theorem \ref{basic euler system thm} we give a natural interpretation of this map in terms of Bockstein homomorphisms arising from Galois descent. For details see \S\ref{section bock}.
\end{remark}

\begin{remark} After the authors finished this work, the second author was informed by Ryotaro Sakamoto that he had also independently developed a theory of equivariant Stark systems by a similar method to ours, using exterior biduals over zero dimensional Gorenstein rings, and in this way obtained results similar to Theorem \ref{structure unit 2} and Theorem \ref{theorem higher fitt}(ii) (see \cite{sakamoto}).
\end{remark}

\subsection{Some general notation}
We fix an algebraic closure $\overline \QQ$ of $\QQ$. A number field is a subfield of $\overline \QQ$ which is finite over $\QQ$. For any number field $F$, we denote $G_F:=\Gal(\overline \QQ/F)$. Fix a number field $K$. The ring of integers of $K$ is denoted by $\co_K$. For each place $v$ of $K$, we fix a place $w$ of $\overline \QQ$ lying above $v$. The decomposition group of $w$ in $G_K$ is identified with $\Gal(\overline \QQ_w/ K_v) $. For any finite extension $F/K$, the completion of $F$ at the place lying under $w$ is denoted by $F_v$, and we set $G_{F_v}:=\Gal(\overline \QQ_w/F_v)$. The maximal unramified extension of $F_v$ is denoted by $F_v^{\rm ur}$, and the inertia subgroup of $G_{F_v}$ is defined by $I_{F_v}:=\Gal(\overline \QQ_w/F_v^{\rm ur})$. The set of all archimedean places of $F$ is denoted by $S_\infty(F)$ or simply by $S_\infty$ if there is no danger of confusion. For a prime number $p$, we denote the set of places of $F$ lying above $p$ by $S_p(F)$. The set of places of $K$ which ramify in $F$ is denoted by $S_{\rm ram}(F/K)$.

A finite place of $K$ is regarded as a prime ideal of $\co_K$, so we often call it a prime (of $K$). A prime of $K$ is often denoted by $\fq$. The Frobenius automorphism at $\fq$ is denoted by ${\rm Fr}_\fq \in \Gal(K_\fq^{\rm ur}/K_\fq)$.

For any commutative ring $R$, we write $D(R)$ for the derived category of complexes of $R$-modules and $D^{\rm p}(R)$ for the full triangulated subcategory of $D(R)$ comprising perfect complexes (that is complexes that are isomorphic in $D(R)$ to a bounded complex of finitely generated projective $R$-modules).

Let $N$ be a continuous $G_K$-module. For a place $v$ of $K$, we say that $N$ is ramified (resp. unramified) at $v$ if $I_{K_v}$ acts non-trivially (resp. trivially) on $N$. The set of places of $K$ at which $N$ is ramified is denoted by $S_{\rm ram}(N)$. Suppose that $S_{\rm ram}(N)$ is finite, and take a finite set $U$ of places of $K$ such that
$$S_{\rm min}(N):=S_\infty(K) \cup S_p(K) \cup S_{\rm ram}(N) \subset U.$$
Then $N$ is regarded as $\Gal(K_U/K)$-module, where $K_U/K$ is the maximal Galois extension unramified outside $U$. In general, for a profinite group $\mathcal{G}$ and a continuous $\mathcal{G}$-module $M$, we denote the standard complex of continuous cochains of $M$ by $R\Gamma(\mathcal{G},M)$. If $\mathcal{G}$ is the Galois group of some field extension $E/F$, then we denote $R\Gamma(\mathcal{G},M)$ by $R\Gamma(E/F, M)$, and by $R\Gamma(F,M)$ when $E$ is the separable closure of $F$. We also denote
$$R\Gamma(\co_{K,U},N):=R\Gamma(K_U/K,N).$$
For every place $v$ of $K$, there is a natural localization morphism
$$R\Gamma(\co_{K,U},N) \to R\Gamma(K_v,N)$$
by inflation and restriction with respect to $G_K \to \Gal(K_U/K)$ and $G_{K_v} \subset G_K$ respectively.
The compactly supported cohomology is defined by
$$R\Gamma_c(\co_{K,U}, N):={\rm Cone}\left(R\Gamma(\co_{K,U},N) \to \bigoplus_{v \in U}R\Gamma(K_v,N) \right)[-1].$$
For a place $v \notin S_{\rm min}(N)$, we define
$$R\Gamma_f(K_v, N):=R\Gamma(K_v^{\rm ur}/K_v,N).$$
We use the fact that this complex is represented by
\begin{equation}\label{unram resolution} N \stackrel{1- {\rm Fr}_v^{-1}}{\to} N,\end{equation}
where the first term is placed in degree zero and ${\rm Fr}_v$ denotes the Frobenius automorphism in $\Gal(K_v^{\rm ur}/K_v)$.

There is also a natural `inflation' morphism
$$R\Gamma_f(K_v,N) \to R\Gamma(K_v,N)$$
and we define $R\Gamma_{/f}(K_v,N)$ to be the mapping cone of this morphism.

The cohomologies $H^i(R\Gamma_\ast(\cdot,N))$, where $\ast \in \{\emptyset, c, f, /f \}$, are denoted by $H^i_\ast(\cdot,N)$. For a finite extension $F/K$, the set of places of $F$ lying above places in $U$ is denoted by $U_F$. We denote $R\Gamma_\ast(\co_{F,U_F}, N)$ simply by $R\Gamma_\ast(\co_{F,U},N)$.

For any positive integer $m$, we denote the group of $m$-th roots of unity in $\overline \QQ^\times$ by $\mu_m$. We define the $G_K$-module $\ZZ_p(1)$ by the inverse limit $\varprojlim_n \mu_{p^n}$.

For a commutative ring $R$ and an $R$-module $X$, the dual module $\Hom_R(X,R)$ is often denoted by $X^\ast$. When $X$ is free and $\{b_1,\ldots,b_d\}$ is a basis, we denote the dual of $b_i$, namely, the map
$$X\to R; \ b_j \mapsto \begin{cases}
1 &\text{ if }i=j,\\
0 &\text{ if }i\neq j,
\end{cases}
$$
by $b_i^\ast \in X^\ast$ for each $i$ with $1\leq i \leq d$.


\section{Higher rank Euler systems}\label{section 2}

In this section, we study higher rank Euler systems. The notion of higher rank Euler systems first appeared in Perrin-Riou's work \cite{PR}. She defined a rank $r$ Euler system as a collection of elements in $r$-th exterior powers of Galois cohomologies of degree one, but as remarked by Mazur and Rubin in \cite[Remark 8.7]{MRGm}, we should use the exterior bidual (see Definition \ref{def exterior bidual} and Remark \ref{def ext bidual rem} below) rather than the usual exterior power. One reason for this is that the exterior bidual defines lattices in which conjectural Rubin-Stark elements should exist (see \cite[Conjecture B$'$]{rubinstark}). Rubin-Stark elements give an example of higher rank Euler systems. It is widely believed that for a general motive there is an Euler system which is related to the value of the motivic $L$-function. In the case of untwisted Tate motives (often referred to as the `$\GG_m$ case'), Kurihara and the present authors observed in \cite[Theorem 5.14]{bks1} that conjectural `zeta elements', whose existence is predicted by the equivariant Tamagawa number conjecture (eTNC, \cite[Conjecture 4(iv)]{BFetnc}), naturally provide Rubin-Stark elements as the image under certain `projectors' from determinants of cohomology. In this paper, we develop their theory for general $p$-adic representations, and show that a certain system in determinants of cohomology, which we call a `vertical determinantal system' (see Definition \ref{definition vertical}),  provides a higher rank Euler system (see Theorem \ref{theorem vertical}). This gives further evidence for the idea that higher rank Euler systems should exist in exterior power biduals rather than in exterior powers themselves.


\subsection{Exterior power biduals}

Let $R$ be a commutative ring, and $X$ an $R$-module. For a positive integer $r$ and a homomorphism $\varphi \in X^\ast$, define
$$\varphi^{(r)}: {\bigwedge}_R^r X \to {\bigwedge}_R^{r-1} X$$
by
$$\varphi^{(r)}(x_1\wedge\cdots\wedge x_r) := \sum_{i=1}^r (-1)^{i+1} \varphi(x_i) x_1\wedge\cdots\wedge x_{i-1}\wedge x_{i+1} \wedge \cdots x_r.$$
By abuse of notation, we often simply denote $\varphi^{(r)}$ by $\varphi$.

For non-negative integers $r$ and $s$ with $r \leq s$, this construction induces a natural homomorphism
$${\bigwedge}_R^r (X^\ast) \to \Hom_R\left({\bigwedge}_R^s X, {\bigwedge}_R^{s-r} X\right); \ \varphi_1\wedge \cdots \wedge \varphi_r \mapsto \varphi_r \circ \cdots \circ \varphi_1$$
that can be described more explicitly by
\begin{equation}\label{exterior explicit}
\varphi_1\wedge \cdots \wedge \varphi_r \mapsto \left(x_1\wedge \cdots \wedge x_s \mapsto \sum_{\sigma \in \mathfrak{S}_{s,r}} {\rm sgn}(\sigma) \det(\varphi_i(x_{\sigma(j)}))_{1\leq i,j \leq r} x_{\sigma(r+1)}\wedge \cdots \wedge x_{\sigma(s)} \right),\end{equation}
with $\mathfrak{S}_{s,r}:=\{ \sigma \in \mathfrak{S}_{s} \mid \sigma(1)<\cdots < \sigma (r) \text{ and } \sigma(r+1)<\cdots <\sigma(s)\}.$

We use this map to regard any element of ${\bigwedge}_R^r (X^\ast)$ as an element of $\Hom_R({\bigwedge}_R^s X, {\bigwedge}_R^{s-r} X)$. In particular, for any $a\in {\bigwedge}_R^s X$ and $\Phi \in {\bigwedge}_R^r (X^\ast)$, the element $\Phi(a) \in {\bigwedge}_R^{s-r} X$ is defined.


In the sequel we make much use of the following $R$-modules.

\begin{definition}\label{def exterior bidual}
For any non-negative integer $r$, we define the `$r$-th exterior bidual' of $X$ to be the $R$-module obtained by setting
$${\bigcap}_R^r X:=\left({\bigwedge}_R^r (X^\ast) \right)^\ast.$$
\end{definition}

\begin{remark}\label{def ext bidual rem}The basic theory of these modules is developed in an appendix. At this stage we just note that there is a canonical homomorphism (\ref{can homo}) and that if $O$ is a Dedekind domain with quotient field $Q$ and $R$ is an $O$-order in some finite dimensional semisimple (commutative) $Q$-algebra, then this homomorphism induces an identification of ${\bigcap}_R^r X$ with the lattices `$\wedge_0^r M$' that are introduced by Rubin in \cite[\S 1.2]{rubinstark} (for details see Proposition \ref{rubin prop 2}).\end{remark}

\subsection{The definition of higher rank Euler systems}\label{def hres}

Let $K$ be a number field. Let $p$ be an {\it odd} prime number. Let $Q/\QQ_p$ be a finite extension and $O$ be the ring of integers in $Q$. Let $\cQ$ be a finite dimensional semisimple commutative $Q$-algebra, and $\cR$ be a Gorenstein $O$-order in $\cQ$ (see \S \ref{A order}).
Let $T$ be a $p$-adic representation of $G_K$ with coefficients in $\cR$, namely,  a free $\cR$-module of finite rank with a continuous $\cR$-linear action of $G_K$.
We assume that $S_{\rm ram}(T)$ is finite.
Let $S$ be a finite set of places of $K$ such that
$$S_{\rm min}(T)=S_\infty(K)\cup S_p(K) \cup S_{\rm ram}(T)\subset S. $$
We set $T^\ast(1):=\Hom_\cR(T, \cR(1))$, where $\cR(1):=\cR \otimes_{\ZZ_p} \ZZ_p(1)$. 





Let $r$ be a non-negative integer. Let $\cK/K$ be a (possibly infinite) abelian $p$-extension such that every $v \in S_\infty(K)$ splits completely in $\cK$.
Let $\Omega(\cK/K)$ be the set of intermediate fields of $\cK/K$ that are of finite degree over $K$. For $F\in \Omega(\cK/K)$, set
$$\cG_F:=\Gal(F/K)$$
and
$$S(F):=S \cup S_{\rm ram}(F/K).$$
Also, we denote the induced module ${\rm Ind}_{G_K}^{G_F} (T)$ by $T_F$. Note that $T_F$ is identified with $\cR[\cG_F] \otimes_\cR T$, on which $G_K$ acts by
$$\sigma \cdot (a \otimes t) := a \overline \sigma^{-1} \otimes \sigma t \text{ ($\sigma \in G_K$, $a \in \cR [\cG_F]$, and $t \in T$)},$$
where $\overline \sigma \in \cG_F$ is the image of $\sigma$.

For a prime $\fq \notin S$, define a polynomial
$$P_\fq(x):=\det(1-{\rm Fr}_\fq^{-1}x \mid T) \in \cR[x].$$

In the following, we often use the following fact: for a commutative ring $R$, a finite abelian group $G$, and an $R[G]$-module $X$, there is a  natural isomorphism of $R[G]$-modules
$$\Hom_R(X,R) \stackrel{\sim}{\to}\Hom_{R[G]}(X,R[G])^\#; \ f \mapsto \sum_{\sigma \in G}f(\sigma(\cdot))\sigma^{-1},$$
where $G$ acts on $\Hom_R(X,R)$ by
$$(\sigma \cdot f) (a):=f(\sigma^{-1}  a) \text{ ($\sigma \in G$, $f \in \Hom_R(X,R)$, $a \in X$)},$$
and for any $G$-module $Y$ we denote by $Y^\#$ the module on which $G$ acts via the involution $\sigma \mapsto \sigma^{-1}$.
From this, we often regard
$$X^\ast:=\Hom_R(X,R)^\# =\Hom_{R[G]}(X,R[G]).$$

Note that, for $F, F'\in \Omega(\cK/K)$ with $F\subset F'$, the corestriction map
$${\rm Cor}_{F'/F}: H^1(\co_{F',S(F')},T^\ast(1)) \to H^1(\co_{F,S(F')}, T^\ast(1))$$
induces
$${\bigwedge}_{\cR[\cG_F]}^r (H^1(\co_{F,S(F')}, T^\ast(1))^\ast) \to {\bigwedge}_{\cR[\cG_{F'}]}^r (H^1(\co_{F',S(F')},T^\ast(1))^\ast).$$
Taking the dual of this map, we have the map
$${\bigcap}_{\cR[\cG_{F'}]}^r H^1(\co_{F',S(F')},T^\ast(1)) \to {\bigcap}_{\cR[\cG_F]}^r H^1(\co_{F,S(F')}, T^\ast(1)).$$
We denote this map also by ${\rm Cor}_{F'/F}$.

Also, note that we have an exact sequence
$$0 \to H^1(\co_{F,S(F)}, T^\ast(1)) \to H^1(\co_{F,S(F')}, T^\ast(1)) \to \bigoplus_{w \in (S(F')\setminus S(F))_F}H^1_{/f}(F_w, T^\ast(1)).$$
Since each $H^1_{/f}(F_w, T^\ast(1))$ is $O$-free (see \cite[Lemma 1.3.5(ii)]{R}), we can regard
$${\bigcap}_{\cR[\cG_F]}^r H^1(\co_{F,S(F)}, T^\ast(1)) \subset {\bigcap}_{\cR[\cG_F]}^r H^1(\co_{F,S(F')}, T^\ast(1)),$$
 by Proposition \ref{prop injective}(i) and (\ref{extvanish}) (note that $\cR[\cG_F]$ is also Gorenstein).

\begin{definition} \label{def euler}
An {\it Euler system of rank $r$} for $(T, \cK)$ is a collection
$$\left\{ c_F \in  {\bigcap}_{\cR[\cG_F]}^r H^1 (\co_{F,S(F)},T^\ast(1)) \ \middle| \ F\in \Omega(\cK/K)\right\}$$
that satisfies
$${\rm Cor}_{F'/F}(c_{F'})=\left(\prod_{\fq \in S(F')\setminus S(F)} P_\fq({\rm Fr}_\fq^{-1})\right)c_F \text{ in }{\bigcap}_{\cR[\cG_F]}^r H^1(\co_{F,S(F')}, T^\ast(1)),$$
for every $F$ and $F'$ in $\Omega(\cK/K)$ with $F\subset F'$.

The set of Euler systems of rank $r$ for $(T,\cK)$ has a natural structure of $\cR[[\Gal(\cK/K)]]$-module, which we denote by ${\rm ES}_r(T,\cK)$.
\end{definition}

\begin{remark}
When $r=1$ and $H^1(\co_{F,S(F)},T^\ast(1))$ is $O$-free, we have
$${\bigcap}_{\cR[\cG_F]}^r H^1 (\co_{F,S(F)},T^\ast(1))=H^1(\co_{F,S(F)},T^\ast(1))^{\ast \ast}=H^1(\co_{F,S(F)},T^\ast(1))$$
(see (\ref{2dual})).
So an Euler system of rank one is regarded as a collection
$$\{ c_F \in  H^1 (\co_{F,S(F)},T^\ast(1)) \mid F\in \Omega(\cK/K)\}.$$
Note that our Euler systems for $T$ correspond to those for $T^\ast(1)$ in the sense of Rubin \cite[Definition 2.1.1]{R}.
\end{remark}

\begin{remark} \label{remark rubin stark} An archetypal source of higher rank Euler systems is the (conjectural) family of Rubin-Stark elements. More precisely, if $K$ is a totally real field, then Rubin-Stark elements associated to the extensions $F/K$ for $F$ in $\Omega(\cK/K)$ should constitute an Euler system of rank $[K:\QQ]$ for the representation $T=\ZZ_p$. See \cite[\S 6]{rubinstark}.
\end{remark}

\subsection{$\Sigma$-modification}\label{section sigma}

In the sequel it is necessary to slightly modify the complexes that occur naturally in arithmetic in order to ensure that they have torsion-free cohomology in the lowest degree.

To explain this process we fix a finite set $\Sigma$ of places of $K$ with the property that $\Sigma$ is disjoint from $S$ and that $\cK/K$ is unramified at $\Sigma$.

\subsubsection{}Let $N$ be an $\mathcal{R}$-module endowed with a continuous action of $\Gal(K_{S}/K)$.

Let $F \in \Omega(\cK/K)$. If $U$ is a finite set of places of $K$ which contains $S(F)$ and is disjoint from $\Sigma$, then for any place $w$ in $\Sigma_F$ the natural localisation morphism $R\Gamma(\co_{F,U},N) \to R\Gamma(F_w,N)$ factors through the inflation morphism $R\Gamma_f(F_w,N) \to R\Gamma(F_w,N)$ and hence induces a canonical morphism $\iota_{w,N}: R\Gamma(\co_{F,U},N) \to R\Gamma_f(F_w,N)$.

We define the `$\Sigma$-modified cohomology complex' $R\Gamma_{\Sigma}(\co_{F,U},N)$ to be a complex that lies in an exact triangle in the derived category $D(\mathcal{R}[\cG_F])$ of complexes of $\mathcal{R}[\cG_F]$-modules 
\begin{equation}\label{Sigma-mod def} R\Gamma_{\Sigma}(\co_{F,U},N) \to R\Gamma(\co_{F,U},N) \xrightarrow{(\iota_{w,N})_w}
\bigoplus_{w\in \Sigma_F}R\Gamma_f(F_w,N) \to \end{equation}
and then in each degree $i$ set $H^i_\Sigma(\co_{F,U},N) := H^i(R\Gamma_{\Sigma}(\co_{F,U},N))$.

We note, in particular, the long exact cohomology sequence of this triangle gives an exact sequence
\begin{multline}\label{lesc} 0 \to H^0_\Sigma(\co_{F,U},N) \to H^0(F, N)\to \bigoplus_{w\in \Sigma_F}H^0(F_w,N) \\\to H^1_{\Sigma}(\co_{F,U},N) \to H^1(\co_{F,U},N) \to \bigoplus_{w\in \Sigma_F}H^1_f(F_w,N) \\ \to
H^2_{\Sigma}(\co_{F,U},N) \to H^2(\co_{F,U},N) \to 0.\end{multline}

\begin{remark}\label{injectivity} If $N$ is $O$-free, then the exact sequence (\ref{lesc}) induces an identification of $O$-torsion subgroups
\[ H^1_{\Sigma}(\co_{F,U},N)_{\rm tor} = \ker \left(H^1(\co_{F,U},N)_{\rm tor} \to \bigoplus_{w\in \Sigma_F}H^1_f(F_w,N)_{\rm tor} \right).\]
Here $(\cdot)_{\rm tor}$ denotes the $O$-torsion submodule. In many cases it is easy to choose $\Sigma$ so that this group vanishes. For example, if $O=\cR= \ZZ_p$ and $N = \ZZ_p(1)$ then
 $H_\Sigma^1(\co_{F,U},\ZZ_p(1))$ identifies with the $p$-completion of the $(U,\Sigma)$-unit group $\ker(\co_{F,U}^\times \to \bigoplus_{w\in \Sigma_F}(\co_F/w)^\times)$ of $F$ and is torsion-free for any non-empty set $\Sigma$. \end{remark}

\subsubsection{}\label{Sigma mod O}We define a `$\Sigma$-modification' of the compactly supported cohomology complex $R\Gamma_{c}(\co_{F,U},T)$ of $T$ by setting
\begin{eqnarray}\label{AV}
R\Gamma_{c,\Sigma}(\co_{F,U},T) := R\Hom_\cR(R\Gamma_\Sigma(\co_{F,U}, T^\ast(1)),\cR)[-3] \oplus \left(\bigoplus_{w \in S_\infty(F)}H^0(F_w, T) \right)[-1]. \nonumber
\end{eqnarray}
Then, as $p$ is odd, the Artin-Verdier duality theorem 
 combines with the triangle (\ref{Sigma-mod def}) with $N = T^*(1)$ to give an exact triangle in $D(\mathcal{R}[\cG_F])$
\begin{equation}\label{obvious tri}R\Gamma_{c,\Sigma}(\co_{F, U'},T) \to R\Gamma_{c,\Sigma}(\co_{F, U},T) \to \bigoplus_{w \in (U'\setminus U)_F}  R\Hom_\cR(R\Gamma_{/f}(F_w, T^*(1)),\cR)[-2]\to, \end{equation}
where $U'$ is any finite set of places of $K$ containing $U$ and disjoint from $\Sigma$.
By using this triangle, one shows that there is a canonical exact triangle in $D(\mathcal{R}[\cG_F])$
\begin{equation}\label{compact scalar change} R\Gamma_{c,\Sigma}(\co_{F, U'},T) \to R\Gamma_{c,\Sigma}(\co_{F, U},T) \to \bigoplus_{w \in (U'\setminus U)_F} R\Gamma_f(F_w, T)\to \end{equation}

\subsubsection{}\label{Sigma mod finite} For later purposes we also consider the following situation. Fix a non-zero element $M$ of $O$ and set $\fr:=\cR/M \cR$. We also set $A:=T/MT$ and for each $F$ in $\Omega(\cK/K)$ define a $\Sigma$-modified version of the compactly supported cohomology complex $R\Gamma_{c}(\co_{F,U},A)$ of $A$ by setting
\[ R\Gamma_{c,\Sigma}(\co_{F,U},A) := \fr \otimes_\cR^{\mathbb{L}}R\Gamma_{c,\Sigma}(\co_{F,U},T).\]
Then (\ref{obvious tri}) induces an exact triangle in $D(\fr[\cG_F])$ of the form
\begin{equation*}R\Gamma_{c}(\co_{F, U'},A) \to R\Gamma_{c,\Sigma}(\co_{F, U},A) \to \bigoplus_{w \in (U'\setminus U)_F}
 R\Hom_{\mathfrak{r}}(R\Gamma_{/f}(F_w, A^\ast(1)),\mathfrak{r})[-2]\to \end{equation*}
with $A^\ast(1):=\Hom_\fr(A, \fr(1))$, whilst for any finite set of places $U'$ of $K$ that contains $U$ and is disjoint from $\Sigma$ the triangle (\ref{compact scalar change}) induces an exact triangle in $D(\mathfrak{r}[\cG_F])$ 
\begin{equation}\label{compact scalar change2} R\Gamma_{c,\Sigma}(\co_{F, U'},A) \to R\Gamma_{c,\Sigma}(\co_{F, U},A) \to \bigoplus_{w \in (U'\setminus U)_F} R\Gamma_f(F_w, A)\to .\end{equation}

\begin{remark}
The definitions in the previous subsection will be replaced by those for $\Sigma$-modified cohomologies. In particular, an Euler system of rank $r$ will be an element in $\prod_{F\in \Omega(\cK/K)}{\bigcap}_{\cR[\cG_F]}^r H^1_\Sigma(\co_{F, S(F)}, T^\ast(1))$ which satisfies a norm compatible relation.
\end{remark}

\subsection{Vertical determinantal systems}

In the sequel we fix a finite set of places $\Sigma$ of $K$ as in \S \ref{section sigma} and, since this set is regarded as fixed, we often do not explicitly indicate it in our notation.

For each field $F$ in $\Omega(\cK/K)$ and a finite set $U$ of places of $K$ with $S(F) \subset U$, we set
$$C_{F,U}(T):=R\Hom_\cR(R\Gamma_{c,\Sigma}(\co_{F,U}, T), \cR)[-2] \in D^{\rm p}(\cR[\cG_F]).$$
(This complex is turned out to be perfect, see the proof of Proposition \ref{prop complex} below.)
%

In the sequel we often use, without explicit indication, the canonical identification
$${\det}_{\cR[\cG_F]}(C_{F,U}(T)) \simeq {\det}_{\cR[\cG_F]}^{-1}(R\Gamma_{c,\Sigma}(\co_{F,U},T))^\#.$$

\begin{definition} \label{definition vertical}
We define the $\cR[[\Gal(\cK/K)]]$-module of {\it vertical determinantal systems} for $(T,\cK)$ by setting
$${\rm VS}(T,\cK):={\varprojlim}_{F \in \Omega(\cK/K)} {\det}_{\cR[\cG_F]}(C_{F,S(F)}(T)),$$
where the limit is taken with respect to the following morphisms. For $F, F' \in \Omega(\cK/K)$ with $F \subset F'$, we take the transition morphism to be the composite surjective homomorphism
\begin{eqnarray}
 &&{\det}_{\cR[\cG_{F'}]}(C_{F',S(F')}(T))\nonumber\\
  &\to&  {\det}_{\cR[\cG_F]}(C_{F,S(F')}(T)) \nonumber\\
 &\simeq&  {\det}_{\cR[\cG_F]}(C_{F,S(F)}(T)) \otimes \bigotimes_{\fq \in S(F')\setminus S(F)} {\det}_{\cR[\cG_F]}(R\Gamma_f(K_\fq, T_F))^\# \nonumber\\
 &\simeq&  {\det}_{\cR[\cG_F]}(C_{F,S(F)}(T)), \nonumber
 \end{eqnarray}
 where the arrow denotes the norm map, the first isomorphism is induced by the triangle (\ref{compact scalar change}) with $N = T$, $U' = S(F')$ and $U= S(F)$, and the second by resolving each term $R\Gamma_f(K_\fq, T_F)$ as in (\ref{unram resolution}) and then using both the evaluation map
 $${\det}_{\cR[\cG_F]}(R\Gamma_f(K_\fq, T_F)) = {\det}_{\cR[\cG_F]}(T_F) \otimes {\det}_{\cR[\cG_F]}^{-1}(T_F) \stackrel{\sim}{\to} \cR[\cG_F]$$
 and the involution $\cR[\cG_F]^\# \stackrel{\sim}{\to} \cR[\cG_F]; \ \sigma \mapsto \sigma^{-1} $.

We shall refer to an element of ${\rm VS}(T,\cK)$ as a `vertical determinantal system' for $(T,\cK)$.
\end{definition}

\begin{remark}\label{motivating exam} An important motivating example of vertical determinantal systems is the collection of `zeta elements', whose existence is predicted by the equivariant Tamagawa number conjecture (eTNC, \cite[Conjecture 4]{BFetnc}).
To be more precise, suppose that the $p$-adic representation $T$ comes from a motive $\mathcal{M}$ defined over $K$. For simplicity, suppose that $O=\cR=\ZZ_p$ and $\Sigma=\emptyset$. For $F \in \Omega(\cK/K)$, the ($S(F)$-truncated) equivariant $L$-function of the motive $\mathcal{M}_F:=\mathcal{M} \otimes h^0(\Spec F)$ is defined by
$$L_S(\mathcal{M}_F, s):=\prod_{\fq \notin S(F)} P_\fq({\N}\fq^{-s} {\rm Fr}_\fq^{-1}) \in \CC_p[\cG_F],$$
where ${\N}\fq$ is the order of the residue field of $\fq$. Here we fix an isomorphism $\CC \simeq \CC_p$, and $L_S(\mathcal{M}_F,s)$ is regarded as a function on $\CC_p$. It is conjectured that $L_S(\mathcal{M}_F,s)$ is analytically continued to $s=0$. Assume this, and let $L_S^\ast(\mathcal{M}_F,0)\in \CC_p[\cG_F]^\times$ be the leading term at $s=0$. Then, the eTNC predicts the existence of a unique basis
$$\mathfrak{z}_F \in {\det}_{\ZZ_p[\cG_F]}(C_{F,S(F)}(T))$$
(called `zeta element'), which corresponds to $L_S^\ast(\mathcal{M}_F,0)$ under the `period-regulator isomorphism'
$$\CC_p \otimes_{\ZZ_p}{\det}_{\ZZ_p[\cG_F]}(C_{F,S(F)}(T)) \stackrel{\sim}{\to} \CC_p[\cG_F].$$
If we assume the eTNC for all $F \in \Omega(\cK/K)$, one can show that $\{\mathfrak{z}_F\}_F \in {\rm VS}(T,\cK)$.
\end{remark}

The following result plays an important role in the sequel.

\begin{proposition} \label{vertical rank one}
The $\cR[[\Gal(\cK/K)]]$-module ${\rm VS}(T,\cK)$ is free of rank one. \end{proposition}

\begin{proof} Fix an ordering of the (countably many) non-archimedean places of $K$ and for each natural number $n$ write $K_{(n)}$ for the maximal extension of $K$ in $\cK$ that has exponent dividing $p^n$ and is ramified only at the first $n$ non-archimedean places of $K$.

Then the finiteness of the ideal class group of $K$ implies that each field $K_{(n)}$ has finite degree over $K$. In addition, since for any $L$ in $\Omega(\cK/K)$ the extension $L/K$ is ramified at only finitely many places there exists a natural number $m$ such that $L \subset  K_{(m)}$.

For each $n$ set $\mathcal{G}_n := \cG_{K_{(n)}}$ and $\Xi_n := {\det}_{\cR[\cG_{n}]}(C_{K_{(n)},S(K_{(n)})}(T))$ and write $\tau_n$ for the transition morphism $\Xi_{n+1} \to \Xi_n$ used in the definition of ${\rm VS}(T,\cK)$. Note that each map $\tau_n$ is surjective and has kernel $I_{n+1}\cdot \Xi_{n+1}$ where $I_{n+1}$ denotes the ideal of $\mathcal{R}[\cG_{n+1}]$ that is generated by all elements of the form $h-1$ with $h \in \Gal(K_{(n+1)}/K_{(n)})$.

Assume now that for some fixed natural number $n$ we have fixed for each natural number $m$ with $m \le n$ an $\mathcal{R}[\cG_m]$-basis $x_m$ of
 $\Xi_m$ in such a way that $\tau_{m}(x_{m+1}) = x_m$ for all $m < n$. Let $x_{n+1}$ be any pre-image of $x_n$ under the surjective homomorphism
 $\tau_n: \Xi_{n+1} \to \Xi_n$. Then $\Xi_{n+1} = \mathcal{R}[\cG_{n+1}]\cdot x_{n+1} + I_{n+1}\cdot \Xi_{n+1}$ and hence, since $I_{n+1}$ is contained in the radical of $\mathcal{R}[\cG_{n+1}]$, Nakayama's Lemma implies that $x_{n+1}$ is an $\mathcal{R}[\cG_{n+1}]$-basis of $\Xi_{n+1}$.

Since the fields $\{K_{(n)}\}_{n \ge 1}$ are cofinal in $\Omega(\cK/K)$ (with respect to set-theoretic inclusion), the sequence of elements $(x_n)_{n \ge 1}$ that is inductively constructed in this way constitutes an $\cR[[\Gal(\cK/K)]]$-basis of ${\rm VS}(T,\cK)$. \end{proof}


\subsection{The module of basic Euler systems}\label{basic euler}

We suppose given a finite group $G$ and a Dedekind domain $O$ with quotient field $Q$ of characteristic prime to $\#G$. We denote the separable closure of $Q$ by $\overline Q$ and the set of $\overline Q$-valued characters of $G$ by $\widehat G_{\overline{Q}}$ (which we abbreviate to $\widehat G$ in the case that $\overline{Q}$ is an algebraic closure of $\QQ_p$).

For each $\chi \in \widehat G_{\overline{Q}}$, we define an idempotent
$$e_\chi:=\frac{1}{\#G} \sum_{\sigma \in G}\chi(\sigma)\sigma^{-1} \in \overline Q[G].$$

If $M$ is an $O[G]$-module, then for any idempotent $e$ of $Q[G]$ we define
\[ M[e] := \{ m \in M \mid e\cdot m = 0\} = \{ m \in M \mid m = (1-e)\cdot m\}.\]

\subsubsection{Faithful Euler systems} We fix notation as in \S\ref{def hres}. In the following, we assume

\begin{hypothesis} \label{hyp free}
\[ Y_K(T):=\bigoplus_{v\in S_\infty(K)}H^0(K_v,T)\]
is a free $\cR$-module.
\end{hypothesis}

We set
\[ r_T := {\rm rank}_{\cR}(Y_K(T)).\]

For each field $F$ in $\Omega(\cK/K)$ we set
\[ \widehat {\cG}_{F,T} := \{ \chi \in \widehat \cG_F \mid {\rm dim}_{\overline \QQ_p}(e_\chi(\overline \QQ_p\otimes_{O} H_\Sigma^1 (\co_{F,S(F)},T^\ast(1)))) = r_T\}\]
and then obtain an idempotent of $Q[\cG_F ]$ by setting
\[ e_{F,T}:=\sum_{\chi\in \widehat \cG_{F,T}} e_\chi .\]
%

\begin{definition}
An Euler system $\{c_F\}_F$ for $(T,\mathcal{K})$ is said to be {\it faithful} if it has rank $r_T$ and for every $F \in \Omega(\cK/K)$ one has
$$c_F \in \left({\bigcap}_{\cR[\cG_F]}^{r_T} H_\Sigma^1 (\co_{F,S(F)},T^\ast(1)) \right)[1-e_{F,T}].$$
The collection of faithful Euler systems for $(T,\mathcal{K})$ is an $\mathcal{R}[[\Gal(\mathcal{K}/K)]]$-submodule of ${\rm ES}_{r_T}(T,\cK)$ that we denote by $\mathcal{E}^{\rm f}(T,\cK)$.
\end{definition}

The structure of the module $\mathcal{E}^{\rm f}(T,\cK)$ will be described in Theorem \ref{det euler fitt} below.

\begin{remark}\label{faithful rubin-stark} We consider the case $K$ is totally real and $T = \ZZ_p$ (as already discussed in Remark \ref{remark rubin stark}) and set $r := [K:\QQ]$. Then Hypothesis \ref{hyp free} is automatic and one has $r_T = r$ and for each $F$ in $\Omega(\cK/K)$ the group $H_\Sigma^1 (\co_{F,S(F)},T^\ast(1))$ identifies with $\ZZ_p\otimes_\ZZ\mathcal{O}_{F,S(F),\Sigma}^\times$. In addition, since $S$ contains (all $p$-adic places and hence) at least one non-archimedean place, an analysis of $\Gamma$-factors in the functional equation of Dirichlet $L$-series implies that, after fixing an identification of $\widehat \cG_F$ with the set of complex characters of $\cG_F$, the $S$-truncated $L$-series of each character in $\widehat\cG_F\setminus \widehat\cG_{F,T}$ vanishes to order at least $r+1$ (see, for example, \cite[(6)]{rubinstark}). This implies that the $p$-completion of the lattice `$\Lambda_{S,T}$' defined by Rubin in \cite[\S2.1]{rubinstark} coincides with $\left({\bigcap}_{\ZZ_p[\cG_F]}^r H_\Sigma^1 (\co_{F,S(F)},\ZZ_p(1))\right)[1-e_{F,T}]$ and hence that the (conjectural) Euler system of Rubin-Stark elements should be faithful.\end{remark}

\begin{remark} \label{remark core rank}
The integer $r_T$ defined above provides a natural reinterpretation of the notion of `core rank' introduced by Mazur and Rubin
in \cite[Definition 4.1.11]{MRkoly}. To be specific, if one takes $K=\QQ$ and $\cR=O$, assumes $T$ is unramified outside $p$ and sets $S:=\{p, \infty\}$, then $H^1(\co_{K,S}, T^\ast(1))$ coincides with the group $H^1_{\mathcal{F}_{\rm can}}(\QQ,T^\ast(1))$, where $\mathcal{F}_{\rm can}$ is the `canonical Selmer structure' on $T^\ast(1)$ (see \cite[Definition 3.2.1]{MRkoly}). In this case, the core rank $\chi(T^\ast(1),\mathcal{F}_{\rm can})$ is equal, under suitable hypotheses (see \cite[Theorem 5.2.15]{MRkoly}), to
$${\rm rank}_{O}(T^\ast(1)^-) + {\rm corank}_O(H^0(\QQ_p, T\otimes_{\ZZ_p}\QQ_p/\ZZ_p)).$$
Moreover, one has ${\rm rank}_{O}(T^\ast(1)^-) = {\rm rank}_O (H^0(\RR, T))$ and, in many cases, the term ${\rm corank}_O(H^0(\QQ_p, T\otimes_{\ZZ_p}\QQ_p/\ZZ_p))$ vanishes (such examples are considered in \cite[\S6]{MRkoly}).
 In any such case, therefore, one has $\chi(T^\ast(1),\mathcal{F}_{\rm can}) = {\rm rank}_{O}(H^0(\RR,T))= r_T$.

Also, in more general cases, one of our working hypotheses in \S\ref{section unit system} is satisfied when $r_T$ is equal to the core rank. See Remark \ref{remark core vertex} for the detail.
\end{remark}




\subsubsection{Basic Euler systems} We now introduce the following hypothesis.

\begin{hypothesis} \label{hyp1} For every field $F$ in $\Omega(\cK/K)$ one has both
\begin{itemize}
\item[(i)] the module $H^1_\Sigma(\co_{F,S(F)},T^\ast(1))$ is free over $O$, and
\item[(ii)]  the module $H_\Sigma^0(F,T^\ast(1))$ vanishes.
\end{itemize}
\end{hypothesis}

\begin{remark} If $T$ is given, then one can always choose a set $\Sigma$ so that the condition (i) is satisfied (see Remark \ref{injectivity}). 
\end{remark}

We now state the main result of this section.

\begin{theorem} \label{theorem vertical} If Hypothesis \ref{hyp1} is satisfied, then there exists a canonical homomorphism of $\mathcal{R}[[\Gal(\mathcal{K}/K)]]$-modules
$\theta_{T,\cK}: {\rm VS}(T,\cK) \to \mathcal{E}^{\rm f}(T,\cK).$
\end{theorem}

This result will be proved in the next subsection. For the moment we use it to make the following definition.

\begin{definition} \label{def basic}
Assuming Hypothesis \ref{hyp1} to be satisfied, we define the module of {\it basic Euler systems for $(T,\cK)$} by setting
$$\mathcal{E}^{\rm b}(T,\mathcal{K}):=\im (\theta_{T,\mathcal{K}}) \subset  \mathcal{E}^{\rm f}(T,\mathcal{K}).$$
We shall say that an Euler system for $(T,\cK)$ is {\it basic} if it belongs to $\mathcal{E}^{\rm b}(T,\mathcal{K})$.
\end{definition}


It is clear $\mathcal{E}^{\rm b}(T,\mathcal{K})$ is an $\cR[[\Gal(\mathcal{K}/K)]]$-submodule of
$\mathcal{E}^{\rm f}(T,\mathcal{K})$ and in Theorem \ref{det euler fitt} below we prove that this module is in most cases free of rank one (see also Remark \ref{remark free}).

\subsection{Admissible complexes and the proof of Theorem \ref{theorem vertical}} In this section we introduce a general class of complexes with good properties with respect to exterior power biduals and then use them to prove Theorem \ref{theorem vertical}.

\subsubsection{Admissible complexes} Let $O$ be a discrete valuation ring with quotient field $Q$. In this subsection, we assume to be given a Gorenstein $O$-order $\cR$ in a finite dimensional separable commutative $Q$-algebra $\cQ$.


\begin{definition} \label{definition admissible}
A complex $C$ of $\cR$-modules is called {\it admissible} if the following conditions are satisfied:
\begin{itemize}
\item[(i)] $C$ is a perfect complex of $\cR$-modules;
\item[(ii)] the Euler characteristic of $Q \otimes_O C$ is zero;
\item[(iii)] $C$ is acyclic outside degrees zero and one;
\item[(iv)] $H^0(C)$ is $O$-free.
\end{itemize}
\end{definition}

%
%

The key property of the class of admissible complexes is the following result (which will be proved in Proposition \ref{prop admissible}).

\begin{proposition} \label{prop admissible new}

Let $\{e_1,\ldots,e_s \}$ be the complete set of the primitive orthogonal idempotents of $\cQ$ (so each $\cQ e_i$ is an extension field over $Q$ and we have the decomposition $\cQ=\bigoplus_{i=1}^s \cQ e_i$).
%
%
%
%
%
%
%
%
Let $\pi \in O$ be a uniformizer. Assume $(\cR/\pi\cR)/{\rm rad}(\cR/\pi \cR)$ is separable over $O/\pi O$, where ${\rm rad}$ denotes the Jacobson radical.
Let $C \in D^{\rm p}(\mathcal{R})$ be an admissible complex. Let $X$ be a free $\mathcal{R}$-module of rank $r$, and suppose that there is a surjective homomorphism $f: H^1(C) \to X.$ Let $e_r \in \cQ$ be the sum of the primitive idempotents which annihilate $Q \otimes_O \ker f$.
\begin{itemize}
\item[(i)] 
There exists a canonical homomorphism of $\cR$-modules
\[ \Pi_{C,f}: {\det}_{\cR}(C) \to \left({\bigcap}_{\cR}^r H^0(C)\right)[1-e_r] \otimes {\bigwedge}_{\cR}^r (X^\ast).\]

\item[(ii)] Let $I$ be a non-zero ideal of $O$. Set $C_I := C \otimes_O O/I$, $\cR_I:=\cR/I\cR$ and write $f_I$ for the homomorphism $H^1(C_I)\simeq H^1(C)\otimes_O O/I \to X\otimes_O O/I$ induced by $f$. Then there exists a canonical homomorphism of $\cR_I$-modules
\[ \Pi_{C_I,f_I}: {\det}_{\cR_I}(C_I) \to \left({\bigcap}_{\cR_I}^r H^0(C_I)\right)\otimes {\bigwedge}_{\cR_I}^r (X\otimes_O O/I)^\ast.\]
\end{itemize}
\end{proposition}

\subsubsection{}In this subsection we establish that Proposition \ref{prop admissible new} can be applied in the setting of the complexes introduced in \S\ref{section sigma}. We note that in our setting $(\cR/\pi\cR)/{\rm rad}(\cR/\pi\cR)$ is always separable over $O/\pi O$, since $O/\pi O$ is finite.

To do this we fix notation as \S\ref{def hres} and set $r := r_T$.
%
We fix a direct sum decomposition of free $\mathcal{R}$-modules
\begin{equation}\label{basis iso} Y_K(T)^* = \bigoplus_{i = 1}^{r}\mathcal{R}\cdot \beta_i.\end{equation}




\begin{proposition} \label{prop complex}
For any $F\in \Omega(\cK/K)$ and a finite set $U$ of places of $K$ such that $S(F) \subset U$ and $U\cap \Sigma =\emptyset$, we have the following.
\begin{itemize}
\item[(i)] $H^0(C_{F,U}(T))$ is canonically isomorphic to $H^1_\Sigma(\co_{F,U},T^\ast(1))$, and $H^1(C_{F,U}(T))$ lies in a split exact sequence
%
\begin{eqnarray}\label{c exact 1}
0 \to H^2_\Sigma(\co_{F,U},T^\ast(1)) \to H^1(C_{F,U}(T)) \to Y_K(T)^\ast \otimes_\cR \cR[\cG_F] \to 0
\end{eqnarray}
in which the first map is canonical and the second depends on the choice of a set of representatives of the orbits of $\Gal(\cK/K)$ on $S_\infty(\cK)$.

\item[(ii)] If Hypothesis \ref{hyp1} is satisfied, then $C_{F,U}(T)$ is an admissible complex of $\cR[ \cG_F ]$-modules.
\end{itemize}
\end{proposition}

\begin{proof} 

We note $C_{F,U}(T)$ is a perfect complex of $\cR[\cG_F]$-modules. For example, this follows from the fact that both terms
 $R\Gamma(\co_{F,U}, T^*(1))$ and $\bigoplus_{w\in \Sigma_F}R\Gamma_f(F_w,T^*(1))$ occuring in the relevant case of the exact triangle (\ref{Sigma-mod def}) are perfect complexes of $\cR[\cG_F]$-modules (the former by, for example, \cite[Proposition 1.6.5(2)]{FK} and the latter by virtue of the resolutions (\ref{unram resolution})).

We note that, since every archimedean place of $K$ splits completely in $\cK$, a choice of representatives of the $\Gal(\cK/K)$ on $S_\infty(\cK)$ induces for each $F$ in $\Omega(\cK/K)$ an isomorphism of $\mathcal{R}[\cG_F]$-modules
\begin{equation} \bigoplus_{w \in S_\infty(F)} H^0(F_w,T)\simeq Y_K(T)\otimes_\cR \cR[\cG_F].\end{equation}

The definition of $R\Gamma_{c,\Sigma}(\co_{F,U},T)$ in \S\ref{Sigma mod O} leads directly to an exact triangle in $D(\mathcal{R}[\cG_F])$ of the form
$$(Y_K(T)^\ast \otimes_\cR \cR[\cG_F])[-2] \to R\Gamma_\Sigma(\co_{F,U},T^\ast(1))[1] \to C_{F,U}(T).$$
From this and Hypothesis \ref{hyp1}(ii), we see that $C_{F,U}(T)$ is acyclic outside degrees zero and one, and obtain the canonical isomorphism
$$H^0(C_{F,U}(T)) \simeq H^1_\Sigma(\co_{F,U},T^\ast(1)),$$
and (\ref{c exact 1}) (note that $H^3_\Sigma(\co_{F,U}, T^\ast(1))$ vanishes since $p$ is odd).

\end{proof}

\subsubsection{The proof of Theorem \ref{theorem vertical}} 
%
%
%
For each $F$ in $\Omega(\cK/K)$ we write $\Pi_F$ for the composite homomorphism of $\cR[\cG_F]$-modules
\begin{eqnarray*}\label{defn Pi_F}
{\det}_{\cR[\cG_F]}(C_{F,S(F)}(T))\!\!\!\! &\to&\!\!\!\! \left({\bigcap}_{\cR[\cG_F]}^r H_\Sigma^1(\co_{F,S(F)},T^\ast(1))\right)\![1-e_{F,T}]\! \otimes\! {\bigwedge}_{\cR[\cG_F]}^r\! (Y_K(T)^\ast \otimes_\cR \cR[\cG_F]) \\
&\stackrel{\sim}{\to} & \left({\bigcap}_{\cR[\cG_F]}^r H_\Sigma^1(\co_{F,S(F)},T^\ast(1))\right)[1-e_{F,T}]. \nonumber
\end{eqnarray*}
Here the first arrow is the map $\Pi_{C,f}$ in Proposition \ref{prop admissible new}(i) with $C = C_{F,S(F)}(T)$ and $f$ the homomorphism
 $H^1(C_{F,S(F)}(T)) \to Y_K(T)^\ast \otimes_\cR \cR[\cG_F]$ in (\ref{c exact 1}), and the second arrow is induced by the isomorphism
$$ {\bigwedge}_{\cR[\cG_F]}^r (Y_K(T)^\ast \otimes_\cR \cR[\cG_F])  \stackrel{\sim}{\to} \cR[\cG_F]; \ \beta_1\wedge \cdots \wedge \beta_r \mapsto 1$$
that is induced by the decomposition (\ref{basis iso}).

For each element $z=\{z_F\}_F$ of ${\rm VS}(T,\cK)$ we then set
\begin{equation}\label{proj defn} c(z)_F := \Pi_F(z_F).\end{equation}
%

To show that the collection $c(z)=\{ c(z)_F\}_F$ belongs to $\mathcal{E}^{\rm f}(T,\cK)$ it suffices to show that for each pair of fields $F$ and $F'$ in $\Omega(T,\cK)$ with $F \subset F'$ and each $\chi$ in $\widehat\cG_F$ one has
\begin{equation}\label{es equality} e_\chi({\rm Cor}_{F'/F}(c(z)_{F'}))= e_\chi (P_{F'/F}\cdot c(z)_F)\end{equation}
with $P_{F'/F}:=\prod_{\fq \in S(F')\setminus S(F)}P_\fq({\rm Fr}_\fq^{-1}).$

From Lemma \ref{vanishing euler factor lemma} below it is enough to consider the case $e_\chi (P_{F'/F}) \not=0$. In this case the complex $e_\chi (\overline \QQ_p \otimes_O R\Gamma_f(K_\fq, T_F))$ is acyclic for each place $\fq$ in $S(F')\setminus S(F)$ and the required equality (\ref{es equality}) follows from the fact that the description of Proposition \ref{prop admissible}(ii) combines with Lemma \ref{lemma euler factor} below to imply the commutativity of
$$\xymatrix{
e_\chi(\overline \QQ_p \otimes_O{\det}_{\cR[\cG_{F'}]}(C_{F',S(F')}(T))) \ar[d] \ar[rr]^{\Pi_{F'}  \quad } &  & e_\chi(\overline \QQ_p \otimes_O{\bigcap}_{\cR[\cG_{F'}]}^rH_\Sigma^1(\co_{F',S(F')},T^\ast(1)))  \ar[d]^{{\rm Cor}_{F'/F}} \\
e_\chi(\overline \QQ_p \otimes_O{\det}_{\cR[\cG_F]}(C_{F,S(F)}(T)))  \ar[rr]_{P_{F'/F} \times \Pi_F  \quad } & & e_\chi(\overline \QQ_p \otimes_O{\bigcap}_{\cR[\cG_{F}]}^r H_\Sigma^1(\co_{F,S(F')},T^\ast(1))),
}
$$
where the left vertical arrow is the transition map in Definition \ref{definition vertical}.

Finally we note that, since ${\rm VS}(T,\cK)$ is a free rank one $\mathcal{R}[[\Gal(\mathcal{K}/K)]]$-module (by Proposition \ref{vertical rank one}), the assignment $z \mapsto c(z)$ gives a well-defined homomorphism $\theta_{T,\cK}$ of $\mathcal{R}[[\Gal(\mathcal{K}/K)]]$-modules of the required sort.

In the next result we regard $\widehat\cG_{F}$ as a subset of $\widehat\cG_{F'}$ in the obvious way.

\begin{lemma}\label{vanishing euler factor lemma} For each $\chi$ in $\widehat\cG_{F}$ one has $e_\chi(c(z)_{F'}) \not= 0$ if and only if both
  $e_\chi (P_{F'/F})\not=0$ and $e_\chi (c(z)_{F}) \not= 0$.
\end{lemma}

\begin{proof} The proof of Proposition \ref{prop admissible}(iii) implies that the $\cQ[\cG_{F'}]$-modules generated by $c(z)_{F'}$ and $c(z)_F$ are respectively isomorphic to $\cQ[\cG_{F'}]e_{F',T}$ and $\cQ[\cG_{F}]e_{F,T}$.

In particular, for each $\chi$ in $\widehat\cG_{F}$, and with $E$ denoting either $F$ or $F'$, one has
\begin{equation}\label{equiv} e_\chi(c(z)_{E}) \not= 0 \Longleftrightarrow \chi\in \widehat \cG_{E,T}\Longleftrightarrow  H^2_\Sigma(\co_{E,S(E)},T^\ast(1))_\chi = 0\end{equation}
where we write $N_\chi$ for the `$\chi$-component' $e_\chi(\overline \QQ_p \otimes_O N)$ of an $O[\cG_E]$-module $N$ and the last equivalence is a consequence of the exact sequence (\ref{c exact 1}) and the natural projection formula isomorphism $H^2_\Sigma(\co_{F',S(F')},T^\ast(1))_{\Gal(F'/F)} \simeq H^2_\Sigma(\co_{F,S(F')},T^\ast(1))$.

To study this condition we note that the $\cR$-linear dual of the exact triangle (\ref{compact scalar change}) induces, upon taking account of the exact sequence (\ref{c exact 1}) in the relevant cases, an exact sequence
 of $\cR[\cG_F]$-modules
\begin{equation*}\label{conv eq} H^2(C_{F'/F}) \xrightarrow{\theta} H^2_\Sigma(\co_{F,S(F)},T^\ast(1)) \to H^2_\Sigma(\co_{F,S(F')},T^\ast(1)) \to H^3(C_{F'/F}) \to 0\end{equation*}
in which we set $C_{F'/F}:=\bigoplus_{w \in (S(F')\setminus S(F))_F} R\Hom_{\cR}(R\Gamma_f(F_w, T),\cR)[-3]$. This sequence implies  $H^2_\Sigma(\co_{F,S(F')},T^\ast(1))_\chi=0$ if and only if both $H^3(C_{F'/F})_\chi=0$ and $\coker(\theta)_\chi=0$.

Next we note that, following on from (\ref{unram resolution}), the complex $C_{F'/F}$ is isomorphic to a complex of the form $M \to M$ for some finitely generated projective $\cR[\cG_F]$-module $M$ and therefore that $H^3(C_{F'/F})_\chi=0$ if and only if $H^2(C_{F'/F})_\chi=0$ in which case one has $e_\chi(P_{F'/F})\not= 0$.

At this stage we know that $e_\chi(c(z)_{F'}) \not= 0$ if and only if both $e_\chi(P_{F'/F})\not= 0$ and
$H^2_\Sigma(\co_{F,S(F)},T^\ast(1))_\chi = \coker(\theta)_\chi=0$ so that $e_\chi(c(z)_{F}) \not= 0$. \end{proof}

\begin{lemma} \label{lemma euler factor}
 Fix a field $F \in \Omega(\cK/K)$, a place $\fq$ of $F$ outside $S(F)$ and a character $\chi$ of $\cG_F$ for which
 the complex $e_\chi (\overline \QQ_p \otimes_O R\Gamma_f(K_\fq, T_F))$ is acyclic.

Then $e_\chi P_\fq({\rm Fr}_\fq)$ is equal to the image of the canonical basis of
$$ {\det}_{\cR[\cG_F]}(R\Gamma_f(K_\fq, T_F))= {\det}_{\cR[\cG_F]}(T_F) \otimes {\det}_{\cR[\cG_F]}^{-1}(T_F) $$
under the map
\[
e_\chi(\overline \QQ_p \otimes_O {\det}_{\cR[\cG_F]}(R\Gamma_f(K_\fq, T_F)))
\simeq {\det}_{e_\chi\overline \QQ_p[\cG_F]}(0) \otimes {\det}_{e_\chi \overline \QQ_p[\cG_F]}^{-1}(0) = e_\chi \overline \QQ_p[\cG_F ]. \nonumber
\]
%
\end{lemma}


\begin{proof} This is well-known (and also easy to verify by explicit computation). \end{proof}

\subsection{Structures of modules of higher rank Euler systems} In this section we study the explicit structures of the various modules of Euler systems introduced in \S\ref{basic euler}. To do this we must first introduce some notation.

Note that if $F$ and $F'$ belong to $\Omega(\mathcal{K}/K)$ and $F \subset  F'$, then Lemma \ref{vanishing euler factor lemma} implies the natural projection map $\pi_{F'/F}: \cQ[\mathcal{G}_{F'}] \to \cQ[\mathcal{G}_{F}]$ maps $\cR[\mathcal{G}_{F'}][e_{F',T}]$ to $\cR[\mathcal{G}_{F}][e_{F,T}]$ and hence induces a natural homomorphism  $\cR[\mathcal{G}_{F'}]e_{F',T} \to \cR[\mathcal{G}_{F}]e_{F,T}$.
%
%

Given an $\cR[\mathcal{G}_F]$-submodule $X_F$ of $\cQ[\mathcal{G}_F]e_{F,T}$ for each $F$ in $\Omega(\cK/K)$ we define an $\cR[[\Gal(\mathcal{K}/K)]]$-submodule of $\prod_{F\in \Omega(\cK/K)}\cQ[\cG_F]$ by setting
\[ {\varprojlim}^\diamond_{F\in \Omega(\cK/K)} X_F := \left\{ (x_F)_F\! \in\!\prod_{F\in \Omega(\cK/K)}\!\! \widetilde X_F \ \middle| \ P_{F'/F}e_{F,T}\cdot (\pi_{F'/F}(x_{F'}) - x_{F})  = 0  \text{ for all } F\subset  F'\right\}.\]
with $\widetilde X_F := \cQ[\cG_F][e_{F,T}] + X_F$.

\begin{remark}\label{funny lim remark} One has
%
\[ {\varprojlim}_{F'}\left({{\bigcap}}_{F \in \Omega(F'/K)}\pi_{F'/F}^{-1}(\widetilde X_{F})\right)  = \varprojlim_{F'} \cQ[\cG_{F'}] \cap \prod_{F'} \widetilde X_{F'} \subset  {\varprojlim}^\diamond_{F'}X_{F'} \]
where in all cases $F'$ runs over $\Omega(\cK/K)$, the inverse limits are taken with respect to the maps $\pi_{F'/F}$ and we write ${\pi}_{F'/F}^{-1}(Y)$ for the set-theoretic pre-image under $\pi_{F'/F}$ of a subset $Y$ of $\cQ[\cG_{F}]$. It is also clear that this inclusion is an equality for any representation $T$ for which for each $F$ and $F'$ in $\Omega(\cK/K)$ with $F\subset  F'$ the product $P_{F'/F}e_{F,T}$ is invertible in $\cQ[\mathcal{G}_{F}]$. We further note that this latter condition is automatically satisfied if $H_\Sigma^2(\mathcal{O}_{F',S(F')},T^*(1))$ is finite (and hence $e_{F',T} =1$) for all $F'$ in $\Omega(\cK/K)$. \end{remark}

We set $r:=r_T$, $W_r(T,\cK) := \prod_{F \in \Omega(\cK/K)}\cQ[\cG_F][e_{F,T}]$ and for each subgroup $\mathcal{X}$ of $\prod_{F \in \Omega(\cK/K)}\cQ[\cG_F]$ define a quotient
\[ \mathcal{X}_{(r)} := \mathcal{X}/(\mathcal{X}\cap W_r(T,\cK)).\]

\begin{remark}\label{tautological id} The tautological map $\pi_{\mathcal{X},r}: \mathcal{X}\to \mathcal{X}_{(r)}$ is obviously bijective if $e_{F,T} = 1$ for all $F$ in $\Omega(\cK/K)$ (so $W_r(T,\cK) = 0$), as is the case for the representations $T$ discussed at the end of Remark \ref{funny lim remark}. There are also other interesting cases in which $\pi_{\mathcal{X},r}$ can be shown to be bijective. For example, if $T = \ZZ_p$ and $\cK$ contains the cyclotomic extension of $K$, then it is straightforward to see that $\pi_{\mathcal{X},r}$ is bijective for any submodule $\mathcal{X}$ of $\ZZ_p[[\Gal(\cK/K)]]$.
\end{remark}

%

Claim (ii) of the following result is a natural generalization of \cite[Theorem 7.5]{bks1}.

\begin{theorem}\label{det euler fitt}\
Assume Hypothesis \ref{hyp1}.
\begin{itemize}
\item[(i)] There exists a commutative diagram of $\mathcal{R}[[\Gal(\mathcal{K}/K)]]$-module homomorphisms
%
\begin{equation}\label{structure diag}\xymatrix{ \mathcal{E}^{\rm f}(T,\mathcal{K}) \ar[r]^{\!\!\!\!\!\!\!\!\!\!\!\!\!\!\!\!\!\!\!\!\!\!\!\!\!\!\!\!\!\!\!\!\!\!\!\!\!\!\!\!\!\!\!\!\!\!\!\!\!\!\!\!\!\!\!\!\!\!\!\!\!\!\!\!\!\!\sim}  & ({\varprojlim}^\diamond_{L\in \Omega(\cK/K)} {\rm Fitt}^0_{\mathcal{R}[\mathcal{G}_{L}]}(H^2_\Sigma(\mathcal{O}_{L,S(L)},T^*(1)))^{-1})_{(r)}\\
\mathcal{E}^{\rm b}(T,\mathcal{K}) \ar@{^{(}->}[u] \ar[r]^{\!\!\!\!\!\!\!\!\!\!\!\!\!\!\sim}   & \mathcal{R}[[\Gal(\mathcal{K}/K)]]_{(r)} \ar@{^{(}->}[u] }
\end{equation}
in which the  horizontal maps are bijective, the left hand vertical map is the tautological inclusion and the right hand vertical map is a natural injective homomorphism. (Here for an $\cR[\cG_L]$-submodule $X$ of $\cQ[\cG_L]$ we define $X^{-1}$ as in \S \ref{A proof}, replacing $\cR$ in loc. cit. by $\cR[\cG_L]$.)

\item[(ii)] For each $F$ in $\Omega(\cK/K)$ one has
$$ {\rm Fitt}_{\cR[\cG_F]}^0(H_\Sigma^2(\co_{F,S(F)},T^\ast(1)))=\langle \im (c_F) \mid c \in \mathcal{E}^{\rm b}(T,\mathcal{K}) \rangle_{\mathcal{R}[\cG_F]},$$
where we use the fact that each element of ${\bigcap}_{\mathcal{R}[\cG_F]}^r H_\Sigma^1(\co_{F,S(F)},T^\ast(1))$ is, by definition, a homomorphism
 with values in $\cR[\cG_F].$
\end{itemize}
\end{theorem}

\begin{proof} To construct the diagram in claim (i) we fix, following Proposition \ref{vertical rank one}, a basis of ${\rm VS}(T,\cK)$ as an $\mathcal{R}[[\Gal(\cK/K)]]$-module and write $\varepsilon = \{ \varepsilon_L \}_{L \in \Omega(\cK/K)}$ for its image in $\mathcal{E}^{\rm f}(T,\cK)$ under the homomorphism $\theta_{T,\cK}$ constructed in Theorem \ref{theorem vertical}.

Then the argument of Proposition \ref{prop admissible}(iii) shows each element $\varepsilon_L$ spans a $\cQ[\cG_L]$-module that is isomorphic to $\cQ[\cG_L]e_{L,T}$ and hence that $\ker(\theta_{T,\cK}) = \cR[[\Gal(\cK/K)]]\cap W_r(T,\cK)$.

We can therefore take the lower horizontal arrow in (\ref{structure diag}) to be the inverse of the isomorphism
\[ \cR[[\Gal(\mathcal{K}/K)]]_{(r)} = \cR[[\Gal(\mathcal{K}/K)]]/\ker(\theta_{T,\cK}) \simeq \im(\theta_{T,\cK}) =: \mathcal{E}^{\rm b}(T,\mathcal{K})\]
induced by $\theta_{T,\cK}$.

%
%

Next we note that the argument of Proposition \ref{prop admissible}(iii) also shows that for any element $\eta = \{\eta_L\}_{L \in \Omega(\cK/K)}$ of $\mathcal{E}^{\rm f}(T,\cK)$ and any $L$ in $\Omega(\cK/K)$ one has
\begin{equation}\label{comparison}\eta_L = x_L\cdot \varepsilon_L \end{equation}
for a unique element $x_L$ of $X_L:= {\rm Fitt}^0_{\mathcal{R}[\mathcal{G}_{L}]}(H^2_\Sigma(\mathcal{O}_{L,S(L)},T^*(1)))^{-1} \subset \cQ[\mathcal{G}_L]e_{L,T}$. In addition, as $\eta$ and $\varepsilon$ both belong to ${\rm ES}_r(T,\cK)$, for each $L' \in \Omega(\cK/K)$ and $L \in \Omega(L'/K)$
one has
\begin{multline*} P_{L'/L}x_{L}\cdot \varepsilon_{L} =  P_{L'/L}\cdot \eta_{L} = {\rm Cor}_{L'/L}(\eta_{L'}) = {\rm Cor}_{L'/L}(x_{L'}\cdot \varepsilon_{L'})\\ = \pi_{L'/L}(x_{L'}){\rm Cor}_{L'/L}(\varepsilon_{L'}) =  P_{L'/L}\pi_{L'/L}(x_{L'})\cdot\varepsilon_{L}\end{multline*}
and hence also $P_{L'/L}x_{L} = P_{L'/L}\pi_{L'/L}(x_{L'})$.

This shows that the assignment $(y_L)_L \mapsto (y_L\cdot \varepsilon_L)_L$ gives a well-defined surjective homomorphism of $\cR[[\Gal(\cK/K)]]$-modules $\theta'_{T,\cK}: {\varprojlim}^\diamond_{L\in \Omega(\cK/K)} X_L \to
\mathcal{E}^{\rm f}(T,\cK).$

The kernel of $\theta'_{T,\cK}$ is $({\varprojlim}^\diamond_{L\in \Omega(\cK/K)} X_L) \cap W_r(T,\cK)$ and we take the lower horizontal arrow in (\ref{structure diag}) to be the inverse of the isomorphism $({\varprojlim}^\diamond_{L\in \Omega(\cK/K)} X_L)_{(r)} \simeq \mathcal{E}^{\rm f}(T,\cK)$  that is induced by $\theta'_{T,\cK}$.

Finally we take the right hand vertical arrow in (\ref{structure diag}) to be the injective homomorphism induced by the obvious inclusion $\cR[[\Gal(\cK/K)]] \subset  {\varprojlim}^\diamond_{L\in \Omega(\cK/K)} X_L$ and note that, with this definition,
commutativity of the diagram is clear.


Claim (ii) follows directly from Lemma \ref{lemma standard}(iii) and the definition of basic Euler systems.\end{proof}

\begin{remark}\label{remark free} Theorem \ref{det euler fitt}(i) implies that for any of the pairs $(T,\cK)$ discussed in Remark \ref{tautological id} the $\mathcal{R}[[\Gal(\cK/K)]]$-module $\mathcal{E}^{\rm b}(T,\cK)$ is free of rank one. \end{remark}

\begin{remark}\label{subsequent article} In a subsequent article we will explain both how Proposition \ref{vertical rank one} and Theorems \ref{theorem vertical} and \ref{det euler fitt} can be extended to the context of non-abelian extensions of number fields and also how, in certain circumstances, they lead to the construction of new Euler systems with prescribed properties.\end{remark}

\section{Equivariant Stark systems} \label{section unit system}

In this section we develop an `equivariant' refinement of the theory of {\it Stark systems} that was introduced by Mazur and Rubin in \cite{MRselmer} and by the second author in \cite{sanojnt} (where the terminology `unit system' is used). To overcome certain algebraic difficulties (which forced previous construction to be restricted to modules over principal ideal rings) we systematically use exterior power biduals in place of exterior powers.

In particular, in this way we describe an explicit, and very natural, construction of Stark systems from certain systems in the determinants of cohomology that we call {\it horizontal determinantal systems} (see Definition \ref{def horizon} and Theorem \ref{theorem HU}). This allows us firstly to prove that, under certain  natural conditions, the module of Stark systems is free of rank one (see Theorems \ref{structure unit} and \ref{structure unit 2}), and then to deduce several explicit results on the structure of Selmer groups, which can each be regarded as natural equivariant generalizations of results of Mazur and Rubin in \cite[\S 8]{MRselmer} (see, in particular, Theorems \ref{det unit fitt} and \ref{theorem higher fitt}). 

In \S\ref{eks} we will then explain how this approach allows one to attach a canonical `regulator-type equivariant Kolyvagin system' from an equivariant Stark system, thus extending the key results of \cite[Theorem 5.7]{sanojnt} and \cite[Proposition 12.3]{MRselmer} (and see also Howard \cite[Appendix B]{MRkoly}).


\subsection{The definition} \label{section def unit}

Let $K$ be a number field. Let $R$ be a self-injective (commutative) ring. Let $A$ be a finite $R[G_K]$-module (endowed with discrete topology). We fix a Selmer structure $\mathcal{F}$ on $A$. Recall that a Selmer structure $\mathcal{F}$ on $A$ is a collection of the following data:
\begin{itemize}
\item a finite set $S(\mathcal{F})$ of places of $K$ such that $S_\infty(K)\cup \{ v \mid \# A\} \cup S_{\rm ram}(A) \subset S(\mathcal{F})$;
\item for every $v \in S(\mathcal{F})$, a choice of $R$-submodule $H_\mathcal{F}^1(K_v, A) \subset H^1(K_v, A)$.
\end{itemize}
We fix a finite set $\Sigma$ of places of $K$ which is disjoint from $S(\mathcal{F})$. The ($\Sigma$-modified) Selmer module attached to $\mathcal{F}$ is defined by
$$H_\mathcal{F}^1(K, A):=\ker\left( H_\Sigma^1(\co_{K,S(\mathcal{F})}, A) \to \bigoplus_{v\in S(\mathcal{F})} H^1(K_v, A)/H^1_{\mathcal{F}}(K_v,A)\right).$$

Let $\mathcal{P}$ be a set of places of $K$ such that
\begin{itemize}
\item $\mathcal{P}$ is disjoint from $S(\mathcal{F})\cup \Sigma$;
\item for every $\fq \in \mathcal{P}$, $H^1_{/f}(K_\fq, A)$ is isomorphic to $R$.
\end{itemize}
We assume that $\mathcal{P}$ is non-empty. Let $\cN=\cN(\cP)$ be the set of square-free products of primes in $\mathcal{P}$. For $\fn \in \cN$, we define the Selmer structure $\mathcal{F}^\fn$ by
\begin{itemize}
\item $S(\mathcal{F}^\fn):=S(\mathcal{F}) \cup \{ \fq \mid \fn\}$;
\item for $v \in S(\mathcal{F}^\fn)$,
$$H^1_\mathcal{F}(K_v,A) := \begin{cases}
H^1_\cF(K_v, A) &\text{ if $v \in S(\mathcal{F})$},\\
H^1(K_v , A) &\text{ if $v\mid \fn$}.
\end{cases}$$
\end{itemize}
The number of primes which divide $\fn$ is denoted by $\nu(\fn)$.

Let $r$ be a non-negative integer. To define Stark systems of rank $r$, we shall define a map
$$v_{\fm,\fn}: {\bigcap}_{R}^{r + \nu(\fm)} H^1_{\cF^\fm} (K, A) \to  {\bigcap}_{R}^{r + \nu(\fn)} H_{\cF^\fn}^1 (K, A)$$
for $\fm, \fn\in \cN$ with $\fn\mid \fm$.
For $\fq \in \cP$, we fix an isomorphism $H^1_{/f}(K_\fq, A)\simeq R$ and define a map $v_\fq$ by the composition
$$v_\fq: H^1(K,A) \to H^1(K_\fq, A) \to H_{/f}^1(K_\fq, A) \simeq R. $$
We have the exact sequence
$$0 \to H_{\cF^\fn}^1(K, A) \to H_{\cF^\fm}^1(K,A) \stackrel{\bigoplus_{\fq \mid \fm / \fn} v_\fq }{\to} R^{\oplus \nu(\fm / \fn)}.$$
By Proposition \ref{prop self injective}, we see that ${\bigwedge}_{\fq \mid \fm/\fn}v_\fq$ induces
$${\bigwedge}_{\fq \mid \fm/\fn}v_\fq : {\bigcap}_R^{r+\nu(\fm)} H_{\cF^\fm}^1(K, A) \to {\bigcap}_R^{r+\nu(\fn)} H_{\cF^\fn}^1(K, A).$$
Here the order of ${\bigwedge}_{\fq \mid \fm/\fn}$ is determined by fixing an injection $\cP \hookrightarrow \ZZ$ and regarding $\cP$ as a totally ordered set $(\cP, \prec)$. We write

$$\fm/\fn \text{ (resp. $\fn$, resp. $\fm$)}=\prod_{i=1}^{\nu(\fm/\fn)} \fq_i \left(\text{resp. $\prod_{i=1}^{\nu(\fn)} \fq'_i$, resp. $\prod_{i=1}^{\nu(\fm)}\fq_i''$}\right)$$
so that $\fq_i\prec \fq_j$ (resp. $\fq_i' \prec \fq_j'$, resp. $\fq_i''\prec \fq_j''$) if $i < j$.
We define
$$v_{\fm,\fn}:={\rm sgn}(\fm, \fn) \cdot {\bigwedge}_{\fq \mid \fm/\fn}v_\fq ={\rm sgn}(\fm, \fn) \cdot v_{\fq_1} \wedge \cdots \wedge v_{\fq_{\nu(\fm/\fn)}} ,$$
where ${\rm sgn}(\fm, \fn)$ is defined to be the sign of the permutation
$$(\fq_1, \ldots, \fq_{\nu(\fm/\fn)}, \fq'_1,\ldots , \fq'_{\nu(\fn)}) \mapsto (\fq''_1,\ldots, \fq''_{\nu(\fm)}).$$
One easily checks that $v_{\fm', \fn}=v_{\fm, \fn} \circ v_{\fm',\fm}$ if $\fn \mid \fm \mid \fm'$.

\begin{definition} \label{def unit}
A {\it Stark system of rank $r$} for $A$ is a collection
$$\left\{ \epsilon_\fn \in {\bigcap}_{R}^{r + \nu(\fn)} H_{\cF^\fn}^1(K, A) \ \middle| \ \fn \in \cN \right\}$$
which satisfies
$$v_{\fm,\fn}(\epsilon_\fm)=\epsilon_{\fn},$$
namely, a Stark system is an element of
$$\varprojlim_{\fn \in \cN} {\bigcap}_{R}^{r + \nu(\fn)} H_{\cF^\fn}^1(K, A),$$
where the inverse limit is taken with respect to the transition map $v_{\fm,\fn}$.

\end{definition}

\subsection{Horizontal determinantal systems} \label{section horizontal}

We set notations, which will be used in the rest of this paper. We keep the notations used in \S \ref{section 2}. As in \S\ref{Sigma mod finite} we fix a non-zero element $M\in O$, and put $\mathfrak{r}:=\cR/M \cR$. Let $\overline M$ be the smallest power of $p$ which is divisible by $M$. We also fix $E \in \Omega(\cK/K)$, which is unramified outside $S$. Put $A:=T/MT$, $\cA:=A_E={\rm Ind}_{G_K}^{G_E}(A)$, $\Gamma:=\cG_E(=\Gal(E/K))$, and $R:=\mathfrak{r}[\Gamma]$. Note that, since $\cR$ is a one-dimensional Gorenstein ring, $R$ is a self-injective ring (namely, a zero-dimensional Gorenstein ring).
We set $\cA^\ast(1):=\Hom_R(\cA, R(1))$, where $R(1):=R \otimes_{\ZZ_p}\ZZ_p(1)$. 
We often use the following fact (the `local duality'): for any prime $\fq $ of $K$ and $i \in \{0,1,2\}$, the cup-product  paring
$$H^i(K_\fq, \cA^\ast(1)) \times H^{2-i}(K_\fq,\cA) \stackrel{\cup}{\to} H^2(K_\fq, R(1))\simeq R$$
induces the isomorphism
\begin{eqnarray} \label{LD}
H^i(K_\fq,\cA^\ast(1)) \simeq H^{2-i}(K_\fq,\cA)^\ast.
\end{eqnarray}
Furthermore, if $\fq\notin S_{\rm min}(\cA)$, then (\ref{LD}) induces the isomorphisms
\begin{eqnarray} \label{LDf}
H^i_f(K_\fq,\cA^\ast(1)) \simeq H^{2-i}_{/f}(K_\fq,\cA)^\ast
\end{eqnarray}
and
\begin{eqnarray} \label{LD/f}
H^i_{/f}(K_\fq,\cA^\ast(1)) \simeq H^{2-i}_{f}(K_\fq,\cA)^\ast
\end{eqnarray}
(See \cite[Theorem 1.4.1 and Proposition 1.4.3(ii)]{R}.)

We denote by $K(1)$ the maximal $p$-extension of $K$ inside the Hilbert class field of $K$.


Let $\mathcal{P}$ be the set of primes $\fq\notin S\cup \Sigma$ of $K$ satisfying the following:
\begin{itemize}
\item $\fq$ splits completely in $E_M:=E(\mu_{\overline M}, (\co_{K}^\times)^{1/\overline M})K(1)$;
\item $\cA/({\rm Fr}_\fq-1)\cA \simeq R$ as $R$-modules.
\end{itemize}

We recall some basic facts from \cite[Lemma 1.2.3]{MRkoly}. First, we have
$$P_\fq(1) \equiv 0 \text{ (mod $M$)}.$$
Here recall that
$$P_\fq(x)=\det(1-{\rm Fr}_\fq^{-1}x \mid T) \in \cR[x].$$
For every $\fq \in \cP$, there exists a unique polynomial $Q_\fq(x) \in \fr[x]$ such that
$$(x-1)Q_\fq(x)\equiv P_\fq(x) \text{ (mod $M$)}.$$
We have the isomorphism
\begin{eqnarray}
Q_\fq({\rm Fr}_\fq): \cA/({\rm Fr}_\fq-1)\cA \stackrel{\sim}{\to} \cA^{{\rm Fr}_\fq=1}; \ \overline a \mapsto Q_\fq({\rm Fr}_\fq)a, \label{q isom 2}
\end{eqnarray}
where $\cA^{{\rm Fr}_\fq=1}:=H^0(K_\fq,\cA)$.
In particular, note that
\begin{eqnarray}
H^2(K_\fq, \cA^\ast(1))\stackrel{(\ref{LD})}{\simeq} (\cA^{{\rm Fr}_\fq=1})^\ast \simeq R.  \label{identifications}
\end{eqnarray}
{\it We fix the last isomorphism, and often regard $H^2(K_\fq, \cA^\ast(1))=(\cA^{{\rm Fr}_\fq=1})^\ast=R$. }The element of $H^2(K_\fq, \cA^\ast(1))$ which corresponds to $1 \in R$ under this identification is denoted by $\gamma_\fq$.

We define a Selmer structure $\cF $ on $\cA^\ast(1)$ by
\begin{itemize}
\item $S(\cF):=S$;
\item for $v \in S(\cF)$, $H_\cF^1(K_v, \cA^\ast(1)):=H^1(K_v, \cA^\ast(1))$.
\end{itemize}
Note that
$$H^1_\cF(K,\cA^\ast(1))=H^1_\Sigma(\co_{K,S}, \cA^\ast(1)). $$

Let $\cN=\cN(\cP)$ be the set of square-free products of primes in $\mathcal{P}$.
For $\fn \in \cN$, we set
$$S_\fn:=S \cup \{ \fq \mid \fn\}. $$
Note that
$$H^1_{\cF^\fn}(K,\cA^\ast(1))=H^1_\Sigma(\co_{K,S_\fn}, \cA^\ast(1)). $$

For $\fq \in \cP$, we define
$$w_\fq: H^1_{/f}(K_\fq, \cA^\ast(1)) \stackrel{(\ref{LD/f})}{\simeq}  H_f^1(K_\fq,\cA)^\ast \simeq (\cA/({\rm Fr}_\fq-1)\cA)^\ast \stackrel{-Q_\fq({\rm Fr}_\fq)^{-1,\ast}}{\simeq} (\cA^{{\rm Fr}_\fq=1})^\ast  \stackrel{(\ref{identifications})}{=} R,$$
where the second isomorphism is the dual of the inverse of the isomorphism
\begin{eqnarray}
H^1_f(K_\fq,\cA) \stackrel{\sim}{\to} \cA/({\rm Fr}_\fq-1)\cA; \ a \mapsto a({\rm Fr}_\fq)\label{frobenius isom}
\end{eqnarray}
($a$ is regarded as a $1$-cocycle, and $a({\rm Fr}_\fq)$ means the evaluation of ${\rm Fr}_\fq$, see \cite[Lemma 1.3.2]{R}), and the third is the dual of $-Q_\fq({\rm Fr}_\fq)^{-1}$ (see (\ref{q isom 2})).
So we see that every $\fq \in \cP$ satisfies $H^1_{/f}(K_\fq, \cA^\ast(1)) \simeq R$. We fix this isomorphism via $w_\fq$, and we can define the module of Stark systems
$${\rm SS}_r(\cA):=\varprojlim_{\fn \in \cN} {\bigcap}_R^{r+\nu(\fn)} H^1_{\cF^\fn}(K,\cA^\ast(1))=\varprojlim_{\fn \in \cN} {\bigcap}_R^{r+\nu(\fn)} H^1_\Sigma(\co_{K,S_\fn},\cA^\ast(1)).$$

For $\fn \in \cN$, we set
\[ C_{K,S_\fn}(\cA):=C_{K,S_\fn}(T_E) \otimes_O (O/MO) \in D^{\rm p}(R).\]

As a variant of vertical determinant systems, we give the following definition.

\begin{definition} \label{def horizon}
A {\it horizontal determinantal system} for $\cA$ is an element of
$${\rm HS}(\cA):=\varprojlim_{\fn \in \cN}{\det}_R(C_{K,S_\fn}(\cA)),$$
where the transition map
\begin{eqnarray*} \label{transition horizon}
{\det}_R(C_{K,S_\fm}(\cA)) \stackrel{\sim}{\to} {\det}_R(C_{K,S_\fn}(\cA))
\end{eqnarray*}
is induced by the exact triangle
$$C_{K,S_\fn}(\cA) \to C_{K,S_\fm}(\cA) \to \bigoplus_{\fq \mid \fm/\fn}R\Gamma_{/f}(K_\fq,\cA^\ast(1))[1]$$
and the isomorphism
\begin{eqnarray}
{\det}_R(R\Gamma_{/f}(K_\fq,\cA^\ast(1)))&\simeq& H_{/f}^1(K_\fq,\cA^\ast(1)) \otimes_R H^2(K_\fq,\cA^\ast(1))^\ast  \nonumber \\
&\stackrel{w_\fq \text{ and (\ref{identifications})}}{\stackrel{\sim}{\to}}& (\cA^{{\rm Fr}_\fq=1})^\ast \otimes_R \cA^{{\rm Fr}_\fq=1} \nonumber \\
&\stackrel{\sim}{\to}& R, \nonumber
\end{eqnarray}
where the last isomorphism is the evaluation.

Note that ${\rm HS}(\cA)$ is isomorphic to $R$, by definition.
\end{definition}

{\it In the following, we assume Hypotheses \ref{hyp free} and \ref{hyp1}.}

We give some remarks, which will be needed in the proof of Theorem \ref{theorem HU} below.
Recall that we defined $\gamma_\fq \in H^2(K_\fq,\cA^\ast(1))$ to be the element which corresponds to $1 \in R$ under the identification (\ref{identifications}). Set $Y_K(\cA):=Y_K(T)\otimes_\cR R$. We denote the image of $\beta_i \in Y_K(T)$ in $Y_K(\cA)$ also by $\beta_i$. By (\ref{c exact 1}), we have a split exact sequence
$$0 \to H_\Sigma^2(\co_{K,S_\fn},\cA^\ast(1)) \to H^1(C_{K,S_\fn}(\cA)) \to Y_K(\cA)^\ast \to 0.$$

\begin{theorem} \label{theorem HU}
Set $r:=r_T$. Then there is a canonical homomorphism of $R$-modules
\begin{eqnarray} \label{horizontal unit}
{\rm HS}(\cA) \to {\rm SS}_r(\cA).
\end{eqnarray}
\end{theorem}

\begin{proof}

Take $\fn \in \cN$, and write $\fn = \prod_{i=1}^{\nu(\fn)}\fq_i$ so that $\fq_1 \prec \cdots \prec \fq_{\nu(\fn)}$.
By Proposition \ref{prop admissible}(iv), there exists a quadratic standard representative $(P_\fn \stackrel{\psi}{\to} P_\fn)$ of $C_{K,S_\fn}(\cA)$ with respect to
$$H^1(C_{K,S_\fn}(\cA)) \to \left( \bigoplus_{\fq \mid \fn}H^2(K_\fq,\cA^\ast(1)) \right) \oplus Y_K(\cA)^\ast=:X_\fn,$$
where $H^1(C_{K,S_\fn}(\cA)) \to \bigoplus_{\fq \mid \fn}H^2(K_\fq,\cA^\ast(1)) $ is induced by the localization morphism $C_{K,S_\fn}(\cA) \to \bigoplus_{\fq \mid \fn}R\Gamma_{/f}(K_\fq,\cA^\ast(1))[1]$.
For a horizontal determinantal system $z = \{z_\fn\}_\fn \in {\rm HS}(\cA)$, we define
$$\epsilon(z)_\fn \in {\bigcap}_R^{r+\nu(\fn)} H_\Sigma^1(\co_{K,S_\fn},\cA^\ast(1))$$
to be the image of $z_\fn \in {\det}_R(C_{K,S_\fn}(\cA))$ under the map
\begin{eqnarray}
\Pi_\fn: {\det}_R(C_{K,S_\fn}(\cA)) &\stackrel{\Pi_\psi}{\to}& {\bigcap}_R^{r+\nu(\fn)}H_\Sigma^1(\co_{K,S_\fn},\cA^\ast(1)) \otimes_R {\bigwedge}_R^{r+\nu(\fn)} (X_\fn^\ast) \nonumber\\
&\stackrel{\sim}{\to}& {\bigcap}_R^{r+\nu(\fn)}H_\Sigma^1(\co_{K,S_\fn},\cA^\ast(1)), \nonumber
\end{eqnarray}
where the first map is well-defined by Lemma \ref{lemma standard}(iv), and the last isomorphism is induced by
$${\bigwedge}_R^{r+\nu(\fn)}(X_\fn^\ast) \stackrel{\sim}{\to} R; \ \gamma_{\fq_1}^\ast \wedge\cdots \wedge \gamma_{\fq_{\nu(\fn)}}^\ast \wedge \beta_1\wedge \cdots \wedge \beta_r\mapsto 1. $$
By Proposition \ref{prop admissible}(v), $\epsilon(z)_\fn$ does not depend on the choice of the standard representative.
One checks that $\epsilon(z):=\{\epsilon(z)_\fn\}_\fn \in \varprojlim_\fn {\bigcap}_R^{r+\nu(\fn)}H_\Sigma^1(\co_{K,S_\fn},\cA^\ast(1))={\rm SS}_r(\cA)$. The homomorphism (\ref{horizontal unit}) is defined by $z \mapsto \epsilon(z)$.
\end{proof}

\begin{definition}
We define the module of {\it basic Stark systems} by
$$\mathcal{S}^{\rm b}(\cA):=\im({\rm HS}(\cA) \to {\rm SS}_r(\cA)).$$
\end{definition}

Similarly to Theorem \ref{det euler fitt}(ii), we deduce the following theorem from Lemma \ref{lemma standard}(iii).

\begin{theorem} \label{det unit fitt}
For every $\fn \in \cN$, we have
$${\rm Fitt}_R^{r+\nu(\fn)}(H^1(C_{K,S_\fn}(\cA)))=\langle \im \epsilon_\fn \mid \epsilon\in \mathcal{S}^{\rm b}(\cA) \rangle_R.$$
\end{theorem}

\subsection{The module of equivariant Stark systems}\label{mess}

In this subsection, we investigate the structures of ${\rm SS}_r(\cA)$ and $\mathcal{S}^{\rm b}(\cA)$. In particular, we give sufficient conditions for the module ${\rm SS}_r(\cA)$ to be  free of rank one (see Theorems \ref{structure unit} and \ref{structure unit 2}).

Let $K(\cA)$ be the minimal Galois extension of $K$ such that $G_{K(\cA)}$ acts trivially on $\cA$. Set $E(\cA)_M:=K(\cA)E_M$.
The following hypothesis is standard in the theory of Euler systems and Kolyvagin systems.

\begin{hypothesis} \label{hyp2}\
\begin{itemize}
\item[(i)] $T\otimes_\cR (\cR/{\rm rad}(\cR))$ is an irreducible $(\cR/{\rm rad}(\cR))[G_K]$-module;
\item[(ii)] there exists $\tau \in G_{E_M}$ such that $\cA/(\tau-1)\cA \simeq R$ as $R$-modules;
\item[(iii)] $H^1(E(\cA)_M/K,\cA)=0.$
\end{itemize}
\end{hypothesis}

\begin{remark}
The hypotheses (i), (ii) and (iii) correspond to (H.1), (H.2) and (H.3) in \cite[\S 3.5]{MRkoly} respectively.
\end{remark}

We also consider the following hypothesis on $\Sigma$.

\begin{hypothesis} \label{hyp surj}
The map $H^1(\co_{K,S}, \cA^\ast(1)) \to \bigoplus_{\fq \in \Sigma}H^1_f(K_\fq, \cA^\ast(1))$ is surjective.
\end{hypothesis}

We set some notations. For any finite set $U$ of places of $K$ such that $S \subset U$, we put
$$\sha^i(\co_{K,U},\cdot):=\ker\left(H^i(\co_{K,U}, \cdot) \to \bigoplus_{v\in U}H^i(K_v, \cdot) \right).$$
Note that the duality
$$\sha^i(\co_{K,U},\cA) \simeq \sha^{3-i}(\co_{K,U},\cA^\ast(1))^\ast$$
holds for $i=1,2$ (see \cite[Chap. I, Theorem 4.10(a)]{milne}).

We set
\begin{eqnarray}
\mathcal{X}_{K,S}(\cA) &:=&  \im \left( H^2(\co_{K,S}, \cA^\ast(1)) \to \bigoplus_{v\in S}H^2(K_v, \cA^\ast(1))  \right)  \nonumber \\
&\simeq& H^2(\co_{K,S},\cA^\ast(1))/\sha^2(\co_{K,S},\cA^\ast(1)). \nonumber
\end{eqnarray}
By the canonical exact sequence
$$0 \to \sha^2(\co_{K,S},\cA^\ast(1)) \to H^2(\co_{K,S},\cA^\ast(1)) \to \bigoplus_{v \in S}H^2(K_v,\cA^\ast(1)) \to H^0(K,\cA)^\ast \to 0$$
(see \cite[Chap. I, Theorem 4.10]{milne}), we see that
\begin{eqnarray}\label{X cokernel}
\mathcal{X}_{K,S}(\cA)^\ast \simeq \coker\left( H^0(K,\cA) \hookrightarrow \bigoplus_{v \in S \setminus S_\infty} H^0(K_v,\cA) \right).
\end{eqnarray}
Here note that $H^2(K_v,\cA^\ast(1))=0$ for any archimedean place $v$, since $p$ is odd.

We also consider the following hypothesis.

\begin{hypothesis} \label{hyp3} $\mathcal{X}_{K,S}(\cA)=0.$
\end{hypothesis}

\begin{remark}
By (\ref{X cokernel}), this hypothesis is satisfied when $H^0(K_v, \cA)=0$ for all $v \in S\setminus S_\infty$, for example.
\end{remark}

\begin{theorem} \label{structure unit}
Assume that Hypotheses \ref{hyp2} and \ref{hyp surj} are satisfied.
\begin{itemize}
\item[(i)] We have an isomorphism of $R$-modules
$$\mathcal{S}^{\rm b}(\cA) \simeq R/{\rm Ann}_R({\rm Fitt}_R^0(\mathcal{X}_{K,S}(\cA))).$$
\item[(ii)] If Hypothesis \ref{hyp3} is also satisfied, then the homomorphism
$${\rm HS}(\cA) \stackrel{ (\ref{horizontal unit})}{\to} {\rm SS}_r(\cA) $$
is an isomorphism of $R$-modules. In particular, we have
$${\rm SS}_r(\cA)=\mathcal{S}^{\rm b}(\cA) \simeq R. $$
\end{itemize}
\end{theorem}

The following terminology will be useful.

\begin{definition}
$\fn\in \cN$ will be said to be {\it large} if $\sha^2(\co_{K,S_\fn},\cA^\ast(1))$ vanishes.
\end{definition}

\begin{proof}[Proof of Theorem \ref{structure unit}]
We first remark that, if Hypothesis \ref{hyp surj} is satisfied, then we have
$$H^2_\Sigma(\co_{K,U},\cA^\ast(1)) =H^2(\co_{K,U},\cA^\ast(1)) $$
for any finite set $U$ of places of $K$ which contains $S$ and is disjoint from $\Sigma$.

We first show (i). Suppose that $\fn \in \cN$ is large. We claim that
\begin{eqnarray}
&&\ker \left(H^1(C_{K,S_\fn}(\cA)) \to Y_K(\cA)^\ast \oplus \bigoplus_{\fq \mid \fn}H^2(K_\fq,\cA^\ast(1)) \right) \nonumber \\
&\simeq& \ker \left(H^2(\co_{K,S_\fn},\cA^\ast(1)) \to \bigoplus_{\fq\mid \fn}H^2(K_\fq,\cA^\ast(1)) \right) \nonumber \\
&\simeq&\mathcal{X}_{K,S}(\cA). \nonumber
\end{eqnarray}
The first isomorphism follows by noting that
\begin{eqnarray}\label{H1decomposition}
H^1(C_{K,S_\fn}(\cA))\simeq Y_K(\cA)^\ast \oplus H^2(\co_{K,S_\fn},\cA^\ast(1)).
\end{eqnarray}
The second isomorphism follows from the exact sequences
\begin{eqnarray}
\bigoplus_{\fq \mid \fn} H_{/f}^1(K_\fq, \cA^\ast(1)) \to \sha^2(\co_{K,S},\cA^\ast(1)) \to \sha^2(\co_{K,S_\fn},\cA^\ast(1)) \to 0, \label{exact 1}
\end{eqnarray}
which is obtained by taking the dual of
$$0 \to \sha^1(\co_{K,S_\fn},\cA) \to \sha^1(\co_{K,S},\cA) \to \bigoplus_{\fq \mid \fn}H^1_f(K_\fq,\cA),$$
and
\begin{equation}\label{H2exact}
\bigoplus_{\fq \mid \fn} H_{/f}^1(K_\fq, \cA^\ast(1)) \to H^2(\co_{K,S},\cA^\ast(1)) \to H^2(\co_{K,S_\fn},\cA^\ast(1)) \to \bigoplus_{\fq \mid \fn}H^2(K_\fq ,\cA^\ast(1)) \to 0.
\end{equation}

By the above claim and Lemma \ref{lemma standard}(ii), we have
$${\rm Ann}_R(\im \Pi_\fn)={\rm Ann}_R({\rm Fitt}_R^0(\mathcal{X}_{K,S}(\cA))) ,$$
where
$$\Pi_\fn: {\det}_R(C_{K,S_\fn}(\cA)) \to {\bigcap}_R^{r+\nu(\fn)} H_\Sigma^1(\co_{K,S_\fn},\cA^\ast(1))$$
is the map defined in the proof of Theorem \ref{theorem HU}. The assertion follows from this and Lemma \ref{key lemma 2} below.

Next, we show (ii). If $\mathcal{X}_{K,S}(\cA)=0$, then the map
$$H^1(C_{K,S_\fn}(\cA)) \to Y_K(\cA)^\ast \oplus \bigoplus_{\fq \mid \fn}H^2(K_\fq,\cA^\ast(1))$$
is an isomorphism for large $\fn$. Hence, the assertion is a consequence of Lemma \ref{lemma standard}(i) and Lemma \ref{key lemma 2} below.
\end{proof}

\begin{lemma} \label{key lemma 2}
Assume that Hypothesis \ref{hyp2} is satisfied. Then, for any $\fn \in \cN$, there exists large $\fm \in \cN$ such that $\fn \mid \fm$.
\end{lemma}

\begin{proof}
The argument is the standard application of Tchebotarev's density theorem, as used in the theory of Euler systems (see \cite{R}).

By the exact sequence (\ref{exact 1}), it is sufficient to show the existence of $\fm$ such that $\fn \mid \fm$ and that the map
$$\bigoplus_{\fq \mid \fm}H^1_{/f}(K_\fq,\cA^\ast(1)) \to \sha^2(\co_{K,S},\cA^\ast(1))\simeq \sha^1(\co_{K,S},\cA)^\ast$$
is surjective. We show that, for any $a \in \sha^1(\co_{K,S},\cA)^\ast$, there exist $\fq \in \cP$ and $e \in H^1_{/f}(K_\fq,\cA^\ast(1))$ such that the map
\begin{eqnarray}
H^1_{/f}(K_\fq ,\cA^\ast(1)) \to \sha^1(\co_{K,S},\cA)^\ast \label{dual local}
\end{eqnarray}
sends $e$ to $a$.

Define the `evaluation homomorphism'
$${\rm Ev}: G_{E(\cA)_M} \to \Hom_R(\sha^1(\co_{K,S},\cA) , \cA/(\tau-1)\cA)$$
by
$${\rm Ev}(\sigma)(b):= b(\sigma),$$
where $\sigma \in G_{E(\cA)_M}$, $b \in \sha^1(\co_{K,S},\cA)$, and $b(\sigma)$ is defined by regarding $b$ as a $1$-cocycle. We show that ${\rm Ev}$ is surjective. Note that $\cA/(\tau-1)\cA\simeq R$ by Hypothesis \ref{hyp2}(ii). We identify these modules. In particular, we identify $\Hom_R(\sha^1(\co_{K,S},\cA) , \cA/(\tau-1)\cA)$ with $\sha^1(\co_{K,S},\cA)^\ast$. Consider the dual of the evaluation homomorphism:
$${\rm Ev}^\ast: \sha^1(\co_{K,S},\cA) \to \Hom(G_{ E(\cA)_M}, R) = \Hom(G_{ E(\cA)_M}, \cA/(\tau-1)\cA) . $$
This homomorphism coincides with the composition map
$$\sha^1(\co_{K,S},\cA) \stackrel{{\rm Res}_{ E(\cA)_M/K}}{\to} \Hom(G_{ E(\cA)_M},\cA)^{G_K} \to \Hom(G_{ E(\cA)_M}, \cA/(\tau-1)\cA),$$
where the second map is induced by the natural map $\cA \to \cA/(\tau-1)\cA$. The first map is injective by Hypothesis \ref{hyp2}(iii). The kernel of the second map is
$$\Hom(G_{ E(\cA)_M},\cA)^{G_K} \cap \Hom(G_{ E(\cA)_M},(\tau-1)\cA).$$
If this set has a non-zero element $f$, then $f(G_{ E(\cA)_M})$ is a non-trivial $G_K$-stable submodule of $(\tau-1)\cA$. This contradicts the irreducibility of $T\otimes_\cR (\cR/{\rm rad}(\cR))$ (see Hypothesis \ref{hyp2}(i)). Hence, we have proved that ${\rm Ev}^\ast$ is injective, and therefore ${\rm Ev}$ is surjective. Note that, letting $L$ be the field corresponding to $\ker {\rm Ev}$, we have the isomorphism (of abelian groups)
$$\Gal(L/ E(\cA)_M) \stackrel{\sim}{\to} \sha^1(\co_{K,S},\cA)^\ast $$
induced by ${\rm Ev}$.

Note that ${\rm Ev}$ can be extended naturally to the homomorphism from $G_{ E(\cA)_M^{\tau=1}}$, where $ E(\cA)_M^{\tau=1}$ is the subfield of $ E(\cA)_M$ fixed by $\tau$. For any $a \in \sha^1(\co_{K,S},\cA)^\ast$, the surjectivity of ${\rm Ev}$ implies the existence of $\alpha \in G_{ E(\cA)_M}$ with ${\rm Ev}(\tau \alpha)=a. $

Let $\mathscr{L}$ be the Galois closure of $L$ over $K$. By Tchebotarev's density theorem, there exists a prime $\fq \notin S\cup \Sigma$ of $K$ such that
$${\rm Fr}_\fq=\tau\alpha|_\mathscr{L} \in \Gal(\mathscr{L}/K)$$
(for some choice of place above $\fq$).
Then we see that $\fq \in \cP$. Define $e \in H^1_{/f}(K_\fq,\cA^\ast(1))$ to be the element which corresponds to the the isomorphism $\cA/(\tau-1)\cA\stackrel{\sim}{\to} R$ we fixed above, under the isomorphism
$$H^1_{/f}(K_\fq,\cA^\ast(1)) \stackrel{(\ref{LD/f})}{\simeq} H_f^1(K_\fq,\cA)^\ast \stackrel{(\ref{frobenius isom})}{\simeq} (\cA/({\rm Fr}_\fq-1)\cA)^\ast =\Hom_R(\cA/(\tau-1)\cA,R).$$
Then we see that the map (\ref{dual local}) sends $e$ to $a$. Thus we have proved the proposition.
\end{proof}

\begin{remark}\label{chebotarev remark}
The proof above shows that, under Hypothesis \ref{hyp2}, any submodule of $\sha^2(\co_{K,S},\cA^\ast(1))$ generated by $i$ elements can be realized as the image of the map
$$\bigoplus_{\fq \mid \fn} H^1_{/f}(K_\fq ,\cA^\ast(1)) \to \sha^2(\co_{K,S},\cA^\ast(1))$$
with some $\fn \in \cN$ with $\nu(\fn)=i$.
\end{remark}




In order to prove that the module of Stark systems ${\rm SS}_r(\cA)$ is free of rank one, it is unnecessary to show that the map (\ref{horizontal unit}) is an isomorphism. In fact,
we can prove that ${\rm SS}_r(\cA)$ is free of rank one under a weaker condition than the assumptions in Theorem \ref{structure unit}(ii) (see Theorem \ref{structure unit 2} below). Consider the following hypothesis.

\begin{hypothesis} \label{hyp3'}
For any $\fn \in \cN$, there exists large $\fm \in \cN$ such that $\fn \mid \fm$ and
$$H_\Sigma^1(\co_{K,S_\fm},\cA^\ast(1)) \simeq R^{\oplus r+\nu(\fm)}.$$
\end{hypothesis}

\begin{proposition}\label{hyp prop}
Hypotheses \ref{hyp2}, \ref{hyp surj}, and \ref{hyp3} imply Hypothesis \ref{hyp3'}.
\end{proposition}

\begin{proof}
From Lemma \ref{key lemma 2} and the proof of Theorem \ref{structure unit}(ii), we see that, for any $\fn \in \cN$, there exits large $\fm \in \cN$ such that $\fn \mid \fm$ and $H^1(C_{K,S_\fm}(\cA))$ is free of rank $r+\nu(\fm)$. Since $C_{K,S_\fm}(\cA)$ is represented by the complex of the form $(F \to F)$ with $F$ free of finite rank (see Proposition \ref{prop admissible}(iv)), we see that $H^0(C_{K,S_\fm}(\cA))=H^1_\Sigma(\co_{K,S_\fm},\cA^\ast(1))$ is also free of rank $r+\nu(\fm)$.
\end{proof}


\begin{theorem} \label{structure unit 2}
If Hypotheses \ref{hyp surj} and \ref{hyp3'} are satisfied, then we have an isomorphism
$${\rm SS}_r(\cA) \simeq R.$$
\end{theorem}

\begin{proof}
For any $\fn, \fm \in \cN$ such that $\fn \mid \fm$, we have the exact sequence (`the global duality')
\begin{eqnarray*}
&&0 \to H^1(\co_{K,S_\fn},\cA^\ast(1)) \to H^1(\co_{K,S_\fm},\cA^\ast(1)) \stackrel{\bigoplus_{\fq\mid \fm/\fn} v_\fq}{\to} R^{\oplus \nu(\fm/\fn)} \nonumber \\
&&\ \to \sha^2(\co_{K,S_\fn},\cA^\ast(1)) \to \sha^2(\co_{K,S_\fm},\cA^\ast(1))\to 0.
\end{eqnarray*}
(See \cite[Chap. I, Theorem 4.10]{milne} or \cite[Theorem 1.7.3]{R}.) By Hypothesis \ref{hyp surj}, we have
$$\ker(\sha^2(\co_{K,S_\fn},\cA^\ast(1)) \to \sha^2(\co_{K,S_\fm},\cA^\ast(1))) = \ker (H^2_\Sigma(\co_{K,S_\fn}, \cA^\ast(1)) \to H^2_\Sigma(\co_{K,S_\fm}, \cA^\ast(1))).$$
So we obtain the exact sequence
\begin{eqnarray} \label{GD}
&&0 \to H_\Sigma^1(\co_{K,S_\fn},\cA^\ast(1)) \to H_\Sigma^1(\co_{K,S_\fm},\cA^\ast(1)) \stackrel{\bigoplus_{\fq\mid \fm/\fn} v_\fq}{\to} R^{\oplus \nu(\fm/\fn)} \nonumber \\
&&\ \to \sha^2(\co_{K,S_\fn},\cA^\ast(1)) \to \sha^2(\co_{K,S_\fm},\cA^\ast(1))\to 0.
\end{eqnarray}

Suppose that $\fn$ is large and we have isomorphisms
$$H_\Sigma^1(\co_{K,S_\fn},\cA^\ast(1)) \simeq R^{\oplus r+\nu(\fn)} \text{ and } H_\Sigma^1(\co_{K,S_\fm},\cA^\ast(1)) \simeq R^{\oplus r+\nu(\fm)}.$$
Then by (\ref{GD}) we have a split short exact sequence
$$0 \to H_\Sigma^1(\co_{K,S_\fn},\cA^\ast(1)) \to H_\Sigma^1(\co_{K,S_\fm},\cA^\ast(1)) \stackrel{\bigoplus_{\fq\mid \fm/\fn} v_\fq}{\to} R^{\oplus \nu(\fm/\fn)} \to 0,$$
and we see that ${\bigwedge}_{\fq \mid \fm/\fn} v_\fq$ induces an isomorphism
$${\bigcap}_R^{r+\nu(\fm)}H_\Sigma^1(\co_{K,S_\fm},\cA^\ast(1)) \stackrel{\sim}{\to}  {\bigcap}_R^{r+\nu(\fn)} H_\Sigma^1(\co_{K,S_\fn},\cA^\ast(1)).$$
These modules are both free of rank one. So, by Hypothesis \ref{hyp3'}, we see that ${\rm SS}_r(\cA)$ is free of rank one.
\end{proof}

\begin{remark} \label{remark core vertex}
An ideal $\fm \in \cN$ satisfying the condition in Hypothesis \ref{hyp3'} corresponds to a `core vertex' in the theory of Mazur and Rubin in \cite{MRkoly} and \cite{MRselmer}. To be more precise, when $R$ is principal, one can consider the notion of `core rank' (see \cite[Definition 4.1.11]{MRkoly}), and if the core rank for a suitable Selmer structure $\mathcal{F}$ is equal to $r$, then any core vertex $\fm \in \cN$ satisfies $H^1_{\mathcal{F}^{\fm}}(K,\cA^\ast(1)) \simeq R^{\oplus r+\nu(\fm)}$ (see \cite[Theorem B.2]{MRkoly}). Also, Mazur and Rubin proved that, for any $\fn \in \cN$, there exists a core vertex $\fm \in \cN$ such that $\fn \mid \fm$ (under their running hypotheses, see \cite[Corollary 4.1.9(iii)]{MRkoly}). So in their theory the condition corresponding to Hypothesis \ref{hyp3'} is ensured by their running hypotheses. We remark that Theorem \ref{structure unit 2} is a direct generalization of \cite[Theorem 6.10]{MRselmer}.
\end{remark}

\subsection{Higher Fitting ideals of Selmer groups}


For each commutative noetherian ring $R$, finitely generated $R$-module $X$, $R$-submodule $Y$ of $X$ and pair of non-negative integers $a$ and $b$ we now use the notion of `higher relative Fitting ideal' ${\rm Fitt}_R^{(a,b)}(X,Y)$ introduced by Kurihara and the present authors in \cite[\S7.1]{bks1}. We recall, in particular, that ${\rm Fitt}_R^{(a,b)}(X,Y)$ is contained in the classical higher Fitting ideal ${\rm Fitt}_R^{a+b}(X,Y)$, with equality if $X=Y$ (for a more precise result see \cite[Lemma 7.2]{bks1}).

In particular, in this subsection, we prove the following theorem.

\begin{theorem}\label{theorem higher fitt}
Assume Hypotheses \ref{hyp2} and \ref{hyp surj}.
\begin{itemize}
\item[(i)] Let $\epsilon \in \mathcal{S}^{\rm b}(\cA)$ be a generator. Then, for every non-negative integer $i$, we have
$${\rm Fitt}_R^{(0,i)}(H^2(\co_{K,S}, \cA^\ast(1)), \sha^2(\co_{K,S},\cA^\ast(1))) = \langle \im \epsilon_\fn \mid \fn \in \cN \text{ with }\nu(\fn)=i \rangle_R. $$
\item[(ii)] Assume also Hypothesis \ref{hyp3'} (so ${\rm SS}_r(\cA) \simeq R$ by Theorem \ref{structure unit 2}).
Let $\epsilon=\{\epsilon_\fn\}_\fn \in {\rm SS}_r(\cA)$ be a basis. Then, for every non-negative integer $i$, we have
$${\rm Fitt}_R^i(\sha^2(\co_{K,S},\cA^\ast(1)))={\rm Fitt}_R^i(\sha^1(\co_{K,S},\cA)^\ast)=\langle \im \epsilon_\fn \mid \fn \in \cN \text{ with }\nu(\fn)=i \rangle_R.$$
\end{itemize}
\end{theorem}

\begin{proof}
We first note that, for every $\fn \in \cN$, we have
\begin{eqnarray}
{\rm Fitt}_R^{r+\nu(\fn)}(H^1(C_{K,S_\fn}(\cA))) &=& {\rm Fitt}_R^{\nu(\fn)}(H^2(\co_{K,S_\fn}, \cA^\ast(1))) \nonumber \\
&=&{\rm Fitt}_R^0\left(\coker \left(\bigoplus_{\fq \mid \fn} H^1_{/f}(K_\fq,\cA^\ast(1)) \to H^2(\co_{K,S},\cA^\ast(1)) \right)\right). \nonumber
\end{eqnarray}
In fact, the first equality follows from (\ref{H1decomposition}), and the second from (\ref{H2exact}). Hence, by Theorem \ref{det unit fitt}, we have
$${\rm Fitt}_R^0\left(\coker \left(\bigoplus_{\fq \mid \fn} H^1_{/f}(K_\fq,\cA^\ast(1)) \to H^2(\co_{K,S},\cA^\ast(1)) \right)\right) = \im \epsilon_\fn. $$

Now note that, for a finitely generated $R$-module $X$ and its submodule $Y$, we have
$${\rm Fitt}_R^{(0,i)}(X,Y)=\langle {\rm Fitt}_R^0(X/Z) \mid Z \subset Y \text{ is generated by $i$ elements}\rangle_R.$$
(See \cite[Lemma 7.1]{bks1}.) Claim (i) follows from this and Remark \ref{chebotarev remark}.

Next, we show claim (ii). Take $\fn \in \cN$ with $\nu(\fn)=i$. By Hypothesis \ref{hyp3'}, we can take large $\fm\in \cN$ such that $\fn \mid \fm$ and  $H_\Sigma^1(\co_{K,S_\fm},\cA^\ast(1))$ is free of rank $r+\nu(\fm)$. For such $\fm$, we see that $\epsilon_\fm$ is a basis of ${\bigcap}_R^{r+\nu(\fm)}H_\Sigma^1(\co_{K,S_\fm},\cA^\ast(1))$. Note that $\epsilon_\fn=\pm ({\bigwedge}_{\fq \mid \fm/\fn}v_\fq) (\epsilon_\fm)$. By applying Proposition \ref{prop injective}(ii) to the exact sequence (\ref{GD}), we have
\begin{eqnarray} \label{fitt0 epsilon}
 {\rm Fitt}_R^0(\sha^2(\co_{K,S_\fn},\cA^\ast(1)))=\im \epsilon_\fn.
 \end{eqnarray}
Hence, it is sufficient to show
$${\rm Fitt}_R^i(\sha^2(\co_{K,S},\cA^\ast(1)))=\langle {\rm Fitt}_R^0(\sha^2(\co_{K,S_\fn},\cA^\ast(1))) \mid \fn \in \cN \text{ with }\nu(\fn)=i \rangle_R.$$
For a finitely generated $R$-module $X$, we have
$${\rm Fitt}_R^{i}(X)={\rm Fitt}_R^{(0,i)}(X,X)=\langle {\rm Fitt}_R^0(X/Z) \mid Z \subset X \text{ is generated by $i$ elements}\rangle_R.$$
Noting this, the assertion follows from the exact sequence (\ref{exact 1}) and Remark \ref{chebotarev remark}.
\end{proof}

\begin{remark}
Assume Hypotheses \ref{hyp2}, \ref{hyp surj}, and \ref{hyp3}. Then we see that Hypothesis \ref{hyp3'} is also satisfied by Proposition \ref{hyp prop}, and we have $H^2(\co_{K,S},\cA^\ast(1))=\sha^2(\co_{K,S},\cA^\ast(1))$ by Hypothesis \ref{hyp3} and $\mathcal{S}^{\rm b}(\cA) ={\rm SS}_r(\cA)$ by Theorem \ref{structure unit}(ii). One sees that claims (i) and (ii) in Theorem \ref{theorem higher fitt} are identical in this case. In fact, we have
\begin{eqnarray}
&&{\rm Fitt}_R^{(0,i)}(H^2(\co_{K,S}, \cA^\ast(1)), \sha^2(\co_{K,S}, \cA^\ast(1))) \nonumber \\
&=& {\rm Fitt}_R^{(0,i)}(\sha^2(\co_{K,S}, \cA^\ast(1)), \sha^2(\co_{K,S}, \cA^\ast(1))) \nonumber \\
&=& {\rm Fitt}_R^i(\sha^2(\co_{K,S},\cA^\ast(1))) . \nonumber
\end{eqnarray}
\end{remark}

\begin{remark}
Theorem \ref{theorem higher fitt}(ii) constitutes a natural generalization of \cite[Proposition 8.5]{MRselmer} to more general coefficient rings. (Our $\sha^1(\co_{K,S},\cdot)$ corresponds to $H^1_{\cF^\ast}(K, \cdot)$ in loc. cit.) Indeed, the coefficient ring is assumed to be principal in loc. cit., whereas our ring $R=\fr[\Gamma]$ is not principal in general.
\end{remark}

\section{Equivariant Kolyvagin systems}\label{eks}

In this section, we study higher rank Kolyvagin systems with equivariant coefficients. Unlike Mazur and Rubin \cite{MRselmer}, we always use the exterior bidual instead of the exterior power, as in the previous section. In \S \ref{def kolysys}, we give a definition of higher rank Kolyvagin systems (see Definition \ref{def koly}). In \S \ref{unit kolysys}, we review the construction of Kolyvagin systems from Stark systems (see Proposition \ref{unit koly}). In \S \ref{euler kolysys} we use the theory of vertical and horizontal determinantal systems to give a canonical, and very natural, construction of Kolyvagin systems from basic Euler systems (see Theorem \ref{main}).

\subsection{The definition} \label{def kolysys}

Although it is possible to define Kolyvagin systems over a general self-injective ring and for a general Selmer structure as in \S \ref{section def unit}, we restrict the setting to the arithmetic one as in \S \ref{section horizontal} (which is, of course, an interesting and important case). In particular, we keep the notations used in \S \ref{section horizontal}. We remark that we do not need to assume Hypotheses \ref{hyp free} and \ref{hyp1} for the moment.

For a prime $\fq$ of $K$, we denote by $K(\fq)$ the maximal $p$-extension of $K$ inside the ray class field modulo $\fq$. For a square-free product $\fn$ of primes of $K$, we denote by $K(\fn)$ the composite of $K(\fq)$'s with $\fq \mid \fn$.

For $\fn \in \cN$, define
$$G_\fn:=\bigotimes_{\fq\mid \fn}\Gal(K(\fq)_\fq/K_\fq).$$
Note that, since $\fq \in \cP$ splits completely in $K(\mu_{\overline M},(\co_K^\times)^{1/\overline M})K(1)$, the order of $G_\fq=\Gal(K(\fq)_\fq/K_\fq)$ is divisible by $ M$ in $O$ (see \cite[Lemma 4.1.2(i)]{R}). In particular, we have
$$R \otimes_\ZZ G_\fn \simeq R.$$
In what follows, we often write $\otimes_\ZZ$ simply as $\otimes$.

For $\fq\in \cP$, we recall the definition of the `finite-singular comparison map'
$$\varphi_\fq^{\rm fs}: H^1_f(K_\fq,\cA^\ast(1)) \stackrel{\sim}{\to} H^1_{/f}(K_\fq,\cA^\ast(1))\otimes G_\fq$$
(see \cite[Lemma 1.2.3]{MRkoly}).
Note first that we have a canonical isomorphism
\begin{eqnarray}
H^1_f(K_\fq,\cA^\ast(1)) \simeq \cA^\ast(1)/({\rm Fr}_\fq-1)\cA^\ast(1); \ a \mapsto a({\rm Fr}_\fq). \label{fs 1}
\end{eqnarray}
Next, we have a natural identification
$$H^1_{/f}(K_\fq,\cA^\ast(1))\simeq \Hom(G_\fq,\cA^\ast(1)^{{\rm Fr}_\fq=1}),$$
so we have a canonical isomorphism
\begin{eqnarray}\label{fs 2}
H^1_{/f}(K_\fq,\cA^\ast(1))\otimes G_\fq \stackrel{\sim}{\to} \cA^\ast(1)^{{\rm Fr}_\fq=1}; \ f\otimes \sigma \mapsto f(\sigma).
\end{eqnarray}
The finite-singular comparison map is defined by composing (\ref{fs 1}), the isomorphism
$$Q_\fq({\rm Fr_\fq^{-1}}): \cA^\ast(1)/({\rm Fr}_\fq-1)\cA^\ast(1)  \stackrel{\sim}{\to} \cA^\ast(1)^{{\rm Fr}_\fq=1},$$
and the inverse of (\ref{fs 2}).

Define the `transverse subgroup' by
$$H^1_{\rm tr}(K_\fq, \cA^\ast(1)):=H^1(G_\fq, \cA^\ast(1)^{G_{K(\fq)_\fq}}) \subset H^1(K_\fq, \cA^\ast(1)).$$
Note that
$$H^1_{\rm tr}(K_\fq, \cA^\ast(1)) \stackrel{\sim}{\to} H^1_{/f}(K_\fq,\cA^\ast(1))$$
and
$$H^1(K_\fq,\cA^\ast(1)) = H_f^1(K_\fq,\cA^\ast(1)) \oplus H_{\rm tr}^1(K_\fq,\cA^\ast(1)).$$
(See \cite[Lemma 1.2.4]{MRkoly}.) In particular, we have a canonical projection
$$H^1(K_\fq,\cA^\ast(1)) \to H^1_f(K_\fq,\cA^\ast(1)) .$$

Following \cite[\S 2.1]{MRkoly}, we use the following notation. For $\fa, \fb, \fc \in \cN$, set
$$H^1_{\cF_\fa^\fb(\fc)}(K, \cdot):=\ker \left( H_\Sigma^1(\co_{K,S_{\fa\fb\fc}}, \cdot) \to \left( \bigoplus_{\fq \mid \fa} H^1(K_\fq, \cdot) \right) \oplus\left( \bigoplus_{\fq \mid \fc} H^1_f(K_\fq, \cdot)\right) \right).$$
If any of $\fa, \fb,\fc$ are equal to $1$, we omit them from the notation. For example, we have
$$H^1_\cF(K,\cA^\ast(1))=H_\Sigma^1(\co_{K,S},\cA^\ast(1)),$$
$$H^1_{\cF^{\fn}}(K,\cA^\ast(1))=H_\Sigma^1(\co_{K,S_\fn},\cA^\ast(1)),$$
$$H^1_{\cF(\fn)}(K,\cA^\ast(1))=\ker \left( H^1_{\cF^\fn}(K,\cA^\ast(1)) \to  \bigoplus_{\fq \mid \fn} H^1_f(K_\fq, \cA^\ast(1)) \right).$$





We often denote the composition map
$$H^1(K,\cA^\ast(1)) \to H^1(K_\fq,\cA^\ast(1)) \to H^1_f(K_\fq,\cA^\ast(1)) \stackrel{\varphi_\fq^{\rm fs}}{\to} H^1_{/f}(K_\fq,\cA^\ast(1))\otimes G_\fq \stackrel{v_\fq }{\to}  R \otimes G_\fq,$$
also by $\varphi_\fq^{\rm fs}$. Let $r$ be a positive integer.
Note that, if $\fq \mid \fn$, then by Proposition \ref{prop self injective} we see that $v_\fq$ and $\varphi_\fq^{\rm fs}$ induce
$$v_\fq: {\bigcap}_R^{r}H^1_{\cF(\fn)}(K,\cA^\ast(1)) \otimes G_\fn \to {\bigcap}_{R}^{r-1} H^1_{\cF_\fq(\fn)}(K,\cA^\ast(1))\otimes G_\fn,$$
$$\varphi_\fq^{\rm fs}: {\bigcap}_R^{r}H^1_{\cF(\fn/\fq)}(K,\cA^\ast(1))\otimes G_{\fn/\fq} \to {\bigcap}_{R}^{r-1} H^1_{\cF_\fq(\fn)}(K,\cA^\ast(1))\otimes G_\fn$$
respectively.

\begin{definition}\label{def koly}
A {\it Kolyvagin system of rank $r$} for $\cA^\ast(1)$ is a collection
$$\left\{ \kappa_\fn \in {\bigcap}_{R}^r H^1_{\cF(\fn)}(K,\cA^\ast(1))\otimes G_\fn \ \middle| \ \fn\in \cN \right\}$$
which satisfies
$$v_\fq(\kappa_\fn)=\varphi_\fq^{\rm fs}(\kappa_{\fn/\fq}) \text{ in }{\bigcap}_{R}^{r-1} H^1_{\cF_\fq(\fn)}(K,\cA^\ast(1))\otimes G_\fn$$
for any $\fn \in \cN$ and $\fq \mid \fn$.

The module of Kolyvagin systems of rank $r$ for $\cA^\ast(1)$ is denoted by ${\rm KS}_r(\cA)$.
\end{definition}

\begin{remark} \label{difference kolyvagin}
Note that, when $r>1$, our definition of Kolyvagin systems is different from that by Mazur and Rubin \cite{MRselmer} (even when $R$ is principal), since we use the exterior bidual ${\bigcap}_R^r$ instead of the usual exterior power ${\bigwedge}_R^r$. When $r=1$, we can identify ${\bigcap}_R^1 X = X^{\ast \ast}=X$ with $X$ for any $R$-module $X$, so in this case our definition is the same as that by Mazur and Rubin \cite{MRkoly}.
\end{remark}

\subsection{Kolyvagin systems from Stark systems} \label{unit kolysys}

We review the construction of Kolyvagin systems from Stark systems (cf. \cite[\S 5]{sanojnt} and \cite[\S 12]{MRselmer}).


Suppose that we have a Stark system
$$\epsilon=\left\{ \epsilon_\fn \in {\bigcap}_{R}^{r + \nu(\fn)} H^1_{\cF^\fn} (K, \cA^\ast(1)) \ \middle| \ \fn \in \cN \right\} \in {\rm SS}_r(\cA).$$
Note that, by Proposition \ref{prop self injective}, ${\bigwedge}_{\fq\mid \fn}\varphi_\fq^{\rm fs}$ induces a map
$${\bigwedge}_{\fq \mid \fn}\varphi_\fq^{\rm fs} : {\bigcap}_{R}^{r + \nu(\fn)} H^1_{\cF^\fn} (K, \cA^\ast(1)) \to {\bigcap}_{R}^r H^1_{\cF(\fn)}(K,\cA^\ast(1))\otimes G_\fn.$$
We define
$$\kappa(\epsilon)_\fn:=(-1)^{\nu(\fn)}\left({\bigwedge}_{\fq \mid \fn} \varphi_\fq^{\rm fs} \right)(\epsilon_\fn) \in  {\bigcap}_{R}^r H^1_{\cF(\fn)}(K,\cA^\ast(1))\otimes G_\fn.$$

\begin{proposition}[{\cite[Theorem 5.7]{sanojnt}, \cite[Proposition 12.3]{MRselmer}}] \label{unit koly}
$\kappa(\epsilon):=\{ \kappa(\epsilon)_\fn\}_\fn$ is a Kolyvagin system.
\end{proposition}

\begin{proof}
For $\fn \in \cN$ and $\fq\mid \fn$, we compute
\begin{eqnarray}
v_\fq(\kappa(\epsilon)_\fn)&=&(-1)^{\nu(\fn)}v_\fq\left( \left({\bigwedge}_{\fq' \mid \fn}\varphi_{\fq'}^{\rm fs}\right)(\epsilon_\fn) \right) \nonumber \\
&=&  {\rm sgn}(\fn,\fn/\fq)\left({\bigwedge}_{\fq' \mid \fn}\varphi_{\fq'}^{\rm fs} \right)(\epsilon_{\fn/\fq}) \nonumber \\
 &=&(-1)^{\nu(\fn/\fq)}\varphi_\fq^{\rm fs}\left(\left({\bigwedge}_{\fq'\mid \fn/\fq}\varphi_{\fq'}^{\rm fs}\right)(\epsilon_{\fn/\fq})\right) \nonumber \\
&=&\varphi_\fq^{\rm fs}(\kappa(\epsilon)_{\fn/\fq}). \nonumber
\end{eqnarray}
\end{proof}

Proposition \ref{unit koly} gives the homomorphism
\begin{eqnarray*} \label{US to KS}
{\rm Reg}_r: {\rm SS}_r(\cA) \to {\rm KS}_r(\cA); \ \epsilon \mapsto \kappa(\epsilon).
\end{eqnarray*}
We think of this map as a `regulator' (as in \cite[\S 5]{sanojnt}), and give the following definition.


\begin{definition} \label{def stub}
We define the module of {\it regulator Kolyvagin systems} by
$${\rm KS}_r^{\rm reg}(\cA):=\im({\rm Reg}_r) \subset {\rm KS}_r(\cA).$$
\end{definition}

\begin{remark}
In the non-equivariant setting, Mazur and Rubin proved in \cite[Theorem 12.4]{MRselmer} that, under their working hypotheses, the image of the module of Stark systems under the map ${\rm Reg}_r$ is the submodule of `stub-Kolyvagin systems' (though our definition of Kolyvagin systems and Stark systems are slightly different from theirs, see Remark \ref{difference kolyvagin}). So the module of regulator Kolyvagin systems can be regarded as an interpretation of stub-Kolyvagin systems by Mazur and Rubin.
\end{remark}

\begin{theorem} \label{bounding selmer}
Assume that Hypotheses \ref{hyp surj} and \ref{hyp3'} are satisfied. Then we have
$${\rm Fitt}_R^0(\sha^2(\co_{K,S}, \cA^\ast(1))) =\langle \im \kappa_1 \mid \kappa \in {\rm KS}_r^{\rm reg}(\cA) \rangle_R . $$
\end{theorem}

\begin{proof}
Since any regulator Kolyvagin system $\kappa$ is written as $\kappa=\kappa(\epsilon)$ with some Stark system $\epsilon \in {\rm SS}_r(\cA)$, and $\kappa(\epsilon)_1=\epsilon_1$ by construction, the theorem follows from (\ref{fitt0 epsilon}). (Note that Hypothesis \ref{hyp2} is not used.)
\end{proof}

\begin{remark}
Theorem \ref{bounding selmer} can be regarded as an equivariant version of \cite[Theorem 4.5.6]{MRkoly} and \cite[Corollary 13.1]{MRselmer}.
\end{remark}

\subsection{Kolyvagin systems from Euler systems} \label{euler kolysys}


In this subsection, we assume the following
\begin{hypothesis}\label{hyp last}\
\begin{itemize}
\item[(i)] Hypothesis \ref{hyp free} and $r:=r_T={\rm rank}_\cR(Y_K(T))\geq 1$;
\item[(ii)] Hypothesis \ref{hyp1};
\item[(iii)]$\cK$ contains $E(\fn):=E\cdot K(\fn)$ for all $\fn \in \cN$;
\item[(iv)]$E$ contains $K(1)$, i.e. $E(1)=E$.
\end{itemize}
\end{hypothesis}


The aim of this section is to give a construction of Kolyvagin systems from Euler systems via the classical construction of Kolyvagin's `derivative classes' (see Theorem \ref{main}).

We set some notations. For $\fn \in \cN$, put $\cG_\fn:=\Gal(E(\fn)/K)$ and $\cH_\fn:=\Gal(E(\fn)/E)$.
By Hypothesis \ref{hyp last}(iv), we have
$$\cH_\fn \simeq \Gal(K(\fn)/K(1))\simeq \prod_{\fq \mid \fn} G_\fq.$$
(Recall that $E/K$ is unramified outside $S$.) Fix a generator $\sigma_\fq \in G_\fq$ for each $\fq \in \cP$. Since each $\fq \in \cP$ splits completely in $E$, we can regard ${\rm Fr}_\fq \in \cH_{\fn}$ when $\fq \not{\mid} \fn$. In particular, if $\fq \neq \fq'$, we can regard ${\rm Fr}_\fq \in G_{\fq'}$. These observations will be used later.

For an Euler system $c=\{c_F\}_F \in {\rm ES}_r(T,\cK)$, we set
$$c_\fn:=c_{E(\fn)} \in {\bigcap}_{\cR[\cG_\fn]}^r H_\Sigma^1(\co_{E(\fn),S_\fn},T^\ast(1)).$$
We often regard $c_\fn \in {\bigcap}_{\fr[\cG_\fn]}^r H_\Sigma^1(\co_{E(\fn),S_\fn},A^\ast(1))$ via the map induced by the natural surjection $T\to A(=T/MT)$.

\subsubsection{Kolyvagin's derivative construction}

Suppose that we have a rank $r$ Euler system $c = \{c_F\}_F \in {\rm ES}_r(T,\cK)$, where $r$ is an arbitrary positive integer (we do not require $r=r_T$ for the moment). For $\fn \in \cN$, Kolyvagin's derivative operator is defined by
$$D_\fn:=\prod_{\fq \mid \fn} \left( \sum_{i=1}^{\# G_\fq -1} i \sigma_\fq^i  \right) \in \ZZ[\cH_\fn].$$
Note that we have a natural isomorphism
$$H_\Sigma^1(\co_{E(\fn),S_\fn},A^\ast(1))^{\cH_\fn} \simeq H_\Sigma^1(\co_{E,S_\fn},A^\ast(1)) =H_\Sigma^1(\co_{K,S_\fn}, \cA^\ast(1)) = H^1_{\cF^\fn}(K,\cA^\ast(1)).$$
(This isomorphism follows from $H_\Sigma^0(E(\fn),A^\ast(1))=0$, which is a consequence of Hypothesis \ref{hyp1}.)
By \cite[Lemma 4.4.2]{R}, we have
$$\kappa'(c)_\fn:= D_\fn  c_\fn \text{ mod $M$} \in \left({\bigcap}_{\fr[\cG_\fn]}^r H_\Sigma^1(\co_{E(\fn),S_\fn},A^\ast(1))\right)^{\cH_\fn} \simeq {\bigcap}_{R}^r H_{\cF^\fn}^1(K,\cA^\ast(1)),$$
where the last isomorphism follows from Proposition \ref{bidual invariant}.

Let $\cI_\fq$ be the augmentation ideal of $\ZZ[G_\fq]$. We have a natural isomorphism
$$\rho_\fq: G_\fq \stackrel{\sim}{\to} \cI_\fq/\cI_\fq^2; \ \sigma \mapsto \sigma-1. $$

We introduce the following definition.
\begin{definition}[{\cite[Definition 4.7]{sanojnt}}] \label{def DKS}
A collection
$$\left\{ \kappa_\fn' \in {\bigcap}_R^r H^1_{\cF^\fn}(K,\cA^\ast(1)) \otimes G_\fn \ \middle| \ \fn \in \cN  \right\}$$
is called a {\it derived Kolyvagin system of rank $r$} if the following conditions are satisfied:
\begin{itemize}
\item[(i)] if $\fn \in \cN$ and $\fq \mid \fn$, then
$$\varphi_\fq^{\rm fs}(\kappa_{\fn/\fq}')=v_\fq(\kappa_\fn') \text{ in }{\bigcap}_R^{r-1} H^1_{\cF^{\fn/\fq}}(K,\cA^\ast(1)) \otimes G_\fn;$$
\item[(ii)] if $\fn \in \cN$, then the element
$$\sum_{\tau \in \mathfrak{S}(\fn)}{\rm sgn}(\tau) \kappa_{\fd_\tau}'\otimes \bigotimes_{\fq \mid \fn/\fd_\tau} \rho_\fq^{-1}(P_{\tau(\fq)}({\rm Fr}_{\tau(\fq)}^{-1}))  \in  {\bigcap}_R^r H^1_{\cF^\fn}(K,\cA^\ast(1))\otimes G_\fn$$
lies in ${\bigcap}_R^r H^1_{\cF(\fn)}(K,\cA^\ast(1))\otimes G_\fn$, where $\mathfrak{S}(\fn)$ is the set of permutations of prime divisors of $\fn$, and $\fd_\tau:=\prod_{\tau(\fq)=\fq}\fq$. ($P_{\tau(\fq)}({\rm Fr}_{\tau(\fq)}^{-1})$ is regarded as an element of $\fr \otimes \cI_\fq/\cI_\fq^2$, so $\rho_\fq^{-1}(P_{\tau(\fq)}({\rm Fr}_{\tau(\fq)}^{-1}))$ lies in $\fr \otimes G_\fq$.)
\end{itemize}
The module of derived Kolyvagin systems of rank $r$ is denoted by ${\rm DKS}_r(\cA)$.
\end{definition}

\begin{remark} \label{remark r1}
For an Euler system $c$, we expect that the collection $\widetilde \kappa'(c):=\{\kappa'(c)_\fn \otimes \bigotimes_{\fq \mid \fn} \sigma_\fq\}_\fn$ lies in ${\rm DKS}_r(\cA)$. When $r=1$, one can show that this is the case under mild conditions. In fact, suppose that $r=1$ and that $E=K(1)=K$. In this case, we can prove that $\widetilde \kappa'(c)$ satisfies Definition \ref{def DKS}(i)
under the assumption that $\cK$ contains maximal abelian $p$-extension of $K$ which is unramified outside $S_p(K) \cup \cP$ (see \cite[Theorem 4.5.4]{R}).
Also, Mazur and Rubin proved the following claim (see \cite[Theorem 3.2.4]{MRkoly}), which implies that $\widetilde \kappa'(c) \in {\rm DKS_1}(\cA)$. Assume the following:
\begin{itemize}
\item[(a)] $\cK$ contains maximal abelian $p$-extension of $K$ which is unramified outside $S_p(K) \cup \cP$;
\item[(b)] ${\rm Fr}_\fq^{p^k} -1$ is injective on $T^\ast(1)$ for all $\fq \in \cP$ and $k \geq0$.
\end{itemize}
For every $\fn \in \cN$, define the element
$$\kappa(c)_\fn:= \sum_{\tau \in \mathfrak{S}(\fn)}{\rm sgn}(\tau) \kappa'(c)_{\fd_\tau} \otimes \left( \bigotimes_{\fq \mid \fd_\tau}\sigma_\fq  \otimes \bigotimes_{\fq \mid \fn/\fd_\tau} \rho_\fq^{-1}(P_{\tau(\fq)}({\rm Fr}_{\tau(\fq)}^{-1})) \right) \in  H^1_{\cF^\fn}(K,\cA^\ast(1))\otimes G_\fn.$$
Then we have
$$\kappa(c)_\fn \in H^1_{\cF(\fn)}(K,\cA^\ast(1))\otimes G_\fn,$$
and $\{ \kappa(c)_\fn\}_\fn \in {\rm KS}_1(\cA)$.
\end{remark}


It is useful to give the following definition.

\begin{definition}\label{def KC}
We define the module of {\it Kolyvagin collections of rank $r$} by
$${\rm KC}_r(\cA):=\prod_{\fn \in \cN} {\bigcap}_R^r H^1_{\cF^\fn}(K,\cA^\ast(1)) \otimes G_\fn.$$
\end{definition}
Obviously, we have
$${\rm DKS}_r(\cA) \subset {\rm KC}_r(\cA) \text{ and } {\rm KS}_r(\cA) \subset {\rm KC}_r(\cA).$$

We shall define a natural operator which shifts derived Kolyvagin systems to Kolyvagin systems.
Define an endomomorphism $\Psi_r \in {\rm End}_R({\rm KC}_r(\cA))$ by
$$\Psi_r(\{\kappa_\fn'\}_\fn):=\left\{ \sum_{\tau \in \mathfrak{S}(\fn)}{\rm sgn}(\tau) \kappa_{\fd_\tau}'\otimes \bigotimes_{\fq \mid \fn/\fd_\tau} \rho_\fq^{-1}(P_{\tau(\fq)}({\rm Fr}_{\tau(\fq)}^{-1})) \right\}_\fn.$$

\begin{theorem}[{\cite[Theorem 4.17]{sanojnt}}] \label{isom dks ks}
The endomorphism $\Psi_r$ induces an isomorphism
$$\Psi_r: {\rm DKS}_r(\cA) \stackrel{\sim}{\to} {\rm KS}_r(\cA).$$
\end{theorem}

\begin{proof}
It is straightforward to check that $\Psi_r({\rm DKS}_r(\cA)) \subset {\rm KS}_r(\cA)$. We construct the inverse of $\Psi_r$. For a given $\kappa=\{\kappa_\fn\}_\fn \in {\rm KS}_r(\cA)$, define $\kappa'=\{\kappa'_\fn\}_\fn$ inductively by
$$\kappa_\fn= \sum_{\tau \in \mathfrak{S}(\fn)}{\rm sgn}(\tau) \kappa_{\fd_\tau}'\otimes \bigotimes_{\fq \mid \fn/\fd_\tau} \rho_\fq^{-1}(P_{\tau(\fq)}({\rm Fr}_{\tau(\fq)}^{-1})) .$$
Then one checks that $\kappa' \in {\rm DKS}_r(\cA)$, and that the map $\kappa \mapsto \kappa'$ is the inverse of $\Psi_r$.
\end{proof}

\begin{definition}\label{derivable}
Let
$$\mathcal{D}_r: {\rm ES}_r(T,\cK) \to {\rm KC}_r(\cA)$$
be the homomorphism defined by
$$\mathcal{D}_r(c):=\Psi_r\left( \left\{\kappa'(c)_\fn \otimes \bigotimes_{\fq \mid \fn}\sigma_\fq \right\}_\fn \right).$$
An Euler system $c \in {\rm ES}_r(T,\cK)$ is said to be {\it Kolyvagin-derivable} (resp. {\it Stark-derivable}) if $\mathcal{D}_r(c)$ lies in ${\rm KS}_r(\cA)$ (resp. ${\rm KS}_r^{\rm reg}(\cA)$).

\end{definition}

\begin{remark}
We expect that all Euler systems are Stark-derivable. As explained in Remark \ref{remark r1}, Mazur and Rubin proved that, under mild conditions, all rank one Euler systems are Kolyvagin-derivable. Furthermore, they also proved that ${\rm KS}_1^{\rm reg}(\cA)={\rm KS}_1(\cA)$ under their working hypotheses (see \cite[Theorem 4.4.1]{MRkoly} and \cite[Theorem 12.4]{MRselmer}), so in the rank one case all Euler systems are Stark-derivable.
\end{remark}

We first give the following result.

\begin{theorem}\label{stark der thm 2}
Assume Hypotheses \ref{hyp2}, \ref{hyp surj}, \ref{hyp3'} and \ref{hyp last}. Fix a system $c$ in ${\rm ES}_r(T,\mathcal{K})$ that is Stark-derivable.
Then, there exists a canonical ideal $J=J_{c,T,E,M}$ of $R$ for which
 $ \langle \mathcal{D}_r( c) \rangle_{R}= J\cdot {\rm KS}_r^{{\rm reg}}(\mathcal{A})$ and the associated module of Stark systems
\[ \{ \epsilon \in {\rm SS}_r(\mathcal{A}) \mid {\rm Reg}_r(\epsilon) \in  \langle \mathcal{D}_r( c) \rangle_{R}\}\]
explicitly determines $J\cdot{\rm Fitt}_{R}^i(\sha^2(\co_{K,S},\cA^\ast(1)))$ for every non-negative integer $i$.

\end{theorem}

\begin{proof}
By Hypotheses \ref{hyp surj} and \ref{hyp3'}, we know that ${\rm SS}_r(\cA)\simeq R$ (see Theorem \ref{structure unit 2}). We fix this isomorphism and define $J$ to be the inverse image of $\langle \mathcal{D}_r( c) \rangle_{R}$ under the surjection
$$R \simeq {\rm SS}_r(\cA) \stackrel{{\rm Reg}_r}{\to} {\rm KS}_r^{\rm reg}(\cA).$$
Then clearly we have $ \langle \mathcal{D}_r( c) \rangle_{R}= J\cdot {\rm KS}_r^{{\rm reg}}(\mathcal{A})$. The assertion follows from Theorem \ref{theorem higher fitt}(ii).
\end{proof}

Now we state the main result of this section.

\begin{theorem} \label{main}
Assume Hypothesis \ref{hyp last}. Then every basic Euler system $c \in \mathcal{E}^{\rm b}(T,\cK)$ is Stark-derivable, i.e. $\mathcal{D}_r(\mathcal{E}^{\rm b}(T,\cK)) \subset {\rm KS}_r^{\rm reg}(\cA)$.
\end{theorem}

The rest of this section is devoted to the proof of this theorem.

We need some preliminaries for the proof.
In \S \ref{section congruence}, we show that a certain congruence relation holds between a basic Euler system and a basic Stark system via `Bockstein maps' (see Theorem \ref{MRS}). In \S \ref{section bock}, we relate Bockstein maps with finite-singular comparison maps (see Theorem \ref{theorem bockstein}). Finally in \S \ref{last subsection}, we complete the proof of Theorem \ref{main}.

\subsubsection{Euler-Stark congruences} \label{section congruence}

Suppose that we have a vertical determinantal system $z \in {\rm VS}(T,\cK)$. By Theorem \ref{theorem vertical}, we can produce an Euler system $c(z) \in \mathcal{E}^{\rm b}(T,\cK)$ from $z$. Also, since we have a natural surjection
$${\rm VS}(T,\cK) \stackrel{z \mapsto z_E}{\to} {\det}_\cR(C_{E,S}(T)) \to {\rm det}_R(C_{K,S}(\cA)) \simeq {\rm HS}(\cA),$$
we can also produce a Stark system $\epsilon(z) \in {\rm SS}_r(\cA)$ from $z$, by Theorem \ref{theorem HU}.

In this subsection, we shall describe an explicit relation between
$$c(z)_{\fn} \in {\bigcap}_{\cR[\cG_\fn]}^r H_\Sigma^1(\co_{E(\fn),S_\fn}, T^\ast(1)) \text{ and }\epsilon(z)_\fn \in {\bigcap}_R^{r+\nu(\fn)} H_\Sigma^1(\co_{K,S_\fn},\cA^\ast(1))$$
for every $\fn \in \cN$ (see Theorem \ref{MRS} below).

We need some lemmas.

\begin{lemma}\label{cyclic lemma}
For any $\fn \in \cN$ and $\fq \mid \fn$, there is a natural surjection
\begin{eqnarray}\label{TA surjection}
H^2(K_\fq,T^\ast(1)_{E(\fn)}) \to H^2(K_\fq, \cA^\ast(1)), \nonumber
\end{eqnarray}
where $T^\ast(1)_{E(\fn)}$ denotes the induced module ${\rm Ind}_{G_K}^{G_{E(\fn)}}(T^\ast(1))$. Furthermore, the $\cR[\cG_\fn]$-module $H^2(K_\fq,T^\ast(1)_{E(\fn)})$ is cyclic, and generated by a lift $\widetilde \gamma_\fq$ of the generator $\gamma_\fq \in H^2(K_\fq, \cA^\ast(1))$.
\end{lemma}

\begin{proof}
We have
$$H^2(K_\fq, T^\ast(1)_{E(\fn)}) \otimes_{\cR[\cG_\fn]} R \simeq H^2(K_\fq, T^\ast(1)_{E(\fn)} \otimes_{\cR[\cG_\fn]} R) \simeq H^2(K_\fq,\cA^\ast(1))\stackrel{(\ref{identifications})}{=} R,$$
where the first isomorphism follows from the fact that the functor $H^2(K_\fq, \cdot)$ is right-exact. The lemma follows from Nakayama's lemma.
%
\end{proof}

\begin{lemma}
For $\fn \in \cN$, there is a quadratic standard representative of $C_{E(\fn),S_\fn}(T)$
$$P \stackrel{\psi}{\to} P$$
with respect to
$$f: H^1(C_{E(\fn), S_\fn}(T)) \to Y_K(T)^\ast \otimes_\cR \cR[\cG_\fn],$$
so that
$$\overline P \stackrel{\overline \psi}{\to} \overline P,$$
where we set $\overline P:=(P/MP)^{\cH_\fn} $ and $\overline \psi$ is the map induced by $\psi$, is a standard representative of $C_{K,S_\fn}(\cA)$ with respect to
$$\overline f: H^1(C_{K,S_\fn}(\cA)) \to \left(\bigoplus_{\fq \mid \fn} H^2(K_\fq,\cA^\ast(1)) \right)\oplus Y_K(\cA)^\ast.$$
\end{lemma}
\begin{proof}
Write $\fn =\prod_{i=1}^{\nu(\fn)} \fq_i$.
By the proof of Proposition \ref{prop admissible}(i), we can choose a quadratic standard representative $(P, \psi, \{b_1,\ldots,b_d\})$ of $C_{E(\fn),S_\fn}(T)$ with respect to the surjection
\begin{eqnarray}
H^1(C_{E(\fn),S_\fn}(T))&\simeq&H_\Sigma^2(\co_{E(\fn),S_\fn},T^\ast(1)) \oplus (Y_K(T)^\ast\otimes_\cR \cR[\cG_\fn]) \nonumber\\
&\to& \left(\bigoplus_{\fq \mid \fn} H^2(K_\fq,T^\ast(1)_{E(\fn)}) \right)\oplus (Y_K(T)^\ast \otimes_\cR \cR[\cG_\fn]), \nonumber
\end{eqnarray}
so that the image of $b_i$ under the map
$$P \to H^1(C_{E(\fn),S_\fn}(T))  \to \left(\bigoplus_{\fq \mid \fn} H^2(K_\fq,T^\ast(1)_{E(\fn)}) \right)\oplus (Y_K(T)^\ast \otimes_\cR \cR[\cG_\fn])$$
is $\beta_i$ (resp. $\widetilde \gamma_{\fq_{i-r}}$) for $1\leq i \leq r$ (resp. $r< i \leq r+\nu(\fn)$), where $\widetilde \gamma_\fq$ is as in Lemma \ref{cyclic lemma}.
Then we see that this representative satisfies the properties in the assertion.
\end{proof}

Fix $\fn \in \cN$, and let $(P,\psi, \{b_1,\ldots,b_d\})$ be as in the proof of the above lemma. Note that we can regard
$${\bigcap}_{\cR[\cG_\fn]}^r H_\Sigma^1(\co_{E(\fn),S_\fn},T^\ast(1))\subset {\bigwedge}_{\cR[\cG_\fn]}^r P \text{ and } {\bigcap}_R^{r+\nu(\fn)} H_\Sigma^1(\co_{K,S_\fn},\cA^\ast(1)) \subset {\bigwedge}_R^{r+\nu(\fn)} \overline P.$$

For each $i$ with $1\leq i \leq d$, we set
$$\psi_i:=b_i^\ast \circ \psi \in \Hom_{\cR[\cG_\fn]}(P,\cR[\cG_\fn]).$$
We often regard $\psi_i \in \Hom_{\fr[\cG_\fn]}(P/MP,\fr[\cG_\fn])$. Let $\cI_\fn$ be the augmentation ideal of $\ZZ[\cH_\fn]$. 

\begin{lemma} \label{mrs lemma}
\begin{itemize}
\item[(i)] For $r< i \leq r+\nu(\fn)$, we have
$$\im \psi_i  \subset  \cI_\fn \cdot \fr[\cG_\fn].$$
\item[(ii)] Regard $c(z)_\fn \in {\bigcap}_{\cR[\cG_\fn]}^r H^1(\co_{E(\fn),S_\fn},T^\ast(1))\subset {\bigwedge}_{\cR[\cG_\fn]}^r P$ as an element of ${\bigwedge}_{\fr[\cG_\fn]}^r (P/MP)$. Then we have
$$c(z)_\fn \in \cI_\fn^{\nu(\fn)} {\bigwedge}_{\fr[\cG_\fn]}^r (P/MP).$$
\end{itemize}
\end{lemma}

\begin{proof}
The assertion (i) follows from the choice of $(P ,\psi, \{b_1,\ldots,b_d\})$ (compare \cite[Lemma 5.20]{bks1}).
By the construction of $c(z)_\fn$, we have
$$c(z)_\fn \in \im \left( {\bigwedge}_{r < i \leq d}\psi_i: {\bigwedge}_{\fr[\cG_\fn]}^d (P/MP) \to {\bigwedge}_{\fr[\cG_\fn]}^r (P/MP) \right),$$
so the assertion (ii) follows from (i).
\end{proof}

\begin{lemma}
$X$ be a free $\fr[\cG_\fn]$-module, and $s$ a non-negative integer. Then there is a canonical homomorphism
\begin{eqnarray} \label{resolvent}
\cI_\fn^s X \to X^{\cH_\fn} \otimes \cI_\fn^s /\cI_\fn^{s+1}
\end{eqnarray}
such that the composition map
$$\cI_\fn^s X \to X^{\cH_\fn} \otimes \cI_\fn^s /\cI_\fn^{s+1} \to X \otimes \cI_\fn^s/\cI_\fn^{s+1} \to X \otimes \ZZ[\cH_\fn]/\cI_\fn^{s+1}$$
coincides with the map defined by
$$x \mapsto \sum_{\sigma\in \cH_\fn} \sigma x \otimes \sigma^{-1}.$$
\end{lemma}
\begin{proof}
Considering component-wise, we may assume $X$ is rank one. Let $e \in X$ be a basis.
Then each $x \in \cI_\fn^s X$ is uniquely written as $x=g e$ with $g \in \cI_\fn^s \cdot \fr[\cG_\fn]$. We write $g=\sum_{\gamma \in \Gamma} \widetilde \gamma \cdot g_\gamma$, where $\widetilde \gamma \in \cG_\fn$ is a lift of $\gamma \in \Gamma=\cG_\fn/\cH_\fn$ and $g_\gamma \in \cI_\fn^s \cdot \fr[\cH_\fn]$. Define (\ref{resolvent}) by
$$x \mapsto \sum_{\gamma \in \Gamma}\N_{\cH_\fn}  \widetilde \gamma e \otimes  g_\gamma,$$
where $\N_{\cH_\fn}:=\sum_{\sigma \in \cH_\fn} \sigma$ and we identify $X^{\cH_\fn} \otimes \cI_\fn^s/\cI_\fn^{s+1}$ with $X^{\cH_\fn} \otimes_{\fr} (\cI_\fn^s/\cI_\fn^{s+1} \otimes \fr)$.
It is straightforward to check that this is well-defined and that
$$\sum_{\sigma \in \cH_\fn} \sigma x \otimes \sigma^{-1}= \sum_{\gamma \in \Gamma}\N_{\cH_\fn}  \widetilde \gamma e \otimes  g_\gamma \text{ in }X\otimes \ZZ[\cH_\fn]/\cI_\fn^{s+1}.$$
Thus the lemma is proved.
\end{proof}

By the above two lemmas, we can regard
$$\sum_{\sigma \in \cH_\fn} c(z)_\fn \otimes \sigma^{-1} \in \left({\bigwedge}_{\fr[\cG_\fn]}^r (P/MP) \right)^{\cH_\fn} \otimes \cI_\fn^{\nu(\fn)}/\cI_\fn^{\nu(\fn)+1}= {\bigwedge}_R^r \overline P\otimes \cI_\fn^{\nu(\fn)}/\cI_\fn^{\nu(\fn)+1}.$$

For $\fq \mid \fn$, we define the `Bockstein map'
$$\varphi_\fq^\fn : \overline P \to \cI_\fn \cdot \fr[\cG_\fn]/ \cI_\fn^2 \cdot \fr[\cG_\fn]=R\otimes \cI_\fn/\cI_\fn^2$$
by
$$\varphi_\fq^\fn( {\N}_{\cH_\fn} a ) := \psi_i(a) \text{ mod $\cI_\fn^2\cdot \fr[\cG_\fn]$}\quad (a \in P/MP),$$
where $i$ is determined by $\fq = \fq_{i-r}$ (note that $\im \psi_i \subset  \cI_\fn \cdot \fr[\cG_\fn]$ by Lemma \ref{mrs lemma}(i)). One checks that $\varphi_\fq^\fn$ is well-defined.
Note that Bockstein maps give the map
$${\bigwedge}_{\fq \mid \fn} \varphi_\fq^\fn: {\bigwedge}_{R}^{r+\nu(\fn)} \overline P \to {\bigwedge}_{R}^r \overline P \otimes \cI_\fn^{\nu(\fn)}/\cI_\fn^{\nu(\fn)+1}.$$

\begin{theorem} \label{MRS}
We have an equality
\begin{eqnarray}
\sum_{\sigma\in \cH_\fn}\sigma c(z)_\fn \otimes \sigma^{-1} =  \left({\bigwedge}_{\fq \mid \fn }\varphi_\fq^\fn \right) (\epsilon(z)_\fn) \label{eq mrs}
\end{eqnarray}
in ${\bigwedge}_R^r \overline P \otimes \cI_\fn^{\nu(\fn)}/\cI_\fn^{\nu(\fn)+1}$.
\end{theorem}

\begin{proof}
Let
$$\N_{\cH_\fn}^d: {\bigwedge}_{\fr[\cG_\fn]}^d (P/MP) \to {\bigwedge}_R^d \overline P$$
be the map induced by the norm map $\N_{\cH_\fn}: P/MP \to \overline P$. For $r+\nu (\fn) <i \leq d$, define
$$\overline \psi_i:=(\N_{\cH_\fn} b_i)^\ast \circ \overline \psi \in \Hom_R(\overline P, R),$$
where $\{\N_{\cH_\fn} b_i\}_i$ is regarded as a basis of $\overline P$.
By computation using the formula (\ref{exterior explicit}), one checks that the following diagram is commutative (up to sign).
$$\xymatrix{
{\bigwedge}_{\fr[\cG_\fn]}^d (P/MP) \ar[d]_{\N_{\cH_\fn}^d} \ar[r]^{{\bigwedge}_{r < i \leq d}\psi_i} & \quad \cI_\fn^{\nu(\fn)}   {\bigwedge}_{\fr[\cG_\fn]}^r (P/MP) \ar[r]^{\sum_{\sigma \in \cH_\fn}\sigma (\cdot)\otimes \sigma^{-1}} &\quad \quad ({\bigwedge}_{\fr[\cG_\fn]}^r (P/MP))^{\cH_\fn} \otimes \cI_\fn^{\nu(\fn)}/\cI_\fn^{\nu(\fn)+1}  \\
{\bigwedge}_R^d \overline P \ar[r]_{{\bigwedge}_{r+\nu(\fn) < i \leq d}\overline \psi_i}  &{\bigwedge}_R^{r+\nu(\fn)} \overline P \ar[r]_{{\bigwedge}_{\fq \mid \fn}\varphi_\fq^\fn} &{\bigwedge}_R^r \overline P \otimes \cI_\fn^{\nu(\fn)}/\cI_\fn^{\nu(\fn)+1} \ar[u]_\simeq.
}
$$
(Compare \cite[Lemma 5.22]{bks1}.)
By the constructions of $c(z)_\fn$ and $\epsilon(z)_\fn$, we have
$$c(z)_\fn= \pm\left( {\bigwedge}_{r <i \leq d}\psi_i \right) (z_\fn')$$
and
$$\epsilon(z)_\fn=\pm \left({\bigwedge}_{r+\nu(\fn)< i \leq d} \overline \psi_i \right)(\N_{\cH_\fn}^d z_\fn')$$
respectively,
where $z_\fn' \in {\bigwedge}_{\fr[\cG_\fn]}^d (P/MP)$ is the element satisfying
$$z_{E(\fn)} = z_\fn' \otimes b_1^\ast \wedge \cdots \wedge b_d^\ast \text{ in } {\det}_{\cR[\cG_\fn]}(C_{E(\fn),S_\fn}(T)) \otimes_\cR \fr = {\bigwedge}_{\fr[\cG_\fn]}^d (P/MP) \otimes {\bigwedge}_{\fr[\cG_\fn]}^d (P/MP)^\ast.$$
Hence, the theorem follows from the above commutative diagram and explicit computation of sign.
\end{proof}

\begin{remark} \label{remark mrs}
The relation (\ref{eq mrs}) is regarded as a variant of Darmon's conjecture on cyclotomic units \cite[Conjecture 4.3]{D} for general $p$-adic representations. In particular, (\ref{eq mrs}) is regarded as a variant of the conjecture on Rubin-Stark elements proposed by Mazur-Rubin in \cite[Conjecture 5.2]{MRGm} and by the second author in \cite[Conjecture 3]{sano}.  Also, we remark that the second author established a `non-explicit version' of (\ref{eq mrs}) in \cite[Theorem 3.8]{sanojnt} in the rank one case, so Theorem \ref{MRS} can also be regarded as its improvement.
\end{remark}

\subsubsection{Computation of Bockstein maps} \label{section bock}

In this subsection, we relate Bockstein maps with finite-singular comparison maps.

First, we give a conceptual definition of Bockstein maps. Recall that in the previous subsection we defined the Bockstein map
$$\varphi_\fq^\fn: \overline P \to R \otimes \cI_\fn/\cI_\fn^2.$$
We shall define a map
\begin{eqnarray} \label{local bock}
 H^1(K_\fq,\cA^\ast(1)) \to R \otimes \cI_\fn/\cI_\fn^2
 \end{eqnarray}
such that the composition map
$$H_\Sigma^1(\co_{K,S_\fn},\cA^\ast(1)) \to  H^1(K_\fq,\cA^\ast(1)) \to R \otimes \cI_\fn/\cI_\fn^2$$
coincides with $\varphi_\fq^\fn$ restricted on $H_\Sigma^1(\co_{K,S_\fn},\cA^\ast(1)) (\subset \overline P)$. The map (\ref{local bock}) will also be denoted by $\varphi_\fq^\fn$.

Set $\mathcal{D}_\fq :=\Gal(E(\fn)_\fq/K_\fq)$, and let $\cJ_\fq$ be the augmentation ideal of $\ZZ[\mathcal{D}_\fq]$. Note that $\mathcal{D}_\fq$ is identified with the decomposition group at $\fq $ in $\cG_\fn$, so we regard $\mathcal{D}_\fq \subset \cG_\fn$. Since $\fq$ splits completely in $E$, we have $\mathcal{D}_\fq \subset \cH_\fn$.
The Bockstein map is the following map.

\begin{eqnarray}\label{bock delta}
\varphi_\fq^\fn: H^1(K_\fq, \cA^\ast(1)) &\stackrel{\delta}{\to} & H^2(K_\fq,\cA^\ast(1)\otimes \cJ_\fq)  \\
&\to& H^2(K_\fq, \cA^\ast(1)\otimes \cJ_\fq/\cJ_\fq^2) \nonumber\\
&=& H^2(K_\fq, \cA^\ast(1))\otimes \cJ_\fq/\cJ_\fq^2 \nonumber\\
&\stackrel{(\ref{identifications})}{=}& R\otimes \cJ_\fq/\cJ_\fq^2 \nonumber\\
&\to& R \otimes  \cI_\fn/\cI_\fn^2, \nonumber
\end{eqnarray}
where $\delta $ is the boundary map with respect to the short exact sequence
$$0 \to \cA^\ast(1) \otimes \cJ_\fq \to \cA^\ast(1) \otimes \ZZ[\mathcal{D}_\fq] \to \cA^\ast(1) \to 0,$$
and the last map is the map induced by the inclusion $\mathcal{D}_\fq \hookrightarrow \cH_\fn$.


One checks that this definition coincides with the definition of Bockstein maps given in the previous subsection.

The aim of this subsection is to prove the following theorem.

\begin{theorem} \label{theorem bockstein}
For $\fn \in \cN$ and $\fq \mid \fn$, we have
$$ \varphi_\fq^\fn=\rho_\fq \circ \varphi_\fq^{\rm fs} + P_\fq({\rm Fr}_\fq^{-1}) \cdot v_\fq \text{ in }\Hom_R(H_\Sigma^1(\co_{K,S_\fn},\cA^\ast(1)),R\otimes \cI_\fn/\cI_\fn^2 ). $$
Here $\rho_\fq \circ \varphi_\fq^{\rm fs}$ stands for the map
$$H_\Sigma^1(\co_{K,S_\fn},\cA^\ast(1)) \stackrel{\varphi_\fq^{\rm fs}}{\to} R \otimes G_\fq \stackrel{\rho_\fq}{\stackrel{\sim}{\to}} R \otimes \cI_\fq/\cI_\fq^2 \subset  R \otimes \cI_\fn/\cI_\fn^2,$$
and $P_\fq({\rm Fr}_\fq^{-1}) $ is reagrded as an element of $\fr \otimes \cI_{\fn/\fq}/\cI_{\fn/\fq}^2\subset  \fr \otimes \cI_\fn/\cI_\fn^2$.
\end{theorem}

\begin{remark}
The above theorem shows that the map $\varphi_\ell^n$ defined by the second author in \cite[\S 3]{sanojnt} is a Bockstein map. Indeed, $\varphi_\ell^n$ is defined in loc. cit. by
$$\varphi_\ell^n:=-(\sigma_\ell-1)\cdot u_\ell -P_\ell({\rm Fr}_\ell)\cdot v_\ell,$$
where $-(\sigma_\ell-1)\cdot u_\ell$, $-P_\ell({\rm Fr}_\ell) $, and $v_\ell$ correspond to our $\rho_\fq\circ \varphi_\fq^{\rm fs}$, $P_\fq({\rm Fr}_\fq^{-1})$, and $v_\fq$ respectively.
\end{remark}

To prove Theorem \ref{theorem bockstein}, we need some lemmas.

\begin{lemma} \label{bock lemma}
$\varphi_\fq^\fn$ coincides with
\begin{eqnarray}
H^1(K_\fq,\cA^\ast(1)) &\stackrel{(\ref{LD})}{\simeq}& H^1(K_\fq, \cA)^\ast  \nonumber \\
&{\to}& H^1(\mathcal{D}_\fq, \cA^{G_{E(\fn)_\fq}})^\ast \label{bock 2} \\
&\to& H^1(\mathcal{D}_\fq,\cA^{{\rm Fr_\fq}=1})^\ast  \label{bock 3} \\
&=&\Hom_R(\Hom_\ZZ(\mathcal{D}_\fq, \cA^{{\rm Fr}_\fq=1}), R)\nonumber \\
&\simeq& (\cA^{{\rm Fr_\fq}=1})^\ast \otimes \cJ_\fq/\cJ_\fq^2 \label{bock 4} \\
&\stackrel{(\ref{identifications})}{=} & R \otimes \cJ_\fq/\cJ_\fq^2, \nonumber
\end{eqnarray}
where (\ref{bock 2}) is the dual of the inflation map, (\ref{bock 3}) is induced by the inclusion $\cA^{{\rm Fr_\fq}=1} \hookrightarrow \cA^{G_{E(\fn)_\fq}}$, (\ref{bock 4}) is the inverse of
\begin{eqnarray}
(\cA^{{\rm Fr_\fq}=1})^\ast \otimes \cJ_\fq/\cJ_\fq^2 &\stackrel{\sim}{\to}& \Hom_R(\Hom_\ZZ(\mathcal{D}_\fq, \cA^{{\rm Fr}_\fq=1}), R) \nonumber \\
f \otimes (\sigma -1) &\mapsto& (g \mapsto f(g(\sigma))).\nonumber
\end{eqnarray}
\end{lemma}
\begin{proof}
For simplicity, we put $D:=\mathcal{D}_\fq, \cJ:=\cJ_\fq, \cA^F:=\cA^{{\rm Fr}_\fq=1}$. We also abbreviate $H^i(K_\fq,\cdot)$ and $\widehat H^i(D, \cdot)$ as $H^i(\cdot)$ and $\widehat H^i(\cdot)$ respectively.

By the duality between homology and cohomology, we have
\begin{eqnarray}
H_1(D, (\cA^F)^\ast) \simeq H^1(D,\cA^F)^\ast. \nonumber
\end{eqnarray}
This isomorphism is explicitly given by
\begin{eqnarray}\label{hom cohom}
H_1(D, (\cA^F)^\ast)=\widehat H^{-2}( (\cA^F)^\ast) \stackrel{\sim}{\to} \widehat H^1(\cA^F)^\ast = H^1(D,\cA^F)^\ast
\end{eqnarray}
$$ x \mapsto (y \mapsto x\cup y),$$
where $x \cup y \in \widehat H^{-1}(R)=R$ (see \cite[Chap. VI, Corollary 7.3]{brown}). One checks that, by a similar computation to \cite[Appendix to Chap. XI, Lemma 3]{serre}, (\ref{bock 4}) coincides with the composition map
$$ \widehat H^1(\cA^F)^\ast \stackrel{(\ref{hom cohom})}{\simeq} \widehat H^{-2}((\cA^F)^\ast) \stackrel{\delta'}{\stackrel{\sim}{\to}} \widehat H^{-1}((\cA^F)^\ast \otimes \cJ)=(\cA^F)^\ast \otimes \cJ/\cJ^2,$$
where $\delta'$ is the boundary map with respect to the short exact sequence
$$0 \to (\cA^F)^\ast \otimes \cJ \to (\cA^F)^\ast \otimes \ZZ[D] \to (\cA^F)^\ast \to 0.$$
Hence, it is sufficient to prove that the following diagram is commutative.
$$\xymatrix{
H^1(\cA^\ast(1))  \ar[d]_{(\ref{LD})}^{\simeq} \ar[r]^{\delta} & H^2(\cA^\ast(1) \otimes \cJ) \ar[r] & H^2(\cA^\ast(1)) \otimes \cJ/\cJ^2  \ar[r]^{\stackrel{(\ref{LD})}{\sim}} &(\cA^F)^\ast \otimes \cJ/\cJ^2  \ar[d]^{=}\\
H^1(\cA)^\ast  \ar[r]_{(\ref{bock 2}) \text{ and }(\ref{bock 3})}& \widehat H^1(\cA^F)^\ast  & \ar[l]_{\sim}^{(\ref{hom cohom})} \widehat H^{-2}((\cA^F)^\ast) \ar[r]_{\delta'}^\sim& \widehat H^{-1}((\cA^F)^\ast \otimes \cJ),
}
$$
where $\delta$ is as in (\ref{bock delta}).
One can check the commutativity by explicit computation, using the compatibility between boundary maps and cup products.
\end{proof}

\begin{lemma} \label{lemma fs comparison}
The map
$$\varphi_\fq^{\rm fs}: H^1_f(K_\fq, \cA^\ast(1)) \stackrel{\sim}{\to} H^1_{/f}(K_\fq,\cA^\ast(1)) \otimes G_\fq$$
coincides with the composition of the map
$$H^1_f(K_\fq, \cA^\ast(1)) \stackrel{(\ref{LDf})}{\stackrel{\sim}{\to}} H^1_{/f}(K_\fq,\cA)^\ast \simeq  \Hom(G_\fq, \cA^{{\rm Fr}_\fq=1})^\ast$$
and the inverse of
\begin{eqnarray}
H^1_{/f}(K_\fq,\cA^\ast(1)) \otimes G_\fq &\stackrel{\sim}{\to}&  (\cA^{{\rm Fr}_\fq=1})^\ast \otimes G_\fq \nonumber \\
& \stackrel{\sim}{\to}& \Hom(G_\fq, \cA^{{\rm Fr}_\fq=1})^\ast, \label{a hom}
\end{eqnarray}
where the first isomorphism is induced by $v_\fq$, and (\ref{a hom}) is given by
$$f \otimes \sigma \mapsto (g \mapsto f(g(\sigma))).$$
\end{lemma}

\begin{proof}
For simplicity, we denote $(\cdot)^{{\rm Fr}_\fq=1}$ and $(\cdot)/({\rm Fr}_\fq-1)(\cdot)$ by $(\cdot)^F$ and $(\cdot)_F$ respectively.

It is sufficient to show that the following diagram is commutative.
$$\xymatrix{
H^1_f(K_\fq, \cA^\ast(1))\ar[d]_{(\ref{LDf})}^{\simeq} \ar[r]^{\quad \stackrel{(\ref{fs 1})}{\sim}} & \cA^\ast(1)_F \ar[r]^{\stackrel{Q_\fq({\rm Fr}_\fq^{-1})}{\sim}} & \cA^\ast(1)^F &  \ar[l]_{\stackrel{(\ref{fs 2})}{\sim} \quad}  \ar[d]_{\simeq}^{(\ref{LD/f})} H^1_{/f}(K_\fq, \cA^\ast(1)) \otimes G_\fq  \\
H^1_{/f}(K_\fq, \cA)^\ast & \ar[l]_{\sim}^{(\ref{a hom})} (\cA^F)^\ast \otimes G_\fq \ar[r]_{-Q_\fq({\rm Fr}_\fq)^\ast}^{\sim}& (\cA_F)^\ast \otimes G_\fq  \ar[r]^{\sim}_{(\ref{frobenius isom})^\ast} & H^1_f(K_\fq, \cA)^\ast \otimes G_\fq
}
$$
One sees by \cite[Lemma 1.4.7(ii)]{R} that this diagram commutes. (There is a sign ambiguity in loc. cit., but we can determine the sign.)
\end{proof}

\begin{proof}[Proof of Theorem \ref{theorem bockstein}]
Recall
$$H^1(K_\fq,\cA^\ast(1))= H_f^1(K_\fq,\cA^\ast(1)) \oplus H_{\rm tr}^1(K_\fq,\cA^\ast(1)).$$
Using Lemmas \ref{bock lemma} and \ref{lemma fs comparison}, it is straightforward to see that $\varphi_\fq^\fn$ coincides with $\rho_\fq\circ \varphi_\fq^{\rm fs}$ on $H_f^1(K_\fq,\cA^\ast(1)) $. We show that $\varphi_\fq^\fn$ coincides with $P_\fq({\rm Fr}_\fq^{-1}) \cdot v_\fq$ on $H_{\rm tr}^1(K_\fq,\cA^\ast(1))$.

Consider the map
$$H^1_f(K_\fq,\cA)^\ast = H^1(\langle {\rm Fr}_\fq \rangle, \cA)^\ast \to \Hom(\langle {\rm Fr}_\fq \rangle, \cA^{{\rm Fr}_\fq=1})^\ast$$
induced by the inclusion $\cA^{{\rm Fr}_\fq=1}\hookrightarrow \cA$. If we identify $H^1(\langle {\rm Fr}_\fq \rangle, \cA)$ with $\cA/({\rm Fr}_\fq-1)\cA$ and $\Hom(\langle {\rm Fr}_\fq \rangle, \cA^{{\rm Fr}_\fq=1})$ with $\cA^{{\rm Fr}_\fq=1}$ by evaluating ${\rm Fr}_\fq$ (as in (\ref{frobenius isom})), the above map coincides with the dual of the natural map
$$\cA^{{\rm Fr}_\fq=1} \to \cA/({\rm Fr}_\fq-1)\cA; \ a \mapsto \overline a.$$
Since $a\in \cA^{{\rm Fr}_\fq=1} $ is written as
$$a=Q_\fq({\rm Fr}_\fq)b$$
with some $b \in \cA$ (see (\ref{q isom 2})), we have
$$\overline a = Q_\fq(1 ) \cdot \overline b=Q_\fq(1)\cdot Q_\fq({\rm Fr}_\fq)^{-1}(a) \text{ in } \cA/({\rm Fr}_\fq-1)\cA.$$
Hence, by Lemma \ref{bock lemma} and the definition of $v_\fq$, we see that $\varphi_\fq^\fn$ on $H_{\rm tr}^1(K_\fq,\cA^\ast(1))$ coincides with the map
$$H_{\rm tr}^1(K_\fq,\cA^\ast(1)) \to R \otimes \cI_{\fn/\fq}/\cI_{\fn/\fq}^2; \ e \mapsto -Q_\fq(1)\cdot v_\fq(e) \otimes ({\rm Fr}_\fq -1).$$
Since we have
$$-1\otimes Q_\fq(1)({\rm Fr}_\fq-1) =1\otimes Q_\fq({\rm Fr}_\fq^{-1})({\rm Fr}_\fq^{-1}-1)=1\otimes P_\fq({\rm Fr}_\fq^{-1})$$
in $R \otimes \cI_{\fn/\fq}/\cI_{\fn/\fq}^2$, we have proved the proposition.
\end{proof}

\subsubsection{The proof of Theorem \ref{main}} \label{last subsection}

In this section, we construct a canonical homomorphism
$${\rm SS}_r(\cA) \to {\rm DKS}_r(\cA),$$
and show the compatibility with the map ${\rm Reg}_r$ and $\Psi_r$ (see Theorem \ref{commutative theorem}). The method in this subsection is similar to that of \cite[\S\S 4 and 5]{sanojnt}. Finally, we complete the proof of Theorem \ref{main}.

We begin with some algebraic lemmas.

\begin{lemma}\
\begin{itemize}
\item[(i)] The group $\cI_\fn^{\nu(\fn)}/\cI_\fn^{\nu(\fn)+1}$ decomposes as a direct sum
$$\left\langle \prod_{\fq \mid \fn} (\sigma_\fq-1) \right\rangle \oplus \left\langle \prod_{i=1}^{\nu(\fn)} (\sigma_{\fq_i} -1) \ \middle| \ \fq_i \mid \fn \text{ for every $i$, and $\fq_i=\fq_j$ for some $1\leq i < j \leq \nu(\fn)$}  \right\rangle$$
and there is an isomorphism
\begin{eqnarray} \label{g isom}
\left\langle \prod_{\fq \mid \fn} (\sigma_\fq-1) \right\rangle  \stackrel{\sim}{\to} G_\fn; \ \prod_{\fq \mid \fn} (\sigma_\fq-1) \mapsto \bigotimes_{\fq \mid \fn} \sigma_\fq.
\end{eqnarray}
\item[(ii)]
Let
$$s_\fn: \cI_\fn^{\nu(\fn)}/\cI_\fn^{\nu(\fn)+1} \to \left \langle \prod_{\fq \mid \fn} (\sigma_\fq-1) \right\rangle$$
be the projection map with respect to the decomposition in (i). Then one has 
\begin{eqnarray}\label{definition s}
s_\fn =\sum_{\fd \mid \fn}(-1)^{\nu(\fn/\fd)} \pi_\fd ,
\end{eqnarray}
where
$$\pi_\fd: \cI_\fn^{\nu(\fn)}/\cI_\fn^{\nu(\fn)+1} \to \cI_\fd^{\nu(\fn)}/\cI_\fd^{\nu(\fn)+1} \subset \cI_\fn^{\nu(\fn)}/\cI_\fn^{\nu(\fn)+1}$$
is the map induced by the natural projection $\cH_\fn \to \cH_\fd \subset \cH_\fn$.
\end{itemize}
\end{lemma}

\begin{proof}
The statement (i) is proved in \cite[Proposition 4.2(i) and (iv)]{MR}.

We show (ii). We set $s_\fn':=\sum_{\fd \mid \fn}(-1)^{\nu(\fn/\fd)} \pi_\fd $ and show $s_\fn=s_\fn'$. Since $\pi_\fd(\prod_{\fq \mid \fn}(\sigma_\fq-1)) = 0 $ if $\fd \neq \fn$, we see that $s_\fn'=\id$ on $\langle \prod_{\fq \mid \fn} (\sigma_\fq-1)\rangle$. We shall show that
$$s_\fn' \left( \prod_{i=1}^{\nu(\fn)} (\sigma_{\fq_i} -1) \right)=0$$
if $\fq_i =\fq_j$ for some $1\leq i< j \leq \nu(\fn)$. Let $\fc$ be the square-free product of $\fq \mid \fn$ which does not divide $\prod_{i=1}^{\nu(\fn)} \fq_i$. Then we have
$$s_\fn' \left( \prod_{i=1}^{\nu(\fn)} (\sigma_{\fq_i} -1) \right)= \sum_{\fd \mid \fc} (-1)^{\nu(\fd)}  \prod_{i=1}^{\nu(\fn)} (\sigma_{\fq_i} -1) .$$
The required result then follows since $\fc \neq 1$ and so $\sum_{\fd \mid \fc}(-1)^{\nu(\fd)}=(1-1)^{\nu(\fc)}=0.$
\end{proof}

The endomorphism of $\ZZ[\cH_\fn]$ induced by $\cH_\fn \to \cH_\fd \subset \cH_\fn$ is also denoted by $\pi_\fd$. The right hand side of (\ref{definition s}) gives an endomorphism of $\ZZ[\cH_\fn]$, which we denote also by $s_\fn$.

\begin{lemma} \label{derivative lemma}\
\begin{itemize}
\item[(i)] The image of the map $s_\fn : \ZZ[\cH_\fn] \to \ZZ[\cH_\fn]$ is contained in $\cI_\fn^{\nu(\fn)}$.
\item[(ii)] For each $\ZZ[\cH_\fn]$-module $X$ we regard $s_\fn$ as a map $X \otimes \ZZ[\cH_\fn] \to X \otimes \cI_\fn^{\nu(\fn)}$. Then for each $x \in X$, one has
$$s_\fn\left(\sum_{\sigma\in \cH_\fn} \sigma x \otimes \sigma^{-1} \right) =(-1)^{\nu(\fn)}D_\fn x \otimes \prod_{\fq \mid \fn}(\sigma_\fq-1) \text{ in }X \otimes \cI_\fn^{\nu(\fn)}/\cI_\fn^{\nu(\fn)+1}.$$
\end{itemize}
\end{lemma}

\begin{proof}
(i) For any $\sigma \in \cH_\fn$, we have
$$s_\fn(\sigma)=\sum_{\fd \mid \fn}(-1)^{\nu(\fn/\fd)}\pi_\fd(\sigma) =\prod_{\fq \mid \fn}(\pi_\fq(\sigma) -1) \in \cI_\fn^{\nu(\fn)}.$$
Thus $\im (s_\fn) \subset \cI_\fn^{\nu(\fn)}$.

(ii) Write $\fn=\prod_{i=1}^{\nu(\fn)} \fq_i$, and set
$$\Lambda_\fn:= \{ (i_1,\ldots,i_{\nu(\fn)}) \in \ZZ^{\oplus\nu(\fn)} \mid 0 \leq i_j <\#G_{\fq_j}\text{ for every $1 \leq j \leq \nu(\fn)$} \}.$$
Note that the map
$$\Lambda_\fn \to \cH_\fn; \ (i_1,\ldots,i_{\nu(\fn)}) \mapsto \sigma_{\fq_1}^{i_1}\cdots \sigma_{\fq_{\nu(\fn)}}^{i_{\nu(\fn)}} $$
is bijective.
We then explicitly compute that $s_\fn\left(\sum_{\sigma\in \cH_\fn} \sigma x \otimes \sigma^{-1} \right)$ is equal to

\begin{eqnarray}
& & \sum_{\fd \mid \fn}(-1)^{\nu(\fn/\fd)} \sum_{\sigma \in \cH_\fn} \sigma x \otimes \pi_\fd(\sigma^{-1}) \nonumber\\
&=& \sum_{\sigma \in \cH_\fn}\sigma x \otimes (\pi_{\fq_1}(\sigma^{-1}) -1) \cdots (\pi_{\fq_{\nu(\fn)}}(\sigma^{-1}) -1) \nonumber\\
& =& \sum_{(i_1,\ldots,i_{\nu(\fn)})\in \Lambda_\fn}  \sigma_{\fq_1}^{i_1}\cdots \sigma_{\fq_{\nu(\fn)}}^{i_{\nu(\fn)}} x \otimes (\sigma_{\fq_1}^{-i_1}-1) \cdots (\sigma_{\fq_{\nu(\fn)}}^{-i_{\nu(\fn)}}-1) \nonumber\\
&=&(-1)^{\nu(\fn)}\sum_{(i_1,\ldots,i_{\nu(\fn)}) \in \Lambda_\fn} i_1\cdots i_{\nu(\fn)}\sigma_{\fq_1}^{i_1}\cdots \sigma_{\fq_{\nu(\fn)}}^{i_{\nu(\fn)}} x \otimes (\sigma_{\fq_1}-1) \cdots (\sigma_{\fq_{\nu(\fn)}}-1) \nonumber\\
&=&(-1)^{\nu(\fn)}D_\fn x \otimes \prod_{\fq\mid \fn}(\sigma_\fq-1). \nonumber
\end{eqnarray}

\end{proof}

We shall define a homomorphism
\begin{eqnarray} \label{US to DKS}
{\rm SS}_r(\cA) \to {\rm KC}_r(\cA) \left(=\prod_{\fn \in \cN} {\bigcap}_R^r H_{\cF^\fn}^1(K,\cA^\ast(1)) \otimes G_\fn\right)
\end{eqnarray}
as follows. Note that
$$\Hom_R(\cdot ,R\otimes \cI_\fn/\cI_\fn^2) \simeq \Hom_R(\cdot,R)\otimes \cI_\fn/\cI_\fn^2,$$
since $R \otimes \cI_\fn/\cI_\fn^2 \simeq \bigoplus_{\fq\mid \fn} R \otimes G_\fq \simeq R^{\oplus \nu(\fn)}$. So ${\bigwedge}_{\fq \mid \fn} \varphi_\fq^\fn $ induces a map
$${\bigwedge}_{\fq \mid \fn}\varphi_\fq^\fn: {\bigcap}_R^{r+\nu(\fn)} H^1_{\cF^\fn}(K,\cA^\ast(1)) \to {\bigcap}_R^r H^1_{\cF^\fn}(K,\cA^\ast(1)) \otimes \cI_\fn^{\nu(\fn)}/\cI_\fn^{\nu(\fn)+1}.$$
We denote the map
\begin{eqnarray}
{\bigcap}_R^r H^1_{\cF^\fn}(K,\cA^\ast(1)) \otimes \cI_\fn^{\nu(\fn)}/\cI_\fn^{\nu(\fn)+1} &\stackrel{s_\fn}{\to}& {\bigcap}_R^r H^1_{\cF^\fn}(K,\cA^\ast(1)) \otimes \left \langle \prod_{\fq \mid \fn} (\sigma_\fq-1) \right\rangle \nonumber \\
& \stackrel{(\ref{g isom})}{\simeq}& {\bigcap}_R^r H^1_{\cF^\fn}(K,\cA^\ast(1)) \otimes G_\fn  \nonumber
\end{eqnarray}
also by $s_\fn$.
We define (\ref{US to DKS}) to be the map
$${\rm SS}_r(\cA) \to {\rm KC}_r(\cA); \ \epsilon \mapsto \left\{ (-1)^{\nu(\fn)}s_\fn \left( \left({\bigwedge}_{\fq \mid \fn}\varphi_\fq^\fn\right) (\epsilon_\fn) \right) \right\}_\fn.$$


\begin{theorem}[{\cite[Theorem 5.7]{sanojnt}}] \label{commutative theorem}
The diagram
$$\xymatrix{
 {\rm SS}_r(\cA)   \ar[r]^{(\ref{US to DKS}) } \ar[d]_{{\rm Reg}_r}& \ar[d]^{\Psi_r} {\rm KC}_r(\cA) \\
 {\rm KS}_r(\cA)  \ar[r]^{\subset  } &  {\rm KC}_r(\cA)}
$$
commutes. In particular, by Theorem \ref{isom dks ks}, the image of (\ref{US to DKS}) is contained in ${\rm DKS}_r(\cA)$, and we have the commutative diagram
$$\xymatrix{
 {\rm SS}_r(\cA)   \ar[r]^{(\ref{US to DKS}) } \ar[d]_{{\rm Reg}_r}&  \ar[ld]_\simeq^{\Psi_r} {\rm DKS}_r(\cA)  \\
 {\rm KS}_r(\cA) &
}
$$
\end{theorem}

\begin{proof}
We identify $G_\fn$ with $\langle \prod_{\fq \mid \fn} (\sigma_\fq-1)\rangle$ by (\ref{g isom}) for any $\fn \in \cN$. In particular, we identify $G_\fq $ with $\cI_\fq/\cI_\fq^2$ via $\rho_\fq$ for any $\fq \in \cP$. By Theorem \ref{theorem bockstein}, note that we can identify $\varphi_\fq^{\rm fs}$ with $\varphi_\fq^\fq$.
For $\fn \in \cN$ and $\fq \in \cP$ with $\fq \not{|} \fn$, we put
$$P_\fq^{\fn}:=P_\fq({\rm Fr}_\fq^{-1}) \in \fr \otimes \cI_{\fn}/\cI_{\fn}^2.$$
Take a Stark system $\epsilon=\{\epsilon_\fn\}_\fn \in {\rm SS}_r(\cA)$.
It is sufficient to show that the equality
\begin{multline}\label{theta kolyvagin}
\sum_{\tau \in \mathfrak{S}(\fn)}{\rm sgn}(\tau) (-1)^{\nu(\fd_\tau)}s_{\fd_\tau} \left( \left({\bigwedge}_{\fq \mid \fd_\tau}\varphi_\fq^{\fd_\tau}\right) (\epsilon_{\fd_\tau}) \right) \prod_{\fq \mid \fn/\fd_\tau}P_{\tau(\fq)}^\fq=(-1)^{\nu(\fn)} \left( {\bigwedge}_{\fq\mid \fn}\varphi_\fq^\fq \right)(\epsilon_\fn)
\end{multline}
holds for every $\fn \in \cN$. Note that, for $\fq \mid \fn$, we have
$$\varphi_\fq^\fq=\varphi_\fq^\fn - P_\fq^{\fn/\fq} \cdot v_\fq$$
by Theorem \ref{theorem bockstein}. Using this, we compute
\begin{eqnarray} \label{compute1}
\left( {\bigwedge}_{\fq\mid \fn}\varphi_\fq^\fq \right)(\epsilon_\fn) =\sum_{\fd \mid \fn}(-1)^{\nu(\fd)} \left( {\bigwedge}_{\fq\mid \fn/\fd}\varphi_\fq^\fn \right)(\epsilon_{\fn/\fd}) \prod_{\fq \mid \fd} P_\fq^{\fn/\fq}.
\end{eqnarray}
If $\fq \mid \fn/ \fd$, then using
$$\varphi_\fq^\fn=\varphi_\fq^{\fn/\fd} +P_\fq^{\fd}\cdot v_\fq$$
we compute
$$\left( {\bigwedge}_{\fq\mid \fn/\fd}\varphi_\fq^\fn \right)(\epsilon_{\fn/\fd})=\sum_{\fc \mid \fn/\fd} \left( {\bigwedge}_{\fq\mid \fn/\fc\fd}\varphi_\fq^{\fn/\fd} \right)(\epsilon_{\fn/\fc\fd}) \prod_{\fq \mid \fc}P_\fq^{\fd}.  $$
Using this formula recursively, we see that the right hand side of (\ref{compute1}) is equal to
\begin{eqnarray}
\sum_{\fd \mid \fn}\left( \sum_{(\fc_1,\ldots,\fc_k)\in \Delta(\fd)} (-1)^{\nu(\fc_k)} \prod_{\fq\mid \fc_k}P_\fq^{\fn/\fq}\prod_{\fq \mid \fc_{k-1}}P_\fq^{\fc_k} \cdots \prod_{\fq \mid \fc_1}P_\fq^{\fc_2} \right)\cdot \left( {\bigwedge}_{\fq\mid \fn/\fd}\varphi_\fq^{\fn/\fd} \right)(\epsilon_{\fn/\fd}), \nonumber
\end{eqnarray}
where
$$\Delta(\fd):=\left\{ (\fc_1,\ldots,\fc_k)\  \middle| \  \text{$k \in \ZZ_{\geq 1}$, $1\neq \fc_i \mid \fd$ for every $i$, and $\fd=\prod_{i=1}^k \fc_i$} \right\}.$$

By \cite[Lemma 4.19]{sanojnt}, we have

\begin{eqnarray}
&& \sum_{(\fc_1,\ldots,\fc_k)\in \Delta(\fd)} (-1)^{\nu(\fc_k)} \prod_{\fq\mid \fc_k}P_\fq^{\fn/\fq}\prod_{\fq \mid \fc_{k-1}}P_\fq^{\fc_k} \cdots \prod_{\fq \mid \fc_1}P_\fq^{\fc_2} \nonumber \\
&=&\det\left(
\begin{array}{ccccc}
-P_{\fq_1}^{\fn/\fd} & -P_{\fq_1}^{\fq_2} & \cdots & \cdots& -P_{\fq_1}^{\fq_\nu}\\
-P_{\fq_2}^{\fq_1} & -P_{\fq_2}^{\fn/\fd} & - P_{\fq_2}^{\fq_3} &\cdots & -P_{\fq_2}^{\fq_\nu}\\
\vdots &-P_{\fq_3}^{\fq_2}& \ddots & & \vdots\\
\vdots & \vdots & & \ddots &\vdots\\
-P_{\fq_\nu}^{\fq_1} &-P_{\fq_\nu}^{\fq_2} &\cdots &\cdots & -P_{\fq_\nu}^{\fn/\fd}
\end{array}
\right)=:\mathcal{D}_{\fn,\fd}, \nonumber
\end{eqnarray}

\noindent{}where we write $\fd=\prod_{i=1}^\nu \fq_i$. (When $\fd=1$, set $\mathcal{D}_{\fn,\fd}:=1$.)
Hence we have
$$\left( {\bigwedge}_{\fq\mid \fn}\varphi_\fq^\fq \right)(\epsilon_\fn) = \sum_{\fd \mid \fn} \left( {\bigwedge}_{\fq\mid \fn/\fd}\varphi_\fq^{\fn/\fd} \right)(\epsilon_{\fn/\fd}) \mathcal{D}_{\fn,\fd}.$$
It is easy to see that
\begin{eqnarray}
&&s_\fn \left(\sum_{\fd \mid \fn} \left( {\bigwedge}_{\fq\mid \fn/\fd}\varphi_\fq^{\fn/\fd} \right)(\epsilon_{\fn/\fd}) \mathcal{D}_{\fn,\fd} \right) \nonumber \\
&=& \sum_{\fd \mid \fn} s_{\fn/\fd} \left(\left( {\bigwedge}_{\fq\mid \fn/\fd}\varphi_\fq^{\fn/\fd} \right)(\epsilon_{\fn/\fd})\right) \det\left(
\begin{array}{ccccc}
0& -P_{\fq_1}^{\fq_2} & \cdots & \cdots&  -P_{\fq_1}^{\fq_\nu}\\
-P_{\fq_2}^{\fq_1} & 0& - P_{\fq_2}^{\fq_3} &\cdots &  -P_{\fq_2}^{\fq_\nu}\\
\vdots & -P_{\fq_3}^{\fq_2} & \ddots  & & \vdots\\
\vdots &\vdots& & \ddots &  \vdots\\
-P_{\fq_\nu}^{\fq_1} &-P_{\fq_\nu}^{\fq_2} &\cdots & \cdots&  0
\end{array}
\right) \nonumber \\
&=&\sum_{\tau \in \mathfrak{S}(\fn)}{\rm sgn}(\tau) (-1)^{\nu(\fn/\fd_\tau)}s_{\fd_\tau} \left( \left({\bigwedge}_{\fq \mid \fd_\tau}\varphi_\fq^{\fd_\tau}\right) (\epsilon_{\fd_\tau}) \right) \prod_{\fq \mid \fn/\fd_\tau}P_{\tau(\fq)}^\fq. \nonumber
\end{eqnarray}
Hence we have
\begin{eqnarray}
(-1)^{\nu(\fn)}\left( {\bigwedge}_{\fq\mid \fn}\varphi_\fq^\fq \right)(\epsilon_\fn) &=&  (-1)^{\nu(\fn)} s_\fn\left( \left( {\bigwedge}_{\fq\mid \fn}\varphi_\fq^\fq \right)(\epsilon_\fn) \right) \nonumber \\
&=& \sum_{\tau \in \mathfrak{S}(\fn)}{\rm sgn}(\tau) (-1)^{\nu(\fd_\tau)}s_{\fd_\tau} \left( \left({\bigwedge}_{\fq \mid \fd_\tau}\varphi_\fq^{\fd_\tau}\right) (\epsilon_{\fd_\tau}) \right) \prod_{\fq\mid \fn/\fd_\tau}P_{\tau(\fq)}^\fq. \nonumber
\end{eqnarray}
This is the desired equality (\ref{theta kolyvagin}).
\end{proof}

\begin{proof}[Proof of Theorem \ref{main}]
It is sufficient to show that the image of $c=c(z) \in \mathcal{E}^{\rm b}(T,\mathcal{K})$ under the homomorphism
$$\mathcal{D}_r: \mathcal{E}^{\rm b}(T,\mathcal{K}) \to {\rm KC}_r (\cA)$$
is equal to the image of $\epsilon(z) \in {\rm SS}_r(\cA)$ under
$${\rm Reg}_r:{\rm SS}_r(\cA) \to {\rm KS}_r(\cA).$$
By the definition of $\mathcal{D}_r$, we have
$$\mathcal{D}_r(c)=\Psi_r \left( \left\{ \kappa'(c)_\fn\otimes \bigotimes_{\fq \mid \fn}\sigma_\fq  \right\}_\fn\right).$$
By Lemma \ref{derivative lemma}(ii), we have
$$\left\{ \kappa'(c)_\fn\otimes \bigotimes_{\fq \mid \fn}\sigma_\fq \right\}_\fn= \left\{(-1)^{\nu(\fn)}s_\fn\left( \sum_{\sigma \in \cH_\fn} \sigma c_\fn \otimes \sigma^{-1} \right)\right\}_\fn.$$
By Theorem \ref{MRS}, the right hand side is equal to the image of $\epsilon(z)$ under (\ref{US to DKS}). Hence, by Theorem \ref{commutative theorem}, we see that $\mathcal{D}_r(c)$ coincides with the image of $\epsilon(z)$ under ${\rm Reg}_r$.

\end{proof}

\appendix\section{The basic theory of exterior power biduals}
\subsection{General properties} Let $R$ be a commutative ring and $X$ an $R$-module.

In the sequel, for each non-negative integer $r$ we write $\xi_X^r$ for the canonical homomorphism ${\bigwedge}_R^r X \to {\bigcap}_R^r X$ defined in (\ref{can homo}).

\begin{lemma} \label{rubin prop 1}
If $R$ is noetherian and $X$ is a finitely generated projective $R$-module, then $\xi_X^r$ is an isomorphism.
\end{lemma}

\begin{proof}
Since $R$ is noetherian, the functor $\Hom_R(Y,\cdot)$ is compatible with localization for any finitely generated $R$-module $Y$. Hence, by localization, we may assume that $X$ is free. In this case it is easy to see that the natural map
\begin{eqnarray}
{\bigwedge}_R^r (X^\ast) \to  \left({\bigwedge}_R^r X \right)^\ast; \ \Phi \mapsto (a \mapsto \Phi(a)) \nonumber
\end{eqnarray}
is an isomorphism. Since ${\bigwedge}_R^r X$ is also free, we see that the canonical map ${\bigwedge}_R^r X \to ({\bigwedge}_R^r X)^{\ast \ast}$ is an isomorphism. Hence $\xi_X^r$ is an isomorphism.\end{proof}

For non-negative integers $r$, $s$ with $r\leq s$ and $\Phi \in {\bigwedge}_R^r (X^\ast)$, define a homomorphism
\begin{eqnarray} \label{map bidual}
{\bigcap}_R^s X \to {\bigcap}_R^{s-r} X
\end{eqnarray}
as the dual of
$${\bigwedge}_R^{s-r} (X^\ast) \to {\bigwedge}_R^{s} (X^\ast) ; \ \Psi \mapsto \Phi\wedge \Psi.$$
We denote the map (\ref{map bidual}) simply by $\Phi$, by abuse of notation. It is easy to see that the following diagram commutes
$$\xymatrix{
{\bigwedge}_R^s X \ar[r]^{\Phi} \ar[d]_{\xi_X^s} & {\bigwedge}_R^{s-r}X \ar[d]^{\xi_X^{s-r}} \\
{\bigcap}_R^s X \ar[r]_{\Phi} & {\bigcap}_R^{s-r} X. \\
}
$$

\begin{proposition}\label{prop injective}\
\begin{itemize}
\item[(i)] Each injective homomorphism of $R$-modules $\iota: X \hookrightarrow Y$ for which $\Ext_R^1(\coker \iota, R)$ vanishes induces for each non-negative integer $r$ an injective homomorphism
    $${\bigcap}_R^r X \hookrightarrow {\bigcap}_R^r Y.$$
\item[(ii)] Suppose that we have an exact sequence of $R$-modules
$$ Y \stackrel{\bigoplus_{i=1}^s \varphi_i}{\to} R^{\oplus s} \to Z \to 0.$$
If $Y$ is free of rank $r+s$, then we have
$$\left\langle \im F \ \middle|  \ F \in \im\left( {\bigcap}_R^{r+s} Y \stackrel{{\bigwedge}_{1\leq i \leq s} \varphi_i}{\to } {\bigcap}_R^{r} Y \right)   \right\rangle_R = {\rm Fitt}_R^0(Z),$$
where ${\rm Fitt}_R^0$ denotes the initial Fitting ideal.
\end{itemize}
\end{proposition}

\begin{proof}
(i) Since $\Ext_R^1(\coker \iota,R)=0$, the $R$-linear dual $Y^\ast \to X^\ast$ of $\iota$ is surjective and hence so is the map  ${\bigwedge}_R^r (Y^\ast) \to {\bigwedge}_R^r (X^\ast)$ induced on the exterior powers. Taking the $R$-linear dual of the latter surjection we obtain the claimed injection ${\bigcap}_R^r X \hookrightarrow {\bigcap}_R^r Y.$

(ii) Let $\{b_1,\ldots,b_{r+s}\}$ be a basis of $Y$. Then $b^\ast:=\xi_Y^{r+s}({\bigwedge}_{1\leq i \leq r+s} b_i)$ is a basis of ${\bigcap}_R^{r+s}Y$, by Lemma \ref{rubin prop 1}. By the explicit formula (\ref{exterior explicit}), we have
$$\left({\bigwedge}_{1\leq i \leq s}\varphi_i\right) (b^\ast)=\sum_{\sigma\in \mathfrak{S}_{r+s,s}}{\rm sgn}(\sigma)\det(\varphi_i(b_{\sigma(j)}))_{1\leq i,j\leq s} \cdot \xi_Y^r(b_{\sigma(1+s)}\wedge \cdots \wedge b_{\sigma(r+s)}).$$
Hence we have
$$\im \left(\left({\bigwedge}_{1\leq i \leq s}\varphi_i\right) (b^\ast) \right)=\langle\det(\varphi_i(b_{\sigma(j)}))_{1\leq i,j\leq s} \mid \sigma \in \mathfrak{S}_{r+s,s} \rangle_R={\rm Fitt}_R^0(Z),$$
as claimed.
\end{proof}

\begin{proposition} \label{prop self injective}
Suppose that $R$ is self-injective, i.e. $R$ is injective as an $R$-module, and that we have an exact sequence of $R$-modules
$$0 \to X \to Y \stackrel{\bigoplus_{i=1}^s \varphi_i}{\to} R^{\oplus s},$$
where $s$ is a positive integer. Then, for every non-negative integer $r$, we have
$$\im\left({\bigwedge}_{1\leq i \leq s} \varphi_i: {\bigcap}_R^{r+s} Y \to  {\bigcap}_R^{r} Y \right)\subset {\bigcap}_R^{r} X.$$
Here we regard ${\bigcap}_R^r X \subset {\bigcap}_R^r Y$ by Proposition \ref{prop injective}(i).
In particular, ${\bigwedge}_{1\leq i\leq s} \varphi_i$ induces a homomorphism
$${\bigwedge}_{1\leq i \leq s} \varphi_i: {\bigcap}_R^{r+s}Y \to {\bigcap}_R^r X.$$
\end{proposition}

\begin{proof}
Since the map
$${\bigwedge}_{1\leq i \leq s}\varphi_i : {\bigcap}_R^{r+s} Y \to {\bigcap}_R^r Y$$
is equal to $\varphi_s \circ \cdots \circ \varphi_1$, we may assume $s=1$. (We can inductively prove the statement for general $s$.)

Suppose that $\varphi \in Y^\ast$ and $X=\ker \varphi$. We shall show that the image of
$$\varphi: {\bigcap}_R^{r+1} Y \to {\bigcap}_R^r Y$$
is contained in ${\bigcap}_R^r X$. It is sufficient to show that the kernel of the surjection
\begin{eqnarray} \label{map YX}
{\bigwedge}_R^r (Y^\ast) \to {\bigwedge}_R^r (X^\ast)
\end{eqnarray}
is contained in the kernel of the map
\begin{eqnarray} \label{map YY}
{\bigwedge}_R^r (Y^\ast) \to {\bigwedge}_R^{r+1} (Y^\ast); \ \Phi \mapsto \varphi\wedge \Phi.
\end{eqnarray}
We define a map
\begin{eqnarray} \label{map XY}
{\bigwedge}_R^r (X^\ast) \to {\bigwedge}_R^{r+1}(Y^\ast)
\end{eqnarray}
by $\psi_1\wedge \cdots \wedge \psi_r \mapsto \varphi \wedge \widetilde \psi_1\wedge \cdots \wedge \widetilde \psi_r,$ where $\widetilde \psi_i$ is any lift of $\psi_i$ through $Y^\ast\to X^\ast$. Since $\ker(Y^\ast \to X^\ast) =\langle \varphi \rangle_R$, we see that this map is well-defined. By definition, the composition
$${\bigwedge}_R^r (Y^\ast) \stackrel{(\ref{map YX})}{\to} {\bigwedge}_R^r (X^\ast) \stackrel{(\ref{map XY})}{\to} {\bigwedge}_R^{r+1} (Y^\ast)$$
coincides with (\ref{map YY}). Hence, the kernel of (\ref{map YX}) is contained in that of  (\ref{map YY}).
\end{proof}

\begin{proposition} \label{bidual invariant}
Suppose that $R$ is self-injective. Let $G$ be a finite abelian group and $H$ be a subgroup of $G$.  For an $R[G]$-module $X$ and a non-negative integer $r$, we have a natural identification
$$\left( {\bigcap}_{R[G]}^r X \right)^H = {\bigcap}_{R[G/H]}^r (X^H).$$
\end{proposition}

\begin{proof}
We remark that, since there is a natural isomorphism between $\Hom_{R[G]}(-,R[G])$ and $\Hom_R(- R)$, the notation $X^\ast$ does not make any confusion.
For any $R$-module $X$, we have
$$(X^\ast)^H=\Hom_{\ZZ[H]}(\ZZ, X^\ast)=\Hom_{R[G]}(X\otimes_{\ZZ[H]}\ZZ, R[G] )=(X_H)^\ast.$$
Since $R$ is self-injective, we have $X^{\ast \ast}=X$. So we have
$$(X^H)^\ast=((X^{\ast \ast})^H)^\ast=((X^\ast)_H)^{\ast \ast}= (X^\ast)_H.$$
Hence, we have
\begin{eqnarray*}
\left( {\bigcap}_{R[G]}^r X \right)^H &=&\left( \left({\bigwedge}_{R[G]}^r(X^\ast)\right)^\ast\right)^H \\
&=& \left(\left( {\bigwedge}_{R[G]}^r (X^\ast) \right)_H \right)^\ast \\
&=& \left( {\bigwedge}_{R[G]}^r (X^\ast)_H \right)^\ast \\
&=&  \left( {\bigwedge}_{R[G]}^r (X^H)^\ast \right)^\ast \\
&=& {\bigcap}_{R[G/H]}^r (X^H).
\end{eqnarray*}

\end{proof}

\subsection{Standard representatives}

\begin{definition} \label{definition standard}
Let $C \in D^{\rm p}(R)$ be a perfect complex that is acyclic outside degrees zero and one. Suppose that we have a surjection
$$f: H^1(C) \to X,$$
where $X$ is an $R$-module generated by $r$ elements. A representative of $C$ of the form
$$P' \stackrel{\psi}{\to} P,$$
where $P'$ is placed in degree zero, is called a {\it standard representative} of $C$ with respect to $f$ if the following conditions are satisfied:
\begin{itemize}
\item[(i)] $P'$ is a finitely generated projective module, and $P$ is a free module of finite rank;
\item[(ii)] $P$ and $P'$ are locally isomorphic, i.e. $P_\mathfrak{p} \simeq P'_\mathfrak{p}$ as $R_\mathfrak{p}$-modules for every prime ideal $\mathfrak{p}$ of $R$;
\item[(iii)] there exists a basis $\{b_1,\ldots,b_d\}$ of $P$ (with $d \geq r$) such that the sequence
$$\langle b_{r+1},\ldots, b_d \rangle_R \to H^1(C) \stackrel{f}{\to} X \to 0$$
is exact, where the first map is induced by the natural map
\begin{eqnarray} \label{surjection f}
P \to \coker \psi =H^1(C). \nonumber
\end{eqnarray}
\end{itemize}
When $P=P'$, we call the standard representative {\it quadratic}.

A standard representative of $C$ with respect to $f$ is expressed as $\psi$, $(P' \stackrel{\psi}{\to} P)$,
or $(P',P,\psi, \{b_1,\ldots,b_d\})$ when we specify the basis in the condition (iii). When the standard representative is quadratic, then we abbreviate it as $(P,\psi, \{b_1,\ldots,b_d\})$.

Suppose $X$ is free of rank $r$. In this case, for a standard representative $(P',P,\psi, \{b_1,\ldots ,b_d\})$ of $C$ with respect to $f$, the map
\begin{eqnarray}\label{X isom}
\langle b_1,\ldots,b_r \rangle_R \to H^1(C) \stackrel{f}{\to} X
\end{eqnarray}
must be an isomorphism.
Define a map
\begin{eqnarray}
\Pi_\psi: {\det}_R(C)={\bigwedge}_R^d P' \otimes {\bigwedge}_R^d (P^\ast) \to {\bigwedge}_R^r P' \otimes {\bigwedge}_R^r (X^\ast) \nonumber
\end{eqnarray}
by
$$a \otimes b_1^\ast \wedge \cdots \wedge b_d^\ast \mapsto (-1)^{r(d-r)}\left({\bigwedge}_{r<i \leq d}\psi_i \right)(a) \otimes x_1^\ast \wedge \cdots \wedge x_r^\ast,$$
where $b_i^\ast \in P^\ast$ is the dual of $b_i$, $\psi_i:=b_i^\ast\circ\psi \in (P')^\ast$, and $x_i \in X$ is the image of $b_i$ under the map
$$P \to H^1(C) \stackrel{f}{\to} X.$$
One checks that the map $\Pi_\psi$ does not depend on the choice of $\{b_1,\ldots,b_d\}$.
\end{definition}

\begin{lemma} \label{lemma standard}
Let $C$ and $f: H^1(C)\to X$ be as in Definition \ref{definition standard}, and assume that $X$ is free of rank $r$. Suppose that there exists a standard representative $(P' \stackrel{\psi}{\to} P)$ of $C$ with respect to $f$.
\begin{itemize}
\item[(i)] If $H^1(C)$ is a free $R$-module of rank $r$, $X=H^1(C)$, and $f={\rm id}$, then $H^0(C)=\ker \psi \subset P'$ is a direct summand, locally isomorphic to $H^1(C)$ (in particular, $H^0(C)$ is projective). Furthermore, we have
$$\im \Pi_\psi = {\bigwedge}_R^r H^0(C) \otimes {\bigwedge}_R^r (H^1(C)^\ast),$$
and $\Pi_\psi$ coincides with the canonical isomorphism
\begin{eqnarray}\label{det C isom}
{\det}_R(C) \stackrel{\sim}{\to} {\det}_R(H^0(C)) \otimes {\det}_R^{-1}(H^1(C)).
\end{eqnarray}
\item[(ii)] We have
$${\rm Ann}_R(\im \Pi_\psi)={\rm Ann}_R({\rm Fitt}_R^r(H^1(C)))={\rm Ann}_R({\rm Fitt}_R^0(\ker f)),$$
where ${\rm Fitt}_R^i$ denotes the $i$-th Fitting ideal.
In particular, if $\ker f$ has a non-zero free direct summand, then $\Pi_\psi=0$.
\item[(iii)] For $\Phi \in {\bigwedge}_R^r (P')^\ast$, let
$${\rm ev}_\Phi: {\det}_R(C) \to R$$
be the composition of $\Pi_\psi$ with the map
$${\bigwedge}_R^r P' \otimes {\bigwedge}_R^r (X^\ast) \to R; \ b\otimes x_1^\ast\wedge\cdots\wedge x_r^\ast \mapsto \Phi(b).$$
Then we have
$${\rm Fitt}_R^r(H^1(C))={\rm Fitt}_R^0(\ker f)=\left\langle \im {\rm ev}_\Phi  \ \middle| \ \Phi \in {\bigwedge}_R^r (P')^\ast\right\rangle_R.$$
\item[(iv)] Assume that $R$ is self-injective. Then we have
$$\im \Pi_\psi \subset {\bigcap}_R^r H^0(C) \otimes_R {\bigwedge}_R^r (X^\ast).$$
In particular, $\Pi_\psi$ induces a map
\begin{eqnarray}
\Pi_\psi: {\det}_R(C) \to {\bigcap}_R^r H^0(C) \otimes_R {\bigwedge}_R^r (X^\ast). \nonumber
\end{eqnarray}

\end{itemize}
\end{lemma}

\begin{proof}
(i) Since $H^1(C)$ is free, the exact sequence
$$0 \to \im \psi \to P \to H^1(C)\to 0$$
splits, and hence $\im \psi$ is projective. This implies that the exact sequence
$$0 \to H^0(C) \to P' \to \im (\psi) \to 0$$
also splits, and so $H^0(C)$ is a direct summand of $P'$, and hence projective. Since $P$ and $P'$ are locally isomorphic, $H^0(C)$ and $H^1(C)$ must be locally isomorphic too. Hence, ${\det}_R(H^0(C))$ is identified with ${\bigwedge}_R^r H^0(C)$. It is straightforward to check that $\Pi_\psi$ coincides with the map (\ref{det C isom}).

(ii) By localization, we may assume $P' = P$. By the formula (\ref{exterior explicit}), we have
\begin{eqnarray}
&&\Pi_\psi(b_1\wedge\cdots \wedge b_d \otimes b_1^\ast \wedge \cdots \wedge b_d^\ast) \nonumber\\
&=&(-1)^{r(d-r)} \sum_{\sigma \in \mathfrak{S}_{d,r}}{\rm sgn}(\sigma) \det(\psi_i(b_{\sigma(j)}))_{r < i,j \leq d} b_{\sigma(1)}\wedge\cdots\wedge b_{\sigma(r)} \otimes x_1^\ast \wedge\cdots \wedge x_r^\ast.\nonumber
\end{eqnarray}
Since the map (\ref{X isom}) is an isomorphism, we see that $\psi_i=0$ for $1 \leq i \leq r$, and so ${\rm Fitt}_R^r(H^1(C))$ is generated by $\{\det(\psi_i(b_{\sigma(j)}))_{r < i,j \leq d} \mid \sigma \in \mathfrak{S}_{d,r}  \}$. Hence we have
$${\rm Ann}_R(\im \Pi_\psi)={\rm Ann}_R({\rm Fitt}_R^r(H^1(C)) ).$$
Since $H^1(C) \simeq \ker f \oplus X$, we see that
$${\rm Fitt}_R^r(H^1(C))={\rm Fitt}_R^0(\ker f).$$
If $\ker f$ has a non-zero free direct summand, then ${\rm Fitt}_R^0(\ker f)=0$, and so ${\rm Ann}_R(\im \Pi_\psi)=R$. This implies $\Pi_\psi=0$.

(iii) follows from the proof of (ii).

(iv) Suppose that $R$ is self-injective. Since we have an exact sequence
$$0 \to H^0(C) \to P' \stackrel{\bigoplus_{r<i \leq d}\psi_i}{\to} R^{\oplus d-r},$$
we see by Proposition \ref{prop self injective} that
$\im \Pi_\psi \subset {\bigcap}_R^r H^0(C) \otimes_R {\bigwedge}_R^r (X^\ast).$
\end{proof}

\subsection{The case of orders}  \label{A order}
Let $O$ be a Dedekind ring (which is not a field) with quotient field $Q$. Let $\cQ$ be a finite dimensional semisimple commutative $Q$-algebra.

\begin{proposition} \label{rubin prop 2}
Let $\mathcal{R}$ be an $O$-order in $\cQ$, and $X$ a finitely generated $\mathcal{R}$-module. Then, for any non-negative integer $r$, $\xi_X^r$ induces an isomorphism
$$\left\{ a \in Q \otimes_O {\bigwedge}_{\mathcal{R}}^r  X \ \middle| \ \Phi(a) \in \mathcal{R} \text{ for all $\Phi \in {\bigwedge}_{\mathcal{R}}^r (X^\ast)$} \right\} \stackrel{\sim}{\to} {\bigcap}_{\mathcal{R}}^r X.$$
\end{proposition}

\begin{proof}
By definition, we have
$${\bigcap}_{\mathcal{R}}^r X=\ker\left(  \Hom_{\mathcal{R}}\left({\bigwedge}_{\mathcal{R}}^r (X^\ast), \cQ \right) \to \Hom_{\mathcal{R}}\left({\bigwedge}_{\mathcal{R}}^r (X^\ast), \cQ/\mathcal{R} \right) \right).$$
Since $\cQ$ is semisimple, we see that $Q\otimes_O X$ is a finitely generated projective $\cQ$-module. So, by Lemma \ref{rubin prop 1}, we have
$$\xi_{Q\otimes_O X}^r: Q \otimes_O {\bigwedge}_{\mathcal{R}}^r X \stackrel{\sim}{\to}  {\bigcap}_{\cQ}^r (Q\otimes_O X)= \Hom_{\mathcal{R}}\left({\bigwedge}_{\mathcal{R}}^r (X^\ast), \cQ \right).  $$
The proposition follows by noting that
\begin{eqnarray}
&&\left\{ a \in Q \otimes_O {\bigwedge}_{\mathcal{R}}^r  X \ \middle| \ \Phi(a) \in \mathcal{R} \text{ for all $\Phi \in {\bigwedge}_{\mathcal{R}}^r (X^\ast)$}\right\} \nonumber \\
&=& \ker \left( Q \otimes_O {\bigwedge}_{\mathcal{R}}^r X  \to \Hom_{\mathcal{R}}\left({\bigwedge}_{\mathcal{R}}^r (X^\ast), \cQ/\mathcal{R} \right) \right). \nonumber
\end{eqnarray}
\end{proof}

By the above proposition, we often regard ${\bigcap}_{\mathcal{R}}^r X$ as a sublattice of $Q \otimes_O {\bigwedge}_{\mathcal{R}}^r X$.

In the following, we assume that $\cQ$ is separable over $Q$. For an $O$-order $\mathcal{R}$ in $\cQ$, the following are equivalent:
\begin{itemize}
\item[(a)] $\mathcal{R}$ is a (one-dimensional) Gorenstein ring;
\item[(b)] the injective dimension of $\mathcal{R}$ as an $\mathcal{R}$-module is one;
\item[(c)] $\Hom_O(\mathcal{R}, O)$ is projective as an $\mathcal{R}$-module;
\item[(d)] every short exact sequence of the form
$$0 \to \mathcal{R} \to X \to Y \to 0,$$
where $X$ and $Y$ are finitely generated $\mathcal{R}$-modules which are $O$-torsion-free, splits.
\end{itemize}
(See \cite[\S 37]{CR} and \cite[Proposition 6.1]{DKR}.)
$\mathcal{R}$ is said to be a Gorenstein order if one of the equivalent conditions above is satisfied.

Group rings are typical examples of Gorenstein orders. Letting $G$ be a finite abelian group whose order is prime to the characteristic of $Q$, the group ring $O[G]$ is a Gorenstein order in $Q[G]$ (see \cite[Corollary 10.29]{CR}).


We note some properties of modules over Gorenstein orders. Suppose $\cR$ is a Gorenstein order. For a finitely generated $\cR$-module $X$ which is $O$-torsion-free, we have
\begin{eqnarray}\label{extvanish}
\Ext_\mathcal{R}^1(X, \mathcal{R})=0 .
\end{eqnarray}
In fact, condition (d) above implies that any 1-extension of $X$ by $\cR$ splits. Also, note that $X$ is reflexive, i.e. the natural map
\begin{eqnarray} \label{2dual}
X \to X^{**}
\end{eqnarray}
is an isomorphism. In fact, by \cite[Theorem 6.2]{bassgorenstein}, every finitely generated torsionless module over a one-dimensional Gorenstein ring is reflexive. Here $X$ is said to be torsionless if the map (\ref{2dual}) is injective. (In our case, note that any finitely generated $\cR$-module which is $O$-torsion-free is torsionless as an $\cR$-module).

For an idempotent $e \in \cQ$ and an $\mathcal{R}$-module $X$, we set
$$X[e]:=\{ a \in X \mid e\cdot (1\otimes a)=0 \text{ in }Q\otimes_O X\}.$$

In the following, we assume $\cR$ is a Gorenstein order.

\begin{lemma} \label{projective lemma}
Let $X$ be a finitely generated $\mathcal{R}$-module that is $O$-torsion-free and $P$ a finitely generated projective $\mathcal{R}$-module. Suppose that there is a surjection $\pi: P^\ast \to X^\ast$, and regard $X \subset P$ via the dual of $\pi$. Then, for any non-negative integer $r$, and any idempotent $e$ of $\cQ$ one has
$$\left({\bigcap}_{\mathcal{R}}^r X\right)[e] = (1-e)\!\!\left( Q\otimes_O {\bigwedge}_{\mathcal{R}}^r X \right) \cap \left({\bigwedge}_{\mathcal{R}}^r P \right).$$
Here we regard
$${\bigcap}_{\mathcal{R}}^r X \subset {\bigcap}_{\mathcal{R}}^r P={\bigwedge}_{\mathcal{R}}^r P$$
by Proposition \ref{prop injective}(i) and Lemma \ref{rubin prop 1}.
\end{lemma}

\begin{proof} If one omits the idempotent $e$, then the proof of this result is the same as that of  \cite[Lemma 4.7(ii)]{bks1} (since the property (\ref{extvanish}) holds). The general case then follows directly upon noting that an element of ${\bigcap}_{\mathcal{R}}^r X$ belongs to $\left({\bigcap}_{\mathcal{R}}^r X\right)[e]$ if and only if it is stable under multiplication by $(1-e)$.
%
%
\end{proof}

\subsection{The proof of Proposition \ref{prop admissible new}} \label{A proof}
In this section we use the notation and hypotheses of Proposition \ref{prop admissible new}.

Let $\{e_1,\ldots,e_s \}$ be the complete set of the primitive orthogonal idempotents of $\cQ$ (so each $\cQ e_i$ is an extension field over $Q$ and we have the decomposition $\cQ=\bigoplus_{i=1}^s \cQ e_i$).

For an $\mathcal{R}$-submodule $X$ of $\cQ$ we also use the $\mathcal{R}$-module obtained by setting
\[ X^{-1} := \{c \in \cQ \mid c\cdot X \subset \mathcal{R} \text{ and for every $1\leq i \leq s$ one has }\,\, e_i( Q \otimes_O X) = 0 \Longrightarrow e_i \cdot c = 0\}.\]

\begin{proposition} \label{prop admissible}
Assume that $O$ is local, i.e. a discrete valuation ring. Let $\pi \in O$ be a uniformizer. Assume $(\cR/\pi\cR)/{\rm rad}(\cR/\pi \cR)$ is separable over $O/\pi O$, where ${\rm rad}$ denotes the Jacobson radical.
Let $C \in D^{\rm p}(\mathcal{R})$ be an admissible complex (see Definition \ref{definition admissible}). Let $X$ be an $\mathcal{R}$-module generated by $r$ elements, and suppose that there is a surjective homomorphism $f: H^1(C) \to X.$ Let $e_r \in \cQ$ be the sum of the primitive idempotents which annihilate $Q \otimes_O \ker f$.
\begin{itemize}
\item[(i)] There exists a quadratic standard representative of $C$ with respect to $f$. (We let $R$ in  Definition \ref{definition standard} be $\mathcal{R}$.)
\item[(ii)] If $X$ is free of rank $r$, then the map $\Pi_\psi$ in Definition \ref{definition standard} coincides with the composite
\begin{eqnarray}
 {\det}_{\mathcal{R}}(C)&\hookrightarrow&Q \otimes_O {\det}_{\mathcal{R}}(C) \nonumber\\
&=& {\det}_{\cQ}(Q \otimes_O C) \nonumber \\
&\simeq& {\det}_{\cQ} (Q\otimes_O H^0(C)) \otimes {\det}_{\cQ}^{-1} (Q \otimes_O H^1(C))   \nonumber  \\
&\stackrel{e_r \times}{\to}& e_r \left( {\det}_{\cQ} (Q\otimes_O H^0(C)) \otimes {\det}_{\cQ}^{-1} (Q \otimes_O H^1(C)) \right)   \nonumber \\
&\simeq& e_r\left( {\bigwedge}_{\cQ}^r (Q\otimes_O H^0(C))\right) \otimes {\bigwedge}_{\cQ}^r (Q \otimes_O X^\ast).  \nonumber
\end{eqnarray}
In particular, by Lemma \ref{projective lemma}, we see that $\Pi_\psi$ induces a map
\begin{eqnarray}
\Pi_{C,f}: {\det}_{\mathcal{R}}(C) &\to&  \left({\bigcap}_{\mathcal{R}}^r H^0(C)\right)[1-e_r] \otimes {\bigwedge}_{\mathcal{R}}^r (X^\ast) , \nonumber
\end{eqnarray}
and this map is independent of the choice of the standard representative.

\item[(iii)] Assume $X$ is free of rank $r$. Then $\mathcal{R} e_r\subset {\rm Fitt}^0_{\mathcal{R}}(\ker(f))^{-1} \subset \cQ e_r$ and for the map $\Pi_{C,f}$ constructed in (ii) there is an isomorphism of $\mathcal{R}$-modules
\[{\rm coker}(\Pi_{C,f}) \simeq \frac{{\rm Fitt}^0_{\mathcal{R}}(\ker(f))^{-1}}{\mathcal{R} e_r}\]

\item[(iv)] Let $I$ be an ideal of $O$. Then there exists a quadratic standard representative of $C \otimes_O O/I \in D^{\rm p}(\cR/I \cR)$ with respect to $ f_I :  H^1(C \otimes_O O/I)=H^1(C) \otimes_O O/I \stackrel{f}{\to} X \otimes_O O/I$.
\item[(v)] Let $O$ and $I$ be as in (iv). Assume $I$ is non-zero and $X$ is free of rank $r$. Then the map $\Pi_\psi$ in Lemma \ref{lemma standard}(iv) is independent of the choice of the standard representative and so can be denoted $\Pi_{C \otimes_O O/I, f_I}$.
\end{itemize}
\end{proposition}

\begin{proof}
(i) Let $P_1$ be a free $\mathcal{R}$-module of rank $r$, and take a surjection $g: P_1 \to X$. Let $\pi_1: P_1 \to H^1(C)$ be a lift of $g$, namely, a map satisfying $f \circ \pi_1=g$. Let $P_2$ be another free $\mathcal{R}$-module, whose rank is sufficiently large so that we can take a surjection $\pi_2 : P_2 \to \ker f$. Put $P:=P_1 \oplus P_2$, and $\pi:=\pi_1\oplus \pi_2: P \to H^1(C)$. By construction, we see that $\pi$ is surjective. Let $\{b_1,\ldots,b_r\}$ be a basis of $P_1$, and $\{b_{r+1},\ldots , b_d\}$ a basis of $P_2$. Then $\{b_1,\ldots,b_d\}$ is a basis of $P$, and we have an exact sequence
$$\langle b_{r+1},\ldots, b_d\rangle_R \stackrel{\pi}{\to} H^1(C) \stackrel{f}{\to} X \to 0.$$

Since $C$ is acyclic outside degrees zero and one, it corresponds to a unique Yoneda extension class in $\Ext_{\mathcal{R}}^2(H^1(C), H^0(C))$, which we denote by $\tau_C$. By considering a projective resolution of $H^1(C)$ of the form
$$\cdots \to P \stackrel{\pi}{\to} H^1(C) \to 0,$$
we can take a representative of $\tau_C$ of the form
$$0 \to H^0(C) \to P' \stackrel{\psi}{\to} P \stackrel{\pi}{\to} H^1(C) \to 0 ,$$
where $P'$ is a finitely generated $\mathcal{R}$-module. Since $C$ is perfect, we see that $P'$ has finite projective dimension. Since $H^0(C)$ is $O$-torsion-free, the above exact sequence implies that $P'$ is also $O$-torsion-free. Hence, by Lemma \ref{fproj} below, we see that $P'$ is a projective $\mathcal{R}$-module.

Since the Euler characteristic of $Q\otimes_O C$ is zero, we must have $Q \otimes_O P' \simeq Q \otimes_O P$.
By Hattori's theorem \cite[Theorem 32.5]{CR}, we have $P' \simeq P$. (Here the separability of $(\cR/\pi\cR)/{\rm rad}(\cR/\pi\cR)$ over $O/\pi$ is needed.) Thus, we may assume $P'=P$ and we have proved $(P \stackrel{\psi}{\to} P)$ is a quadratic standard representative of $C$ with respect to $f$.

(ii) follows from Lemma \ref{lemma standard}(i) and (ii).

(iii) The inclusions $\mathcal{R} e_r\subset  {\rm Fitt}^0_{\mathcal{R}}(\ker(f))^{-1} \subset \cQ e_r$ follows directly from the definition of ${\rm Fitt}^0_{\mathcal{R}}(\ker(f))^{-1}$ and the fact that for each primitive idempotent $e_i$ of $\cQ$ one has
\[
e_i \cdot e_r \neq 0 \Longleftrightarrow
 e_i(Q \otimes_O\ker(f))=0 \Longleftrightarrow e_i(Q\otimes_O{\rm Fitt}^0_{\mathcal{R}}(\ker(f))) \not=0.\]

Next we note that since, by assumption, $X$ is free and $P' = P$, the argument of Lemma \ref{lemma standard}(ii) shows $\im(\Pi_{\psi})$ is generated as an $\mathcal{R}$-module by the element $y \otimes x_1^\ast \wedge\cdots \wedge x_r^\ast$ with

\[ y := \sum_{\sigma \in \mathfrak{S}_{d,r}}c_\sigma \cdot(b_{\sigma(1)}\wedge\cdots\wedge b_{\sigma(r)}) \in \left({\bigcap}_{\mathcal{R}}^r H^0(C)\right)[1-e_r] \]
for a specific set of generators $\{c_\sigma:= {\rm sgn}(\sigma) \det(\psi_i(b_{\sigma(j)}))_{r < i,j \leq d}\}_{\sigma \in \mathfrak{S}_{d,r}}$ of the $\mathcal{R}$-module ${\rm Fitt}_{\mathcal{R}}^0(\ker(f))$.

From \cite[Lemma 4.2]{bks1} we also know that
\[ \cQ\cdot y = \left({\bigwedge}_{\cQ}^r (Q\otimes_O H^0(C)) \right)[1-e_r] = Q\otimes_O\left({\bigcap}_{\mathcal{R}}^r H^0(C)\right)[1-e_r]\]
and hence, by Lemma \ref{projective lemma} (with $X = H^1(C)$ and $e = 1-e_r$) that
\[ \left({\bigcap}_{\mathcal{R}}^r H^0(C)\right)[1-e_r] = \left\{c\cdot y \ \middle| \ c \in \cQ e_r \,\,\,\text{ and }\,\,\, c\cdot y \in {\bigwedge}_{\mathcal{R}}^r P \right\}.\]

In addition, since $\{b_{\sigma(1)}\wedge\cdots\wedge b_{\sigma(r)}\}_{\sigma \in \mathfrak{S}_{d,r}}$ is an $\mathcal{R}$-basis of ${\bigwedge}_{\mathcal{R}}^r P$, one has
\begin{multline*} c\cdot y \in {\bigwedge}_{\mathcal{R}}^r P \Longleftrightarrow c\cdot c_\sigma\in \mathcal{R} \,\,\text{ for all } \,\, \sigma \in \mathfrak{S}_{d,r}\\
\Longleftrightarrow c\cdot {\rm Fitt}_{\mathcal{R}}^0(\ker(f))\subset  \mathcal{R} \Longleftrightarrow c\in {\rm Fitt}_{\mathcal{R}}^0(\ker(f))^{-1}.\end{multline*}
The claimed isomorphism in (iii) is therefore induced by the assignment $c\cdot y \mapsto c$.

(iv) Let $(P \stackrel{\psi}{\to} P)$ be as in the proof of (i).
Clearly, $(P/IP \stackrel{\psi}{\to} P/IP)$ is a quadratic standard representative of $C \otimes_O O/I$ with respect to $ f_I$.

(v) For simplicity, we denote $O/I$, $\cR/I\cR$ and $C\otimes_O O/I$ by $O_I$, $\cR_I$ and $C_I$ respectively.

Note that $\cR_I$ is artinian, so a product of local rings. Considering component-wise, we may assume $\cR_I$ is local. Also note that $\cR_I$ is a zero-dimensional Gorenstein ring, i.e. a self-injective ring.

As in the proof of (i), construct a surjection $\pi=\pi_1\oplus \pi_2: P=P_1\oplus P_2 \to H^1(C_I)$, so that the rank of $P_2$ is minimal.
Let $( P \stackrel{\psi}{\to}  P)$ be the quadratic standard representative of $C_I$ with respect to $ f_I$, which is constructed in the proof of (i) (and (iv)). Let
$$0 \to H^0(C_I) \stackrel{\iota}{\to}  P \stackrel{\psi}{\to} P \stackrel{\pi}{\to} H^1(C_I) \to 0 $$
be the tautological exact sequence. We claim that

\begin{itemize}
\item[(a)] $\pi$ is a projective cover of $H^1(C_I)$;
\item[(b)] $\iota$ is an injective envelope of $H^0(C_I)$.
\end{itemize}

To see (a), we remark that $\pi_2: P_2 \to \ker f$ is a projective cover, since the rank of $P_2$ is minimal. (Note that $\mathcal{R}_I$ is local so every projective module is free.) Noting this, one easily sees that $\pi$ is a projective cover.

We show (b). Since $\mathcal{R}_I$ is self-injective and local, the free modules are exactly the injective modules. Let $\mu: H^0(C_I) \hookrightarrow E$ be an injective envelope. Take a representative of $\tau_{C_I} \in\Ext_{\mathcal{R}_I}^2(H^1(C_I), H^0(C_I)) $ of the form
$$0 \to H^0(C_I) \stackrel{\mu}{\to} E \to E' \to H^1(C_I) \to 0.$$
(This is possible by taking a representative of the image of $\tau_{C_I}$ in $\Ext_{\mathcal{R}_I}^1(H^1(C_I), \coker( \mu)) (\simeq \Ext_{\mathcal{R}_I}^2 (H^1(C_I),H^0(C_I)))$.) There is a quasi-isomorphism

$$(P \stackrel{\psi}{\to} P) \to (E \to E'),$$

\noindent{}and the acyclicity of the mapping cone gives an exact sequence

$$0 \to P \to P\oplus E \to E' \to 0.$$
From this, we have $P \oplus E \simeq P\oplus E'$. Since $\mathcal{R}_I$ is local, we can cancel $P$, and we have $E \simeq E'$ (see \cite{evans}). So we may assume $E=E'$. By the property of injective envelope, there is a split injection $\alpha: E \hookrightarrow P$ such that $\alpha\circ \mu =\iota$. If the rank of $E$ is strictly smaller than that of $P$, this contradicts (a). Hence $\alpha$ must be an isomorphism. This shows the claim (b).

By the above claims (a) and (b), we see that any quadratic standard representative of $C_I$ with respect to $ f_I$ is of the form
$$(P\oplus F \stackrel{\psi\oplus \varphi}{\to} P \oplus F),$$
where $F$ is free of finite rank and $\varphi$ is an automorphism of $F$. Denote this complex by $P^\bullet\oplus F^\bullet$, and set $P^\bullet:=(P \stackrel{\psi}{\to} P)$. The mapping cone of $P^\bullet \to P^\bullet\oplus F^\bullet$ is quasi-isomorphic to $F^\bullet:=(F \stackrel{\varphi}{\to} F)$, so we have a canonical isomorphism

\begin{eqnarray}\label{PF isom}
{\det}_{\mathcal{R}_I}(P^\bullet \oplus F^\bullet)\simeq {\det}_{\mathcal{R}_I}(P^\bullet) \otimes {\det}_{\mathcal{R}_I}(F^\bullet) \simeq {\det}_{\mathcal{R}_I}(P^\bullet),
\end{eqnarray}
where the second isomorphism is induced by
$$ {\det}_{\mathcal{R}_I}(F^\bullet) \simeq {\det}_{\mathcal{R}_I}(0) \otimes {\det}_{\mathcal{R}_I}^{-1}(0) =\mathcal{R}_I.$$
When we identify ${\det}_{\mathcal{R}_I}(C_I)= {\det}_{\mathcal{R}_I}(P^\bullet)$, the identification ${\det}_{\mathcal{R}_I}(C_I)= {\det}_{\mathcal{R}_I}(P^\bullet\oplus F^\bullet)$ must be done via
$${\det}_{\mathcal{R}_I}(C_I)= {\det}_{\mathcal{R}_I}(P^\bullet) \stackrel{(\ref{PF isom})}{\simeq} {\det}_{\mathcal{R}_I}(P^\bullet\oplus F^\bullet).$$
It is now clear that $\Pi_\psi=\Pi_{\psi \oplus \varphi}$. Hence, we have proved that $\Pi_\psi$ does not depend on the choice of the standard representative of $C_I$.
\end{proof}

\begin{lemma}\label{fproj}
Let $X$ be a finitely generated $\cR$-module. Then, $X$ is projective if and only if $X$ is $O$-torsion-free and has finite projective dimension.
\end{lemma}

\begin{proof}
The ``only if" part is clear.
Suppose that $X$ is $O$-torsion-free and has finite projective dimension. Choose a projective resolution
$$0 \to P_n \stackrel{f}{\to} P_{n-1} \to \cdots \to P_0 \to X \to 0.$$
Since $\cR$ is a Gorenstein order, we see that the $O$-dual of any finitely generated projective $\cR$-module is again projective. Taking the $O$-dual of the short exact sequence
$$0 \to P_n \stackrel{f}{\to} P_{n-1} \to \coker(f) \to 0,$$
we see that exactness is preserved since each term is $O$-torsion-free and that $\coker(f)$ is projective. Repeating this argument, we deduce from the resolution above that $X$ is projective.

\end{proof}

\end{document}